\documentclass[12pt]{article}
\usepackage{amssymb}

\newtheorem{definition}{Definition}[section]
\newtheorem{theorem}[definition]{Theorem}
\newtheorem{lemma}[definition]{Lemma}
\newtheorem{corollary}[definition]{Corollary}

\newtheorem{note}[definition]{Note}

\def\K{\mathbb K}

\def\fld{\mathbb K}
\def\Mdf{{\hbox{Mat}}_{d+1}(\K)}
\def\alg{\cal A}

\newcommand{\beast}{\begin{eqnarray*}}
\newcommand{\eeast}{\end{eqnarray*}}

\begin{document}
\newenvironment{proof}{\noindent{\it Proof\/}:}{\par\noindent $\Box$\par}

\title{ \bf Leonard pairs from 24\\
points of view}
\author{Paul Terwilliger  
}
\date{}
\maketitle
\begin{abstract} 
Let $\K$ denote a field, and let $V$ denote a  
vector space over $\K$ with finite positive dimension.
We consider a pair
of linear transformations
$A:V\rightarrow V$ and $A^*:V\rightarrow V$
that satisfy both conditions below:
\begin{enumerate}
\item There exists a basis for $V$ with respect to which
the matrix representing $A$ is diagonal and the matrix
representing $A^*$ is irreducible tridiagonal.
\item There exists a basis for $V$ with respect to which
the matrix representing $A^*$ is diagonal and the matrix
representing $A$ is irreducible tridiagonal.
\end{enumerate}

\medskip
\noindent
We call such a pair a {\it Leonard pair} on $V$. Referring to the above
Leonard pair,
we investigate 24 bases for $V$ on which the 
 action of $A$ and $A^*$ takes an attractive form. 
Our bases are described  as follows.
Let  $\Omega$ denote the set consisting of four symbols 
$0,d,0^*, d^* $.
We identify the symmetric group $S_4$ with the
set of all linear orderings of $\Omega$.
For each element $g$ of $S_4$,
we define an (ordered) basis for $V$, which we denote by 
$\lbrack g \rbrack $. The 24 resulting bases 
are related as follows.
For all elements $wxyz$ in  $S_4$,
the transition matrix from the basis
 $\lbrack wxyz \rbrack $ to the basis 
 $\lbrack xwyz \rbrack $
 (resp. $\lbrack wyxz \rbrack $)
is diagonal (resp. 
lower triangular). 
The basis 
 $\lbrack wxzy \rbrack $ is the basis  
 $\lbrack wxyz \rbrack $ in inverted order. 
The transformations $A$ and $A^*$ act on the 24 bases as follows.
For all $g \in S_4$,
let $A^g$ (resp. $A^{*g}$) denote
the matrix representing  $A$ (resp. $A^*$)
with respect to
$\lbrack g \rbrack $.
To describe 
$A^g$ and $A^{*g}$,
we refer to $0^*,  d^*$ as the {\it starred} elements of $\Omega $.
Writing $g=wxyz$,
if neither of $y,z$ are starred then
 $A^g$ is 
 diagonal
and  $A^{*g}$
is irreducible tridiagonal.
If $y$ is starred but $z$ not,
then $A^g$ is lower bidiagonal and $A^{*g}$ is upper bidiagonal.
If $z$ is starred but $y$ not,
then $A^g$ is upper bidiagonal and $A^{*g}$ is lower bidiagonal.
If both of $y,z$ are starred, then  
 $A^g$ is irreducible tridiagonal 
and $A^{*g}$ is diagonal.

We define a symmetric binary relation
on $S_4$ called adjacency.
An element $wxyz$ of $S_4$ 
is by definition   adjacent to each of
$xwyz$, $wyxz$, $wxzy$ and no other elements of $S_4$.
For all ordered pairs of adjacent elements $g,h $ in $S_4$, we
find the entries of the transition matrix from
the basis $\lbrack g \rbrack $ to the basis 
$\lbrack h \rbrack $. We express these entries in terms
of the eigenvalues of $A$, the eigenvalues of $A^*$, and two sequences
of parameters called the
first split sequence and the second split sequence.
%
For all $g \in S_4$, we compute the entries of
$A^g$ and $A^{*g}$
in terms of the eigenvalues of $A$, the eigenvalues of $A^*$,
the first split sequence 
and the second split sequence. 
%
\end{abstract}

%

\section{Leonard pairs}
\medskip
\noindent Throughout this paper, $\K$ will denote an arbitrary
field, and $\tilde \K$ will denote the algebraic closure
of $\K$.

\medskip
\noindent 
We begin by recalling the notion of a Leonard pair.

\begin{definition} \cite{LS99}
\label{def:lprecall}
Let 
 $V$ denote a  
vector space over $\K$ with finite positive dimension.
By a {\it Leonard pair} on $V$,
we mean an ordered pair $A, A^*$, where
$A:V\rightarrow V$ and $A^*:V\rightarrow V$ are linear transformations 
that 
 satisfy both (i), (ii) below. 
\begin{enumerate}
\item There exists a basis for $V$ with respect to which
the matrix representing $A$ is diagonal and the matrix
representing $A^*$ is irreducible tridiagonal.
\item There exists a basis for $V$ with respect to which
the matrix representing $A^*$ is diagonal and the matrix
representing $A$ is irreducible tridiagonal.

\end{enumerate}
(A tridiagonal matrix is said to be irreducible
whenever all entries immediately above and below the main
diagonal are nonzero).

\end{definition}

\begin{note} According to a common notational convention, for
a linear transformation $A$ the conjugate-transpose of $A$ is denoted
$A^*$. We emphasize we are not using this convention. In a Leonard
pair $A,A^*$, the linear transformations $A$ and $A^*$
are arbitrary subject
to (i),  (ii) above.
\end{note}
\noindent
Here is an example of a Leonard pair.
Set 
$V={\K}^4$ (column vectors), set 
\beast
A = 
\left(
\begin{array}{ c c c c }
0 & 3  &  0    & 0  \\
1 & 0  &  2   &  0    \\
0  & 2  & 0   & 1 \\
0  & 0  & 3  & 0 \\
\end{array}
\right), \qquad  
A^* = 
\left(
\begin{array}{ c c c c }
3 & 0  &  0    & 0  \\
0 & 1  &  0   &  0    \\
0  & 0  & -1   & 0 \\
0  & 0  & 0  & -3 \\
\end{array}
\right),
\eeast
and view $A$ and $A^*$  as linear transformations from $V$ to $V$.
We assume 
the characteristic of $\K$ is not 2 or 3, to ensure
$A$ is irreducible.
Then $A, A^*$ is a Leonard
pair on $V$. 
Indeed, 
condition (ii) in Definition
\ref{def:lprecall}
is satisfied by the basis for $V$
consisting of the columns of the 4 by 4 identity matrix.
To verify condition (i), we display an invertible  matrix  
$P$ such that 
$P^{-1}AP$ is 
diagonal and 
$P^{-1}A^*P$ is
irreducible tridiagonal.
Set 
\beast
P = 
\left(
\begin{array}{ c c c c}
1 & 3  &  3    &  1 \\
1 & 1  &  -1    &  -1\\
1  & -1  & -1  & 1  \\
1  & -3  & 3  & -1 \\
\end{array}
\right).
\eeast
 By matrix multiplication $P^2=8I$, where $I$ denotes the identity,   
so $P^{-1}$ exists. Also by matrix multiplication,    
\begin{equation}
AP = PA^*.
\label{eq:apeq}
\end{equation}
Apparently
$P^{-1}AP$ equals $A^*$ and is therefore diagonal.
By (\ref{eq:apeq}), and since $P^{-1}$ is
a scalar multiple of $P$, we find
$P^{-1}A^*P$ equals $A$ and is therefore irreducible tridiagonal.  Now 
condition (i) of  Definition 
\ref{def:lprecall}
is satisfied
by the basis for $V$ consisting of the columns of $P$. 

\medskip
\noindent The above example is a member of the following infinite
family of Leonard pairs.
For any nonnegative integer $d$, 
the pair
\begin{equation}
A = 
\left(
\begin{array}{ c c c c c c}
0 & d  &      &      &   &{\bf 0} \\
1 & 0  &  d-1   &      &   &  \\
  & 2  &  \cdot    & \cdot  &   & \\
  &   & \cdot     & \cdot  & \cdot   & \\
  &   &           &  \cdot & \cdot & 1 \\
{\bf 0} &   &   &   & d & 0  
\end{array}
\right),
\qquad A^*= \hbox{diag}(d, d-2, d-4, \ldots, -d)
\label{eq:fam1}
\end{equation}
is a Leonard pair on the vector space $\K^{d+1} $,
provided the 
 characteristic of $\K$ is zero or an odd prime greater than $d$.
This can be  proved by modifying the 
 proof for $d=3$ given above. One shows  
$P^2=2^dI$  and $AP= PA^*$, where 
$P$ denotes the matrix with $ij$ entry
\begin{equation}
P_{ij} =  
\Biggl({{ d }\atop {j}}\Biggr) {{}_2}F_1\Biggl({{-i, -j}\atop {-d}}
\;\Bigg\vert \;2\Biggr)
\qquad \qquad (0 \leq i,j\leq d).
\label{eq:ex1}
\end{equation}
We follow the standard notation for
hypergeometric series \cite{gasperrahmanbk}. 
 The details of the above calculations
are given in Section 16 below.

\medskip
\noindent 
To motivate our results we mention some background on Leonard pairs. 
There is a connection between
Leonard pairs and certain orthogonal polynomials contained
in the Askey scheme \cite{KoeSwa}.
Observe the ${{}_2F_1}$ that appears in
(\ref{eq:ex1}) 
is a Krawtchouk polynomial \cite{KoeSwa}.
There exist 
families of Leonard pairs similar to the one above in 
which 
 the Krawtchouk polynomial is replaced by one of the following.

\medskip
\centerline{
\begin{tabular}[t]{c|c}
        type & polynomial \\ \hline 
 ${{}_4F_3}$ & Racah \\ 
 ${{}_3F_2}$ & Hahn, dual Hahn \\ 
 ${{}_2F_1}$ & Krawtchouk \\ 
 ${{}_4\phi_3}$ & $q$-Racah \\ 
 ${{}_3\phi_2}$ & $q$-Hahn, dual $q$-Hahn \\ 
 ${{}_2\phi_1}$ & $q$-Krawtchouk (classical, affine, quantum, dual) 
\end{tabular}}

\medskip
\noindent
The above polynomials are defined in Koekoek and Swarttouw \cite{KoeSwa},
and the connection to Leonard pairs  is given in 
\cite[ch. 15]{LS99} and 
\cite[p. 260]{BanIto}. This connection is also discussed in 
 Section 16 below.
%

\medskip
\noindent Leonard pairs play a role in representation theory.
For instance, Leonard pairs arise naturally
in the representation theory 
 of the Lie algebra $sl_2$ 
 \cite{TD00}, the quantum algebra
$U_q(sl_2)$ 
\cite{Koelink3},
\cite{Koelink1},
\cite{Koelink2}, 
\cite{Koelink4},
\cite{koo3},
\cite[ch. 4]{Hjal},
\cite{Terint}, 
\cite{qSerre},
the Askey-Wilson algebra 
\cite{GYZnature},
 \cite{GYLZmut},
\cite{GYZTwisted},
\cite{GYZlinear},
\cite{GYZspherical},
\cite{Zhidd}, 
\cite{ZheCart},
\cite{Zhidden},
and 
the Tridiagonal algebra 
\cite{TD00},
\cite{qSerre},
\cite{LS99}. 

\medskip
\noindent
Leonard pairs play a role  in combinatorics. For instance,
there is a combinatorial object  
called a  $P$-and $Q$-polynomial
association scheme \cite{BanIto}, \cite{bcn}, \cite{Leopandq},
\cite{Tercharpq}, \cite{Ternew}.
Leonard pairs have been used to describe
certain irreducible modules for the subconstitutent algebra of
 these schemes \cite{TersubI}, 
\cite{TersubII}, \cite{TersubIII}.  
See \cite{Cau}, \cite{CurNom}, \cite{Curspin}, \cite{go}, \cite{HobIto},
\cite{TD00},
\cite{Tan}
for more
information on Leonard pairs and association  schemes.

\medskip
\noindent 
Leonard pairs are closely related to the  work of 
Grunbaum and Haine on the 
 ``bispectral problem''
\cite{GH7},
\cite{GH6}.
See
\cite{GH4},
\cite{GH5},
\cite{GH1}, 
\cite{GH3},
\cite{GH2} 
for related work.

\medskip
\noindent We now give an overview of the present paper.
Let $V$ denote a vector space over $\K$ with finite
positive dimension, and let $A, A^*$ denote a Leonard 
pair on $V$. Using this pair, we  define 
24 bases for $V$  which we find attractive. In our
study of these 24 bases,
we will be concerned with
(i) how these bases are related to each other,
and (ii)  for each basis, the matrices that
represent $A$ and $A^*$.
We will elaborate on these two points below, but
first we 
sharpen our notation.
By a {\it  basis} for 
 $V$, we mean a sequence
of vectors in $V$ that are linearly
independent and span $V$. 
We emphasize the ordering is important.
Let $v_0, v_1, \ldots, v_d$ denote
a  basis for $V$. Then  
the sequence $v_d, v_{d-1}, \ldots, v_0$ is a basis
for $V$, which we call
the {\it inversion} of 
$v_0, v_1, \ldots, v_d$.

\medskip
\noindent When we define our 24 bases, we will find they 
are related to each other according to the  
diagram in Figure 1.
In that diagram, each  vertex represents one  of the
24 bases.
For each pair  
 of bases in the 
diagram that are connected by an arc,
consider the transition matrix from
one of these bases to the other. 
The shading on the arc indicates the nature of
this transition matrix. 
If the  arc is solid, the transition matrix
is diagonal. If the arc  is dashed,
the transition matrix
is lower triangular.
If the arc is dotted, the two bases are the inversion
of one another.

\bigskip
\input psfig.sty
\centerline{\psfig{figure=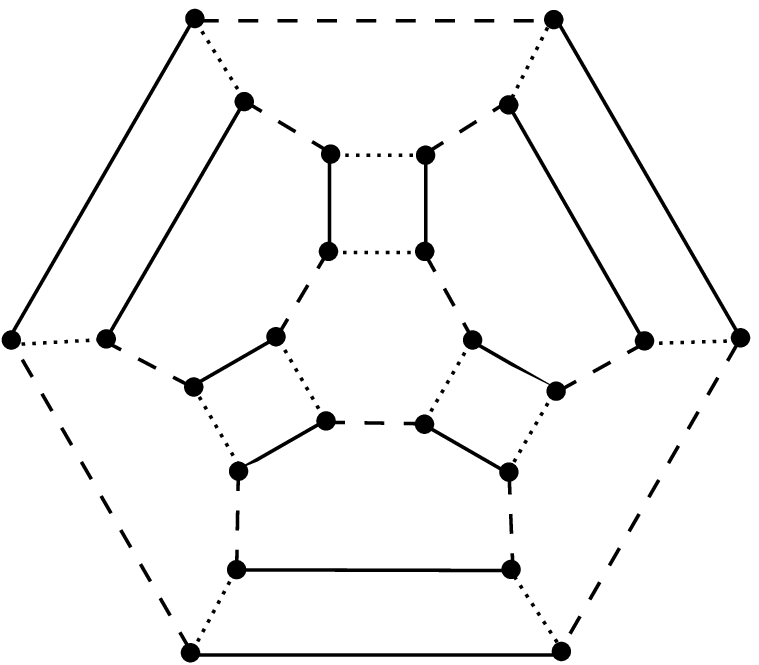,height=4.5cm}} 
\vskip .5cm
\baselineskip 7pt \par{\leftskip 1cm \rightskip 1cm \noindent 
\small Figure 1. How the 24 bases are related. Each vertex represents
one of the 24 bases.
Solid arc: transition matrix
is diagonal. 
Dashed arc:
 transition matrix
is lower triangular. 
Dotted arc: inversion.
\par} \baselineskip=\normalbaselineskip
%
%
%
%
%

\medskip
\noindent  The reader might observe 
the above diagram is 
 a Cayley graph for the 
symmetric group $S_4$. 
Apparently, there is a connection between our 24 bases and
$S_4$.
We now  make this connection explicit.

\medskip
\noindent
Let $\Omega $  denote the set consisting of four symbols 
$0,d,0^*, d^* $.
We identify the symmetric group $S_4$ with the
set of all linear orderings of $\Omega$.
For $i=1,2,3$ we define a symmetric binary relation
on $S_4$ which we call $i$-{\it adjacency}.
Each element $wxyz$ of $S_4$ 
is by definition $1$-adjacent (resp. $2$-adjacent) (resp. $3$-adjacent)
to 
$xwyz$ (resp. $wyxz$) (resp. $wxzy$) and no other elements of $S_4$.
Two elements in $S_4$ will be called {\it adjacent} whenever they are
$i$-adjacent for some $i$ $(1 \leq i \leq 3)$.
If we draw a diagram in which we
represent the elements of $S_4$ by vertices,
and for $i=1,2,3$ we represent
$i$-adjacency by
solid, dashed, and dotted arcs, respectively,
we get the diagram in Figure 1.

\medskip
\noindent 
For each element $g$ of $S_4$,
we will define a certain basis for $V$, which we denote by 
$\lbrack g \rbrack $.
We will find that
for all pairs $g,h$ of adjacent elements in $S_4$,
\begin{enumerate}
\item if $g,h$ are $1$-adjacent the transition matrix
from $\lbrack g\rbrack $ to $\lbrack h\rbrack $ is diagonal,
\item if $g,h$ are $2$-adjacent the transition matrix
from $\lbrack g\rbrack $ to $\lbrack h\rbrack $ is lower triangular,
\item if $g,h$ are $3$-adjacent then 
 $\lbrack g\rbrack $ is the inversion of $\lbrack h\rbrack $.
\end{enumerate}
When we define our 24 bases, we will find
that $A$ and $A^*$ act on them as follows.
For all $g \in S_4$,
let $A^g$ (resp. $A^{*g}$) denote
the matrix representing  $A$ (resp. $A^*$)
with respect to
$\lbrack g \rbrack $.
To describe 
$A^g$ and $A^{*g}$,
 we refer to $0^*,  d^*$ as the {\it starred} elements of $\Omega $.
Writing $g=wxyz$,  we will find
\begin{enumerate}
\item
if neither of $y,z$ are starred then
 $A^g$ is 
 diagonal
and  $A^{*g}$
is irreducible tridiagonal.
\item
if $y$ is starred but $z$ is not,
then $A^g$ is lower bidiagonal and $A^{*g}$ is upper bidiagonal.
\item
if $z$ is starred but $y$ is not,
then $A^g$ is upper bidiagonal and $A^{*g}$ is lower bidiagonal.
\item
if both of $y,z$ are starred, then  
 $A^g$ is irreducible tridiagonal 
and $A^{*g}$ is diagonal.
\end{enumerate}
(A square matrix is said to be {\it lower bidiagonal } whenever
all nonzero entries lie either on or immediately below the main diagonal.
A matrix is said to be {\it upper bidiagonal} whenever the transpose
is lower bidiagonal).

\medskip
\noindent 
For all ordered pairs $g,h$  of adjacent elements in $S_4$, we
find the entries of the transition matrix from
the basis $\lbrack g \rbrack $ to the basis 
$\lbrack h \rbrack $. We express these entries in terms
of the eigenvalues of $A$, the eigenvalues of $A^*$, and two sequences
of scalars called the
first split sequence and the second split sequence.
For all $g \in S_4$, we compute the entries of
$A^g$ and $A^{*g}$
in terms of the eigenvalues of $A$, the eigenvalues of $A^*$,
the first split sequence 
and the second split sequence.

\section{Leonard systems}
\medskip
\noindent When working with a Leonard pair, 
it is often convenient to consider a closely related
and somewhat more abstract object, which we call
a {\it Leonard system}.
In order to define this, we first make an observation about
Leonard pairs.

\begin{lemma} \cite{LS99}
\label{lem:preeverythingtalkS99}
Let $V$ denote a vector space over $\K$ with finite positive
dimension, and 
let $A, A^*$ denote a Leonard pair on $V$. Then
the eigenvalues of $A$ are distinct and contained in $\fld$.
Moreover, the eigenvalues of
$A^*$ are distinct and contained in $\fld$.
\end{lemma}

\noindent To prepare for our definition of  a Leonard system,
we recall a few concepts from elementary linear algebra. 
Let $d$  denote  a  nonnegative
integer, 
and let $\Mdf$ denote the $\fld$-algebra consisting of all
$d+1$ by $d+1$ matrices with entries in $\fld$. We
index the rows and columns by $0,1,\ldots, d$.
Let  $\alg$ 
denote a $\fld$-algebra
isomorphic to 
$\;\hbox{Mat}_{d+1}(\fld)$. 
Let $A$ denote an element of $\cal A$. By an {\it eigenvalue }
of $A$, we mean a root of the minimal polynomial of $A$.
The eigenvalues  of $A$ are contained in the algebraic closure of $\fld$.
The element $A $  will be called {\it
multiplicity-free} whenever it has $d+1$ distinct  eigenvalues,
all of which are 
in $\;\fld$.
Let $A$ denote a  multiplicity-free element of $\alg$.
Let $\theta_0, \theta_1, \ldots, \theta_d$ denote an ordering of 
the eigenvalues
of $A$, and for $0 \leq i \leq d$   put 
\begin{equation}
E_i = \prod_{{0 \leq  j \leq d}\atop
{j\not=i}} {{A-\theta_j I}\over {\theta_i-\theta_j}},
\label{eq:primiddef}
\end{equation}
where $I$ denotes the identity of $\cal A$.
By elementary linear algebra,
\begin{eqnarray}
&&AE_i = E_iA = \theta_iE_i \qquad \qquad  (0 \leq i \leq d),
\label{eq:primid1S99}
\\
&&
\quad E_iE_j = \delta_{ij}E_i \qquad \qquad (0 \leq i,j\leq d),
\label{eq:primid2S99}
\\
&&
\qquad \qquad \sum_{i=0}^d E_i = I.
\label{eq:primid3S99}
\end{eqnarray}
From this, 
one  finds  $E_0, E_1, \ldots, E_d$ is a  basis for the
subalgebra of $\alg$ generated by $A$.
We refer to $E_i$ as the {\it primitive idempotent} of
$A$ associated with $\theta_i$.
It is helpful to think of these primitive idempotents as follows. 
Let $V$ denote the irreducible left $\cal A$-module. Then
\begin{eqnarray}
V = E_0V + E_1V + \cdots + E_dV \qquad \qquad (\hbox{direct sum}).
\label{eq:VdecompS99}
\end{eqnarray}
For $0\leq i \leq d$, $E_iV$ is the (one dimensional) eigenspace of
$A$ in $V$ associated with the 
eigenvalue $\theta_i$, 
and $E_i$ acts  on $V$ as the projection onto this eigenspace.

\begin{definition} \cite{LS99}
\label{def:defls}
Let $d$  denote  a  nonnegative
integer, let $\fld $ denote a field, 
and let $\alg$ 
denote a $\;\fld$-algebra isomorphic to 
$\hbox{Mat}_{d+1}(\fld)$. 
By a {\it Leonard system} in $\;\alg$, we mean a 
sequence 
\begin{equation}
\; \Phi = (A;\,E_0,\,E_1,\,\ldots,
\,E_d;\,A^*;\,E^*_0,\,E^*_1,\,\ldots,\,E^*_d)
\label{eq:ourstartingpt}
\end{equation}
 that satisfies  (i)--(v) below.
\begin{enumerate}
\item $A$,  $\;A^*\;$ are both multiplicity-free elements in $\;\alg$.
\item $E_0,\,E_1,\,\ldots,\,E_d\;$ is an ordering of the primitive 
idempotents of $\;A$.
\item $E^*_0,\,E^*_1,\,\ldots,\,E^*_d\;$ is an ordering of the primitive 
idempotents of $\;A^*$.
\item ${\displaystyle{
E_iA^*E_j = \cases{0, &if $\;\vert i-j\vert > 1$;\cr
\not=0, &if $\;\vert i-j \vert = 1$\cr}
\qquad \qquad 
(0 \leq i,j\leq d)}}$.
\item ${\displaystyle{
 E^*_iAE^*_j = \cases{0, &if $\;\vert i-j\vert > 1$;\cr
\not=0, &if $\;\vert i-j \vert = 1$\cr}
\qquad \qquad 
(0 \leq i,j\leq d).}}$
\end{enumerate}
We refer to $d$ as the {\it diameter} of $\Phi$, and say 
$\Phi$ is {\it over } $\fld$.  We sometimes write
${\cal A} = {\cal A}(\Phi)$, $\fld= \fld(\Phi)$. 
 For notational convenience, we set $E_{-1}=0$, $E_{d+1}=0$, $
E^*_{-1}=0$, $E^*_{d+1}=0$.

\end{definition}

\noindent 
In the two lemmas below, we
explain the relationship between the notions of Leonard pair and 
Leonard system. We will use the  following notation.
Let $V$ denote a vector space over $\K$ with finite positive 
dimension. We let $\hbox{End}(V)$ denote the $\K$-algebra
consisting of all linear transformations from $V$ to $V$.
We recall  
 $\hbox{End}(V)$ is $\K$-algebra isomorphic to
$\hbox{Mat}_{d+1}(\fld)$, where $d+1 = \hbox{dim} V$.

\begin{lemma}
\label{lem:lpgivesls}
Let $V$ denote a vector space over $\K$ with finite positive
dimension.  Let $A, A^*$ denote a Leonard pair on $V$,
and observe each of $A,A^*$ is multiplicity-free by Lemma
\ref{lem:preeverythingtalkS99}.
Let $v_0, v_1, \ldots, v_d$ denote a basis for $V$ that
satisfies Definition
\ref{def:lprecall}(i).
For $0 \leq i \leq d$,
observe $v_i$ 
is an eigenvector
for $A$; 
let $\theta_i$
denote the corresponding eigenvalue, and let
$E_i$ 
denote the primitive   idempotent of $A$
associated with $\theta_i$.
Similarly,
let $v^*_0, v^*_1, \ldots, v^*_d$ denote a basis for $V$ that
satisfies Definition
\ref{def:lprecall}(ii).
For $0 \leq i \leq d$,
observe $v^*_i$ 
is an eigenvector
for $A^*$; 
let $\theta^*_i$
denote the corresponding eigenvalue, and let
$E^*_i$ 
denote the primitive   idempotent
of $A^*$ associated with $\theta^*_i$.
Then the sequence
\beast
(A;\,E_0,\,E_1,\,\ldots,
\,E_d;\,A^*;\,E^*_0,\,E^*_1,\,\ldots,\,E^*_d)
\eeast
is a Leonard system in $End(V)$.

\end{lemma}

\begin{proof} We verify the conditions (i)--(v) of
Definition
\ref{def:defls}. Condition (i) is  immediate
from Lemma
\ref{lem:preeverythingtalkS99} and the definition of
multiplicity-free. Conditions (ii), (iii) are immediate
from the construction.
Condition (iv) holds, since 
by Definition \ref{def:lprecall}(i)  
the matrix representing $A^*$ with respect
to the basis $v_0, v_1, \ldots, v_d$ is
irreducible tridiagonal.
Condition (v) holds, since
by Definition \ref{def:lprecall}(ii)  
the matrix representing $A$ with respect
to the basis $v^*_0, v^*_1, \ldots, v^*_d$ is
irreducible tridiagonal.

\end{proof}

\begin{lemma}
\label{lem:lsgiveslp}
Let $\Phi$ denote the Leonard system in
(\ref{eq:ourstartingpt}), and let $V$ 
denote the irreducible left $\cal A$-module.
For $0 \leq i\leq d$, let $v_i$ denote a
nonzero vector in $E_iV$. Then 
$v_0, v_1, \ldots, v_d$ is a basis for
$V$ with respect to which the matrix representing
$A$ is diagonal and the matrix representing $A^*$
is irreducible tridiagonal. 
For $0 \leq i\leq d$, let $v^*_i$ denote a
nonzero vector in $E^*_iV$. Then 
$v^*_0, v^*_1, \ldots, v^*_d$ is a basis for
$V$ with respect to which the matrix representing
$A^*$ is diagonal and the matrix representing $A$
is irreducible tridiagonal.
Moreover the pair $A,A^*$ is a Leonard pair on $V$.

\end{lemma}

\begin{proof} Routine.

\end{proof}

\medskip
\noindent We mention a few basics concerning Leonard systems.

\medskip
\noindent 
Let $\Phi$ 
denote the Leonard system in
(\ref{eq:ourstartingpt}), 
and let 
$\sigma :\alg \rightarrow {\cal A}'$ denote an isomorphism of
$\fld$-algebras. We write 
\begin{equation}
\Phi^{\sigma}:= 
(A^{\sigma};  E_0^{\sigma},E_1^{\sigma},\ldots, E_d^{\sigma};
A^{*\sigma}; E_0^{*\sigma}, 
 E_1^{*\sigma},\ldots,
E_d^{*\sigma}),
\label{eq:lsisoAS99}
\end{equation}
and observe 
$\Phi^{\sigma}$ 
is a Leonard  system in ${\cal A }'$.

\begin{definition}
 \cite{LS99}
\label{def:isolsS99o}
Let $\Phi$ 
and  
 $\Phi'$ 
denote Leonard systems over $\fld$.
 By an {\it isomorphism of Leonard  systems
 from $\Phi $ to $\Phi'$}, we mean an isomorphism of $\fld $-algebras
$\sigma :{\cal A}(\Phi) \rightarrow {\cal A}(\Phi')$ such 
that  $\Phi^\sigma = \Phi'$. 
The Leonard systems $\Phi $, $\Phi'$
are said to be {\it isomorphic} whenever there exists
an isomorphism of Leonard  systems from $\Phi $ to $\Phi'$. 
\end{definition}

\noindent We finish this section with a remark.
 Let $d$ denote a
nonnegative integer, and let $\cal A$ denote
a $\K$-algebra isomorphic to 
$\hbox{Mat}_{d+1}(\K)$.  
Let  $\sigma :\cal A \rightarrow 
\cal A$ denote  any map.
Then by the Skolem-Noether theorem 
\cite{CR},
$\sigma $ is an isomorphism of $\K$-algebras
if and only if there exists an invertible $S \in \cal A$ such that
$X^\sigma = S X S^{-1}$ for all  $X \in  \cal A$.

\section{The structure of a Leonard system}

\noindent Let $\Phi$ denote the  Leonard system in 
(\ref{eq:ourstartingpt}). In this section, we show
there does not exist an
isomorphism of Leonard systems from $\Phi$ to
itself, other than the identity map.
We begin with a lemma.

\begin{lemma}
\label{eq:lsmatbasis}
Let $\Phi$ denote the Leonard system
in 
(\ref{eq:ourstartingpt}). Then the elements
\begin{equation}
A^rE^*_0A^s \qquad \qquad (0 \leq r,s\leq d)
\label{eq:lpbasis}
\end{equation}
form a basis for $\cal A$.
\end{lemma}

\begin{proof}  
The number of elements in 
(\ref{eq:lpbasis}) equals $(d+1)^2$, and this number is the dimension
of $\cal A$. Therefore
it suffices
to show 
the elements
in (\ref{eq:lpbasis})
are linearly independent. To do this,
we represent 
the elements 
in (\ref{eq:lpbasis}) by matrices. 
Let $V$ denote the irreducible left $\cal A$-module.
For $0 \leq i \leq d$, let $v_i$ denote a nonzero vector in 
$E^*_iV$, and observe $v_0, v_1,\ldots, v_d$ is a basis for
$V$. For the purposes of this proof, let us identify each element
of $\cal A$ with the matrix in 
$\hbox{Mat}_{d+1}(\K)$ that
represents it with respect to the
basis $v_0, v_1, \ldots, v_d$. Adopting this point of view 
$A$ is irreducible tridiagonal and $A^*$ is diagonal. 
For $0 \leq r,s\leq d$ we show  the
entries  of $A^rE^*_0A^s$ satisfy 
\begin{equation}
(A^rE^*_0A^s)_{ij} = 
 \cases{0, 
&$\qquad $if $\quad i>r \quad$ or $\quad j>s $ ;\cr
\not=0,  & $\qquad $if $\quad i=r\quad $ and $\quad j=s$\cr} 
\qquad \qquad (0 \leq i,j\leq d).
\label{eq:lowert}
\end{equation}
Observe that for $0 \leq i,j\leq d$, 
the $ij^{\hbox{th}}$ entry of $E^*_0$ is 
one if both $i=0,j=0$, and zero otherwise. 
From this we find
\begin{equation}
(A^rE^*_0A^s)_{ij} =  A^r_{i0} A^s_{0j} \qquad \qquad (0 \leq i,j\leq d).
\label{eq:redent}
\end{equation}
Since $A$ is irreducible tridiagonal, we find that for 
$0 \leq i \leq d$, the $i0^{\hbox{th}}$ entry of $A^r$ is zero if $i>r$, and  nonzero if
$i=r$. Similarly for $0 \leq j\leq d$, the $0j^{\hbox{th}}$ entry
of $A^s$ is zero 
if $j>s$, and nonzero if $j=s$. Combining these facts with 
(\ref{eq:redent}) we routinely obtain
(\ref{eq:lowert}) and it follows  the elements 
(\ref{eq:lpbasis}) are linearly independent.
Apparently the elements 
(\ref{eq:lpbasis}) form a basis  for $\cal A $, as desired.

\end{proof}

\begin{corollary} 
\label{cor:genset}
Let $\Phi$ denote the Leonard system
in 
(\ref{eq:ourstartingpt}). Then 
the elements 
$A, E^*_0$ together generate $\cal A$. Moreover,
the elements $A,A^*$ together generate $\cal A$.

\end{corollary}

\begin{proof} 
The first assertion is immediate from Lemma 
\ref{eq:lsmatbasis}.  
The second assertion follows from the first and
the observation that $E^*_0$ is a polynomial in $A^*$. 

\end{proof}

\noindent We mention a useful consequence of Corollary
\ref{cor:genset}.

\begin{corollary}
\label{cor:rig}
Let $\Phi$ denote 
the Leonard system
(\ref{eq:ourstartingpt}),
and let $X$ denote an element in $\cal A$ that  commutes
with both $A$ and $A^*$. Then $X$ is a scalar multiple
of the identity.
Put another way, there does not exist 
an isomorphism of Leonard systems
from $\Phi$ to itself, other than the identity map.

\end{corollary}

\begin{proof} Since $A,A^*$ together generate $\cal A$, we find
$X$ commutes with everything in $\cal A$. Now
$X$ is a scalar multiple of the identity by elementary
linear algebra. The last assertion follows in view of
our remark at the end of  Section 2. 

\end{proof}

\noindent We mention an implication of Lemma 
\ref{eq:lsmatbasis} that will be useful later in the paper.

\begin{lemma}
\label{lem:altbase}
Let $\Phi$ denote the Leonard system
in 
(\ref{eq:ourstartingpt}). Let $\cal D$ denote
the subalgebra of $\cal A$ generated by $A$,
and observe $\cal D$ has dimension $d+1$ since
$A$ is multiplicity-free. Let $X_0, X_1, \ldots, X_d$
denote a basis for $\cal D$. Then the elements
\begin{equation}
X_r E^*_0 X_s \qquad \qquad (0 \leq r,s\leq d)
\label{altbase}
\end{equation}
form a basis for $\cal A$.

\end{lemma}

\begin{proof} The number of elements in    
(\ref{altbase}) is $(d+1)^2$, and this number is the dimension
of $\cal A$. Therefore it suffices to show the elements
(\ref{altbase}) span $\cal A$. 
But this is immediate from Lemma
\ref{eq:lsmatbasis}, and since 
each  element 
in (\ref{eq:lpbasis}) is contained in the span of
the elements (\ref{altbase}).

\end{proof}

\begin{corollary}
\label{cor:eibasis}
Let $\Phi$ denote the Leonard system
in 
(\ref{eq:ourstartingpt}).
Then the elements
\begin{equation}
E_r E^*_0 E_s \qquad \qquad (0 \leq r,s\leq d)
\label{eq:eibase}
\end{equation}
form a basis for $\cal A$.

\end{corollary}

\begin{proof}  Immediate from  
Lemma \ref{lem:altbase}, with  $X_i = E_i$ for $0 \leq i \leq d$.

\end{proof}

%
%
%

\section{The relatives of a Leonard system}

\medskip
\noindent A given Leonard system  can be modified in  several
ways to get a new Leonard system. For instance, 
let $\Phi$ 
 denote the Leonard system in
(\ref{eq:ourstartingpt}).
%
Then each of the following three sequences is a Leonard system
in $\cal A$.
\begin{eqnarray}
 \;\Phi^*&:=& (A^*; E^*_0,E^*_1,\ldots,E^*_d;A;E_0,E_1, \ldots,E_d),
\label{eq:lsdualS99}
\\
\Phi^{\downarrow}&:=& (A; E_0,E_1,\ldots,E_d;A^*;E^*_d,E^*_{d-1}, \ldots,E^*_0),
\label{eq:lsinvertS99}
\\
\Phi^{\Downarrow} 
&:=& (A; E_d,E_{d-1},\ldots,E_0;A^*;E^*_0,E^*_1, \ldots,E^*_d).
\label{eq:lsdualinvertS99}
\end{eqnarray}
 We refer to $\Phi^*$
(resp.  
 $\Phi^\downarrow$)   
(resp.  
 $\Phi^\Downarrow$) 
 as the  
{\it dual} 
(resp. {\it first inversion})
(resp.  {\it second inversion}) of  $\Phi$.
Viewing $*, \downarrow, \Downarrow$
as permutations on the set of all Leonard systems,
\begin{eqnarray}
&&\qquad \qquad \qquad  *^2 \;=\;  
\downarrow^2\;= \;
\Downarrow^2 \;=\;1,
\qquad \quad 
\label{eq:deightrelationsAS99}
\\
&&\Downarrow *\; 
=\;
* \downarrow,\qquad \qquad   
\downarrow *\; 
=\;
* \Downarrow,\qquad \qquad   
\downarrow \Downarrow \; = \;
\Downarrow \downarrow.
\qquad \quad 
\label{eq:deightrelationsBS99}
\end{eqnarray}
The group generated by symbols 
$*, \downarrow, \Downarrow $ subject to the relations
(\ref{eq:deightrelationsAS99}),
(\ref{eq:deightrelationsBS99})
is the dihedral group $D_4$.  
We recall $D_4$ is the group of symmetries of a square,
and has 8 elements.
Apparently $*, \downarrow, \Downarrow $ induce an action of 
 $D_4$ on the set of all Leonard systems.
Two Leonard systems will be called {\it relatives} whenever they
are in the same orbit of this $D_4$ action.
The relatives of $\Phi$ are as follows:
\medskip

\centerline{
\begin{tabular}[t]{c|c}
        name &relative \\ \hline 
        $\Phi$ & $(A;E_0,E_1,\ldots,E_d;A^*;E^*_0,E^*_1,\ldots,E^*_d)$   \\ 
        $\Phi^\downarrow$ &
         $(A;E_0,E_1,\ldots,E_d;A^*;E^*_d,E^*_{d-1},\ldots,E^*_0)$   \\ 
        $\Phi^\Downarrow$ &
         $(A;E_d,E_{d-1},\ldots,E_0;A^*;E^*_0,E^*_1,\ldots,E^*_d)$   \\ 
        $\Phi^{\downarrow \Downarrow}$ &
         $(A;E_d,E_{d-1},\ldots,E_0;A^*;E^*_d,E^*_{d-1},\ldots,E^*_0)$   \\ 
	$\Phi^*$ &
        $(A^*;E^*_0,E^*_1,\ldots,E^*_d;A;E_0,E_1,\ldots,E_d)$   \\ 
        $\Phi^{\downarrow *}$ &
	 $(A^*;E^*_d,E^*_{d-1},\ldots,E^*_0;  
         A;E_0,E_1,\ldots,E_d)$ \\
        $\Phi^{\Downarrow *}$ &
	 $(A^*;E^*_0,E^*_1,\ldots,E^*_d;    
         A;E_d,E_{d-1},\ldots,E_0)$ \\ 
	$\Phi^{\downarrow \Downarrow *}$ &
	 $(A^*;E^*_d,E^*_{d-1},\ldots,E^*_0;    
         A;E_d,E_{d-1},\ldots,E_0)$
	\end{tabular}}
\medskip
\noindent 
We remark there may be some isomorphisms among the above Leonard
systems.

\medskip 
\noindent We finish this section by  recalling some parameters
that will help us describe 
a given Leonard system.

\begin{definition}
 \cite{LS99}
\label{def:evseq}
Let $\Phi$ denote the Leonard system in 
(\ref{eq:ourstartingpt}).
For $0 \leq i \leq d$, 
we let $\theta_i $ (resp. $\theta^*_i$) denote the eigenvalue
of $A$ (resp. $A^*$) associated with $E_i$ (resp. $E^*_i$).
We refer to  $\theta_0, \theta_1, \ldots, \theta_d$ as the 
eigenvalue sequence of $\Phi$.
We refer to  $\theta^*_0, \theta^*_1, \ldots, \theta^*_d$ as the 
dual eigenvalue sequence of $\Phi$. We observe 
 $\theta_0, \theta_1, \ldots, \theta_d$ are mutually distinct
 and contained in $\K$. Similarly
  $\theta^*_0, \theta^*_1, \ldots, \theta^*_d$  
 are mutually distinct
 and contained in $\K$. 
\end{definition}

\section{The standard basis and the split basis}

\medskip
\noindent Let $\Phi$ denote the Leonard system
in 
(\ref{eq:ourstartingpt}), and let $V$ denote the irreducible
left $\cal A$-module.
As  we mentioned earlier, we will obtain 24  
bases for $V$. One way to view our construction is as follows.
Using $\Phi$
we define three  bases
for $V$, called the {\it $\Phi$-standard basis},
the {\it $\Phi$-split basis}, and the {\it $\Phi$-inverted split
basis}. In each of the three cases, the basis is 
defined up  
to multiplication of each element by the same nonzero scalar in
$\K$.
Our set of 
 24 bases will consist of 
a $\Psi$-standard basis,
a $\Psi$-split basis,
and a $\Psi$-inverted split basis for each relative $\Psi$ of $\Phi$.

\medskip
\noindent
We now define the notion of a  {\it standard basis}.

\begin{lemma}
\label{lem:stbasis}
Let $\Phi$ denote the Leonard system
in 
(\ref{eq:ourstartingpt}),
and let $V$ denote the irreducible left $\cal A$-module.
Let 
$u$ denote a nonzero element of $E^*_0V$. Then
for $0 \leq i \leq d$, the element $E_iu$ is nonzero
and hence a basis
for $E_iV$. Moreover the sequence
\begin{equation}
E_0u, E_1u, \ldots, E_du
\label{eq:stbasisint}
\end{equation}
is a basis for $V$.

\end{lemma}

\begin{proof} Let the integer $i$ be given.
Recall $E^*_0V$ has dimension 1, and
$u$ is a nonzero vector in 
 $E^*_0V$,  
so $u$ spans  $E^*_0V$. Apparently 
$E_iu$ spans $E_iE^*_0V$. Observe $E_iE^*_0$ is nonzero by
Corollary
\ref{cor:eibasis}
so $E_iE^*_0V$ is nonzero. Apparently 
$E_iu$ is nonzero,
and is therefore a  basis for $E_iV$, as desired.
The sequence
(\ref{eq:stbasisint}) is a basis for $V$ in view of 
(\ref{eq:VdecompS99}).

\end{proof}

\begin{definition}
\label{def:stbasis}
Let $\Phi$ denote the Leonard system
in 
(\ref{eq:ourstartingpt}), and let $V$ denote the 
irreducible left $\cal A$-module.
By a $\Phi$-standard basis for $V$,
we mean
a sequence  
(\ref{eq:stbasisint}), where $u$ is a  nonzero vector
in $E^*_0V$. 
When the identity of $\Phi$ is clear,
we will occasionaly 
 speak of a standard basis
instead of a $\Phi$-standard basis.

\end{definition}

\noindent Let $\Phi$ denote the Leonard system
in 
(\ref{eq:ourstartingpt}), and let $V$ denote the 
irreducible left $\cal A$-module.
With  respect to any  $\Phi$-standard basis for $V$,
the matrix representing $A$ is
\beast
\hbox{diag}(\theta_0, \theta_1, \ldots, \theta_d),
\eeast
where the $\theta_i$ are from
Definition \ref{def:evseq}. 
Moreover, by Lemma 
\ref{lem:lsgiveslp},
the matrix representing $A^*$ 
is irreducible tridiagonal. 
We will work out the entries of this tridiagonal  matrix
in due course, but it is convenient  to wait until
after we have introduced some more bases. 
For those who wish to skip ahead,
the  entries of this tridiagonal matrix
can be found in the second table 
of Theorem
\ref{thm:repa}, row 1.

\medskip
\noindent We now define the notion
of a 
 {\it split basis}.
In the process, we will recall two
sequences of scalars which we will find useful. 
These sequences are called the
{\it first split sequence} of $\Phi$  and the  {\it 
second split sequence} of $\Phi$.

\medskip
\noindent
In order to define a split basis,
 we review some results of \cite{TD00}, \cite{LS99}.
Let $\Phi$ denote the Leonard system
in (\ref{eq:ourstartingpt})
and let $V$ denote the irreducible left $\cal A$-module.
For $0 \leq i \leq d$ we define
\begin{equation}
U_i = 
(E^*_0V + E^*_1V + \cdots + E^*_iV)\cap (E_iV + E_{i+1}V + \cdots + E_dV).
\label{eq:defui}
\end{equation}
We showed in \cite{TD00}
that each of $U_0, U_1, \ldots, U_d$ has dimension 1, and that
\begin{equation}
V = U_0 + U_1 + \cdots + U_d \qquad \qquad (\hbox{direct sum}).
\label{eq:splitdec}
\end{equation}
Moreover,
\begin{eqnarray}
U_0 + U_1 + \cdots + U_i &=& E^*_0V + E^*_1V + \cdots + E^*_iV,
\label{eq:split1}
\\
U_i + U_{i+1} + \cdots + U_d &=& E_iV + E_{i+1}V + \cdots + E_dV
\label{eq:split2}
\end{eqnarray}
for $0 \leq i \leq d$. 
The elements $A$ and $A^*$ act on the $U_i$ as follows.
We showed in \cite{LS99} that
\begin{eqnarray}
(A-\theta_i I)U_i &=& U_{i+1} \qquad (0 \leq i \leq d-1),
\qquad (A-\theta_d I)U_d = 0,
\label{eq:raise}
\\
(A^*-\theta^*_i I)U_i &=& U_{i-1} \qquad (1 \leq i \leq d),
\qquad (A^*-\theta^*_0 I)U_0 = 0,
\label{eq:lower}
\end{eqnarray}
where the $\theta_i, \theta^*_i$ are from
Definition \ref{def:evseq}. 
Pick an integer $i$ $(1 \leq i \leq d)$. By
(\ref{eq:lower}) we find 
$(A^*-\theta^*_i I)U_i = U_{i-1}$ and by
(\ref{eq:raise}) we find 
$(A-\theta_{i-1} I)U_{i-1} = U_i$. Apparently $U_i$
is an eigenspace for
$(A-\theta_{i-1}I)(A^*-\theta^*_i I)$,
and the corresponding eigenvalue is a nonzero element of $\K$.
We denote this 
eigenvalue by $\varphi_i$.
We refer to the sequence $\varphi_1, \varphi_2,
\ldots, \varphi_d$ as the {\it first split sequence} of 
$\Phi$. We let $\phi_1, \phi_2, \ldots, \phi_d$ denote
the first split sequence for $\Phi^{\Downarrow}$, and
call  this the {\it  second split sequence} of $\Phi$. For notational
convenience, we define $\varphi_0=0$, $\varphi_{d+1} = 0$,
 $\phi_0=0$, $\phi_{d+1} = 0$.

\medskip
\noindent
We obtain our split basis as follows.
Setting $i=0$ in 
(\ref{eq:split1}), we 
find $U_0=E^*_0V$. Combining this with
(\ref{eq:raise}), we find  
\begin{equation}
U_i = (A-\theta_0 I)(A - \theta_1 I)\cdots (A-\theta_{i-1} I )E^*_0V
\qquad \qquad (0 \leq i\leq d).
\label{eq:uialt}
\end{equation}
Let $u$ denote a nonzero vector in $E^*_0V$. 
From 
(\ref{eq:uialt}) we find that for 
$0 \leq i \leq d$, the vector
$(A-\theta_0I)\cdots (A-\theta_{i-1}I)u $ is a basis for
$U_i$. From this and 
(\ref{eq:splitdec}) we find the sequence
\begin{equation}
(A-\theta_0 I)(A-\theta_1 I)\cdots (A-\theta_{i-1}I)u  \qquad \qquad
(0 \leq i \leq d)
\label{eq:basis1}
\end{equation}
is a basis for $V$.

\begin{definition}
\label{def:splitbase}
Let $\Phi$ denote the Leonard system
in 
(\ref{eq:ourstartingpt}), and let $V$ denote the 
irreducible left $\cal A$-module.
By a {\it $\Phi$-split basis} for $V$, 
 we mean a sequence
(\ref{eq:basis1}), where $u$ is a nonzero vector in $E^*_0V$. 
When the identity of $\Phi$ is clear,
we will occasionaly 
 speak of a split basis
instead of a $\Phi$-split basis.

\end{definition}

Let $\Phi$ denote the Leonard system
in 
(\ref{eq:ourstartingpt}), and let $V$ denote the 
irreducible left $\cal A$-module. From 
(\ref{eq:basis1}) and 
the lines 
below 
(\ref{eq:lower}),
 we find that 
 with respect to any $\Phi$-split basis for $V$, the matrices
representing $A$ and $A^*$ are 
\begin{equation}
\left(
\begin{array}{c c c c c c}
\theta_0 & & & & & {\bf 0} \\
1 & \theta_1 &  & & & \\
& 1 & \theta_2 &  & & \\
& & \cdot & \cdot &  &  \\
& & & \cdot & \cdot &  \\
{\bf 0}& & & & 1 & \theta_d
\end{array}
\right),
\qquad  \quad 
\left(
\begin{array}{c c c c c c}
\theta^*_0 &\varphi_1 & & & & {\bf 0} \\
 & \theta^*_1 & \varphi_2 & & & \\
&  & \theta^*_2 & \cdot & & \\
& &  & \cdot & \cdot &  \\
& & &  & \cdot & \varphi_d \\
{\bf 0}& & & &  & \theta^*_d
\end{array}
\right),
\label{eq:matrepaastar}
\end{equation}
respectively.

\medskip
\noindent 
We now define the notion of an {\it inverted 
split basis}. As its name implies, an inverted split
basis is nothing but the inversion of a split basis. 
To be concrete, we make the following definition.

\begin{definition}
\label{def:invsplit}
Let $\Phi$ denote the Leonard system
in 
(\ref{eq:ourstartingpt}), and let $V$ denote the 
irreducible left $\cal A$-module.
By a {\it $\Phi$-inverted split basis} for $V$,
we mean 
a sequence
\begin{equation}
(A-\theta_0 I)(A-\theta_1 I)\cdots (A-\theta_{d-i-1}I)u  \qquad \qquad
(0 \leq i \leq d),
\label{eq:invbasis1}
\end{equation}
where $u$ is a nonzero vector in $E^*_0V$.
When the identity of $\Phi$ is clear,
we will occasionaly 
 speak of an  inverted split basis
instead of a $\Phi$-inverted split basis.

\end{definition}
Let $\Phi$ denote the Leonard system
in 
(\ref{eq:ourstartingpt}), and let $V$ denote the 
irreducible left $\cal A$-module. Combining
(\ref{eq:matrepaastar})
with Definition
\ref{def:invsplit}, we find that 
with respect to any $\Phi$-inverted split basis for $V$,
the matrices representing
$A $ and $A^*$ are
\begin{equation}
\left(
\begin{array}{c c c c c c}
\theta_d & 1 & & & & {\bf 0} \\
 & \theta_{d-1} & 1 & & & \\
&  & \theta_{d-2} & \cdot & & \\
& &  & \cdot & \cdot &  \\
& & &  & \cdot & 1 \\
{\bf 0}& & & &  & \theta_0
\end{array}
\right),
\qquad  \quad 
\left(
\begin{array}{c c c c c c}
\theta^*_d & & & & & {\bf 0} \\
\varphi_d & \theta^*_{d-1} &  & & & \\
& \varphi_{d-1} & \theta^*_{d-2} &  & & \\
& & \cdot & \cdot &  &  \\
& & & \cdot & \cdot &  \\
{\bf 0}& & & & \varphi_1 & \theta^*_0
\end{array}
\right),
\label{eq:invmat}
\end{equation}
respectively.


\section{A classification of Leonard systems}

\medskip
\noindent  In the preceeding section, we defined the
first and second split sequence of a Leonard system.
The scalars involved in these sequences are related by many equations.
To describe these relationships,
we recall our classification of Leonard systems. 

\begin{theorem} \cite{LS99}
\label{thm:classls} Let 
$d$ denote a nonnegative integer,  
and let 
\begin{eqnarray}
&&\theta_0, \theta_1, \ldots, \theta_d; \qquad \qquad \; 
\theta^*_0, \theta^*_1, \ldots, \theta^*_d; 
\label{eq:paramlist1S99}
\\
&&\varphi_1, \varphi_2, \ldots, \varphi_d;  \qquad \qquad 
\phi_1, \phi_2, \ldots, \phi_d \qquad \quad
\label{eq:paramlist2S99}
\end{eqnarray}
denote scalars in $\fld$. 
Then there exists  a Leonard system $\Phi$ over $\fld$  with 
eigenvalue sequence $\theta_0, \theta_1, \ldots, \theta_d$, 
dual eigenvalue sequence  
$\theta^*_0, \theta^*_1, \ldots, \theta^*_d$, first split sequence
$ \varphi_1, \varphi_2, \ldots, \varphi_d $, and second split sequence
$\phi_1, \phi_2, \ldots, \phi_d$ if and only if 
(i)--(v) hold below.
\begin{enumerate}
\item $ \varphi_i \not=0, \qquad \phi_i\not=0 \qquad \qquad \qquad\qquad (1 \leq i \leq d)$,
\item $ \theta_i\not=\theta_j,\qquad  \theta^*_i\not=\theta^*_j\qquad $
if $\;\;i\not=j,\qquad \qquad \qquad (0 \leq i,j\leq d)$,
\item $ {\displaystyle{ \varphi_i = \phi_1 \sum_{h=0}^{i-1}
{{\theta_h-\theta_{d-h}}\over{\theta_0-\theta_d}} 
\;+\;(\theta^*_i-\theta^*_0)(\theta_{i-1}-\theta_d) \qquad \;\;(1 \leq i \leq d)}}$,
\item $ {\displaystyle{ \phi_i = \varphi_1 \sum_{h=0}^{i-1}
{{\theta_h-\theta_{d-h}}\over{\theta_0-\theta_d}} 
\;+\;(\theta^*_i-\theta^*_0)(\theta_{d-i+1}-\theta_0) \qquad (1 \leq i \leq d)}}$,
\item The expressions
\begin{equation}
{{\theta_{i-2}-\theta_{i+1}}\over {\theta_{i-1}-\theta_i}},\qquad \qquad  
 {{\theta^*_{i-2}-\theta^*_{i+1}}\over {\theta^*_{i-1}-\theta^*_i}} 
 \qquad  \qquad 
\label{eq:betaplusone}
\end{equation} 
 are equal and independent of $i$ for $\;2\leq i \leq d-1$.  
\end{enumerate}
Moreover, if (i)--(v) hold 
above then $\Phi$ is unique up to isomorphism of Leonard systems.
\end{theorem}
\noindent We view Theorem
\ref{thm:classls} 
as a linear algebraic version
of a theorem of Leonard
 \cite{Leodual}, \cite[p260]{BanIto}. This is discussed in \cite{LS99}.

%
%
%
\medskip
\noindent
One nice feature of 
the parameter sequences
(\ref{eq:paramlist1S99}), 
(\ref{eq:paramlist2S99})
 is that they are modified 
in a simple way as one passes from a given Leonard system
to a relative.  Our result is the following.

\begin{theorem} \cite{LS99}
\label{thm:phimod}
Let $\Phi$ denote a Leonard system, with 
eigenvalue sequence 
$\theta_0, \theta_1, \ldots, \theta_d$,
dual eigenvalue sequence 
$\theta^*_0, \theta^*_1, \ldots, \theta^*_d $,  
first split sequence
 $\varphi_1, \varphi_2,
\ldots, \varphi_d$ and second split sequence
$\phi_1, \phi_2, \ldots, \phi_d$. Then  (i)--(iii) hold below.
\begin{enumerate}
\item  The eigenvalue and dual eigenvalue sequences of $\Phi^*$
are given by
$\theta^*_0, \theta^*_1, \ldots, \theta^*_d $ and
$\theta_0, \theta_1, \ldots, \theta_d$, respectively.
The first  and second split sequences of $\Phi^*$ are
given by
 $\varphi_1, \varphi_2, \ldots, \varphi_d$ and 
$\phi_d, \phi_{d-1}, \ldots, \phi_1$, respectively.

\item 
 The eigenvalue and dual eigenvalue sequences of $\Phi^{\downarrow}$
are given by
$\theta_0, \theta_1, \ldots, \theta_d$ and 
$\theta^*_d, \theta^*_{d-1}, \ldots, \theta^*_0 $,
respectively.
The first  and second split sequences of $\Phi^\downarrow$ are
given by
$\phi_d, \phi_{d-1}, \ldots, \phi_1$ and 
 $\varphi_d, \varphi_{d-1}, \ldots, \varphi_1$,
respectively.

\item  The eigenvalue and dual eigenvalue sequences of $\Phi^\Downarrow$
are given by
$\theta_d, \theta_{d-1}, \ldots, \theta_0$ and
$\theta^*_0, \theta^*_1, \ldots, \theta^*_d $,
respectively.
The first and second split sequences of $\Phi^\Downarrow$ are
given by
$\phi_1, \phi_2, \ldots, \phi_d$  and
 $\varphi_1, \varphi_2,
\ldots, \varphi_d$,
respectively.

\end{enumerate}
\end{theorem}

\section{Four flags for $V$}

\noindent Let $\Phi$ denote the Leonard system in
(\ref{eq:ourstartingpt}),
and let $V$ denote the 
irreducible left ${\cal A}$-module.
We mentioned earlier we will obtain 24 bases
for $V$. 
In Section 5 we described these bases to some extent, but we
stopped short of displaying them. The reason is 
we wish to first introduce our labelling scheme.
As we indicated in Section 1, it is appropriate to
label our bases with elements of $S_4$.
We begin with a definition.

\begin{definition}
\label{def:S4interp}
Let $\Omega $  denote the set consisting of four symbols 
$0,d,0^*, d^* $.
We identify the symmetric group $S_4$ with the
set of all linear orderings of $\Omega$.
For $i=1,2,3$ we define a symmetric binary relation
on $S_4$ which we call $i$-adjacency.
An element $wxyz$ of $S_4$ 
is by definition $1$-adjacent (resp. $2$-adjacent) (resp. $3$-adjacent)
to 
$xwyz$ (resp. $wyxz$) (resp. $wxzy$) and no other elements of $S_4$.
Two elements in $S_4$ will be called adjacent whenever they are
$i$-adjacent for some $i$ $(1 \leq i \leq 3)$.
\end{definition}

\medskip
\noindent Let $\Phi$ denote the Leonard system in
(\ref{eq:ourstartingpt}),
and let $V$ denote the 
irreducible left ${\cal A}$-module.
We recall the notion of a flag on $V$.
By a {\it flag} on $V$, we mean a 
sequence $F_0, F_1, \ldots, F_d$ consisting 
of subspaces of $V$ such that 
$F_{i-1}\subseteq F_i$ for $1 \leq i\leq d$  and  such that
$F_i$ has dimension $i+1$ for
$0 \leq i\leq d$.
We refer to $F_i$ as the $i^{\hbox{th}}$ {\it component} of
the flag. 

\medskip
\noindent The following construction yields a
flag on $V$. To explain the construction, we make a definition. 
 By a {\it  decomposition}  of $V$, we mean
a sequence 
$L_0, L_1,\ldots, L_d$ consisting of  
1-dimensional subspaces of $V$
 such that
\begin{equation}
V = L_0 + L_1 + \cdots + L_d \qquad \qquad (\hbox{direct sum}).
\end{equation}
Let 
$L_0, L_1,\ldots, L_d$ denote a decomposition of $V$, and set
\beast
F_i &=& L_0 + L_1 + \cdots + L_i  
\eeast
for $0 \leq i \leq d$. Then the sequence
 $F_0, F_1, \ldots, F_d$ is a flag on $V$.

\medskip
\noindent We will be concerned with 
the following four flags on $V$.

\begin{definition}
\label{def:fourflags}
Let $\Phi$ denote the Leonard system  in
(\ref{eq:ourstartingpt}), and let $V$ denote the irreducible
left $\cal A$-module. Let the set $\Omega $ be
as in 
Definition
\ref{def:S4interp}.
For each element $z \in \Omega$,
we define a flag on $V$, which we denote by $\lbrack z \rbrack $.
To define this flag, we display its $i^{\hbox{th}}$ component
for $0  \leq i \leq d$.

\bigskip
\centerline{
\begin{tabular}[t]{c|c}
        $z$ &$i^{\hbox{th}}$ component of the flag $\lbrack z \rbrack $\\ \hline 
        $0$ &
	$E_0V+E_1V+\cdots + E_iV$  
	\\ 
        $d$ & 	
        $E_dV + E_{d-1}V + \cdots + E_{d-i}V$ 	
	\\ 
        $0^*$ &
        $E^*_0V+E^*_1V+ \cdots + E^*_iV$  
        \\
	$d^*$ &
	$E^*_dV+ E^*_{d-1}V + \cdots +E^*_{d-i}V$ 
	\end{tabular}}
\medskip
\noindent 

\end{definition}

\medskip
\noindent 
 Let $\Phi$ denote the Leonard system in
(\ref{eq:ourstartingpt}),
and let $V$ denote the 
irreducible left ${\cal A}$-module.
We recall what it means for two flags on $V$ to be opposite.
Suppose we are given two flags on $V$, denoted
 $F_0, F_1, \ldots, F_d$
and $G_0, G_1, \ldots, G_d$. These flags are
said to be {\it opposite} whenever 
\begin{equation}
F_i \cap G_j = 0 \quad {\hbox{if}} \quad  i+j < d, \qquad \qquad (0 \leq i,j\leq d).
\label{eq:defopp}
\end{equation}

\medskip
\noindent Given a decomposition of $V$, the following construction
yields an
ordered pair of opposite
flags on $V$. Let
$L_0, L_1,\ldots, L_d$ denote a decomposition of $V$, and set
\begin{eqnarray}
F_i &=& L_0 + L_1 + \cdots + L_i, 
\nonumber\\
G_i &=& L_d + L_{d-1} + \cdots + L_{d-i}
\label{eq:getg} 
\end{eqnarray}
for $0 \leq i \leq d$. Then the sequences
 $F_0, F_1, \ldots, F_d$
and $G_0, G_1, \ldots, G_d$ are opposite flags on $V$.

\medskip
\noindent We now turn things around.
Given an ordered pair of opposite flags on $V$, 
the following construction yields a decomposition of
$V$.
Suppose we are given an ordered pair of opposite flags on $V$,
denoted
 $F_0, F_1, \ldots, F_d$
and $G_0, G_1, \ldots, G_d$.
Set
\begin{equation}
L_i = F_i \cap G_{d-i} \qquad \qquad (0\leq i \leq d).
\label{eq:fgtoell}
\end{equation}
Then the sequence 
$L_0, L_1,\ldots, L_d$ is a decomposition of $V$. 
%

\medskip
\noindent 
Let $D$ denote the set consisting of all decompositions 
of $V$, and let $F$ denote the set consisting of all  
ordered pairs of opposite flags on $V$.
In the previous two paragraphs, we 
defined a map from $D$ to $F$ and a map from
$F$ to $D$. It is routine to show 
these maps are inverses of one another
 \cite{Ronan}.
In particular, each of these maps is a bijection.
%
%
%

\medskip
\noindent We now return to the Leonard system $\Phi$.

\begin{theorem}
\label{thm:fourf}
The four flags in Definition  
\ref{def:fourflags} are mutually opposite.

\end{theorem}
\begin{proof} It is immediate from 
the construction that
flags $\lbrack 0 \rbrack, \lbrack d \rbrack$ are opposite,
and that flags $\lbrack 0^*\rbrack, \lbrack d^* \rbrack $ are opposite.
We now show the flags $\lbrack 0^*\rbrack , 
\lbrack d \rbrack $ are opposite. For $0 \leq i \leq d$,
let $U_i$ denote the subspace of $V$ from 
(\ref{eq:defui}). By the two lines following
(\ref{eq:defui}), we find the sequence  $U_0, U_1, \ldots, U_d$
is a decomposition of $V$.
By (\ref{eq:split1}), 
(\ref{eq:split2})
 and the line following
(\ref{eq:getg}),
we find the flags 
$\lbrack 0^* \rbrack, \lbrack d \rbrack $ are opposite.
Applying
this fact to the relatives of $\Phi$, we see that
the remaining pairs of flags
in 
Definition \ref{def:fourflags}
are opposite. 

\end{proof}

\section{Twelve decompositions of $V$}

\medskip
\noindent Let $\Phi$ denote the Leonard system in
(\ref{eq:ourstartingpt}),
 let $V$ denote the 
irreducible left ${\cal A}$-module,
and let the set  $\Omega $ be as in Definition
\ref{def:S4interp}. 
In this section, we obtain for each ordered
pair $yz$ of distinct  elements in $\Omega$, a decomposition
of $V$ which we  denote by $\lbrack yz \rbrack $.

\begin{definition}
\label{def:byzb}
 Let $\Phi$ denote the Leonard system in
(\ref{eq:ourstartingpt}),  
 let $V$ denote the 
irreducible left ${\cal A}$-module,
and let the set  $\Omega $ be as in Definition
\ref{def:S4interp}. 
Let $yz$ denote an ordered pair of  distinct elements in $\Omega $.
Set 
\beast
L_i = F_i \cap G_{d-i} \qquad \qquad ( 0 \leq i \leq d),
\eeast
where $F_j$ (resp. $G_j$) denotes the $j^{\hbox{th}}$ component
of the flag $\lbrack y \rbrack $ (resp. $\lbrack z \rbrack $)  for
$0 \leq j \leq d$.
Recall $\lbrack y\rbrack$ and $\lbrack z \rbrack$ are opposite,
so the sequence 
$L_0, L_1,\ldots, L_d$ is a decomposition of $V$. We denote
this decomposition by $\lbrack yz \rbrack $.

\end{definition}

With reference to 
Definition
\ref{def:byzb}, we remark on the  
difference   between $\lbrack yz \rbrack $ and $ \lbrack zy \rbrack $.
To do this, we use the following notation.
Let $L_0, L_1, \ldots, L_d$ denote a decomposition
of $V$. Then the sequence
 $L_d, L_{d-1}, \ldots, L_0$ is a decomposition
of $V$, which we call the {\it inversion} of 
 $L_0, L_1, \ldots, L_d$.

\begin{lemma}
\label{lem:yzinvzy}
 Let $\Phi$ denote the Leonard system in
(\ref{eq:ourstartingpt}),  
 let $V$ denote the 
irreducible left ${\cal A}$-module,
and let the set  $\Omega $ be as in Definition
\ref{def:S4interp}. 
Let $y, z$ denote distinct elements in $\Omega $.
Then each of the decompositions $\lbrack yz \rbrack $,
 $\lbrack zy \rbrack $  is the inversion of the 
 other. 

\end{lemma}

\begin{proof} Immediate from
Definition \ref{def:byzb} and the definition
of inversion.

\end{proof}

\noindent Let $\Phi$ denote the Leonard system in
(\ref{eq:ourstartingpt}),
 let $V$ denote the 
irreducible left ${\cal A}$-module,
and let the set  $\Omega $ be as in Definition
\ref{def:S4interp}.
In Definition 
\ref{def:byzb},
we obtained for each ordered pair
$yz$ of distinct elements in   $\Omega $,
a decomposition of $V$
denoted  $\lbrack yz \rbrack $. This gives 12 decompositions of $V$.
By Lemma 
\ref{lem:yzinvzy},
these  consist
of 6 pairs of inverse decompositions.
To be concrete, we now display these decompositions. 

\begin{theorem}
\label{thm:sixdecp}
Let $\Phi$ denote the Leonard
system in 
(\ref{eq:ourstartingpt}), 
let $V$ denote the irreducible left $\cal A$-module, 
and let the set $\Omega $ be as in
Definition
\ref{def:S4interp}.
Let $yz$  denote an ordered pair of distinct
elements in  $\Omega $, 
and consider the corresponding decomposition
$\lbrack yz \rbrack $ of $V$
from
Definition
\ref{def:byzb}.
For $0\leq i \leq d$, the $i^{\hbox{th}}$ subspace of $\lbrack yz \rbrack $
is given in the following table. 
\medskip

\centerline{
\begin{tabular}[t]{c|c}
        $yz$ &$i^{\hbox{th}}$ subspace of decomposition  $\lbrack yz\rbrack
	$ \\ \hline  \hline
        $0^*d$ &
	$(E^*_0V+\cdots + E^*_iV)\cap (E_iV+\cdots +E_dV)$   \\ 
        $d0^*$ &
	$(E^*_0V+\cdots + E^*_{d-i}V)\cap (E_{d-i}V+\cdots +E_dV)$   \\ 
       \hline 
	$0d^*$ &
	$(E_0V+\cdots + E_iV)\cap (E^*_iV+\cdots +E^*_dV)$   \\ 
	$d^*0$ &
	$(E_0V+\cdots + E_{d-i}V)\cap (E^*_{d-i}V+\cdots +E^*_dV)$   \\ 
        \hline 
	$00^*$ & $(E_0V+\cdots +E_iV)\cap 
	(E^*_{d-i}V+\cdots +E^*_0V)$   \\ 
	$0^*0$ & $(E_0V+\cdots +E_{d-i}V)\cap 
	(E^*_{i}V+\cdots +E^*_0V)$   \\ 
        \hline 
	$dd^*$ & $(E_dV+\cdots +E_{d-i}V)\cap 
	(E^*_{i}V+\cdots +E^*_dV)$   \\ 
	$d^*d$ & $(E_dV+\cdots +E_{i}V)\cap 
	(E^*_{d-i}V+\cdots +E^*_dV)$   \\ 
        \hline 
	$0d$ & $E_iV$   \\ 
	$d0$ & $E_{d-i}V$   \\ 
        \hline 
	$0^*d^*$ & $E^*_iV$   \\
        $d^*0^*$ & $E^*_{d-i}V$  
	\end{tabular}}
\medskip
\noindent 

\end{theorem}
\noindent Describing our 12 decompositions from another point of
view, we have the following.

\begin{theorem}
\label{thm:decsum}
Let $\Phi$ denote the Leonard
system in 
(\ref{eq:ourstartingpt}), 
let $V$ denote the irreducible left $\cal A$-module, 
and let the set $\Omega $ be as in
Definition
\ref{def:S4interp}.
Let $yz$  denote an ordered pair of distinct
elements in  $\Omega $, 
and consider the corresponding decomposition
$\lbrack yz \rbrack $ from 
Definition
\ref{def:byzb}.
Let us denote this decomposition by  
$L_0, L_1, \ldots, L_d$. 
 Then for $0 \leq i \leq d$,
the sums $L_0+L_1+\cdots + L_i$ and $L_i+L_{i+1}+\cdots + L_d$
are given as follows.

\medskip

\centerline{
\begin{tabular}[t]{c|c|c}
        $yz$ &$L_0+\cdots + L_i$ & $L_i+\cdots + L_d$ \\ \hline  \hline
        $0^*d$ &
        $E^*_0V+\cdots + E^*_iV$ & $E_iV+\cdots +E_dV$   \\ 
        $d0^*$ &
        $E_dV+\cdots + E_{d-i}V$ &  $E^*_{d-i}V + \cdots + E^*_0V $ \\
        \hline 
       $0d^*$ & 
	$E_0V+\cdots + E_iV $ & $E^*_iV+\cdots +E^*_dV$   \\ 
       $d^*0$ & 
	$E^*_dV+\cdots + E^*_{d-i}V$ 
&	
	$E_{d-i}V + \cdots +E_0V $ 
	\\ 
        \hline 
	$00^*$ & $E_0V+\cdots +E_iV$ & $ 
	E^*_{d-i}V+\cdots +E^*_0V$   \\ 
	$0^*0$ & 
	$E^*_{0}V+\cdots +E^*_iV$  
	&
	$E_{d-i}V+\cdots +E_{0}V$
	\\ 
        \hline
	$dd^*$ & $E_dV+\cdots +E_{d-i}V $&
	$E^*_{i}V+\cdots +E^*_dV$    \\
	$d^*d$ &
	$E^*_{d}V+\cdots +E^*_{d-i}V$ 
&	
	$E_iV+\cdots +E_{d}V $
	\\
        \hline
	$0d$ &
	$E_0V + \cdots + E_iV$  & $E_iV + \cdots + E_dV $ \\ 
	$d0$ &
	$E_{d}V + \cdots +  E_{d-i}V $
 &	
	$E_{d-i}V + \cdots + E_{0}V$
	\\ 
        \hline
	$0^*d^*$ &
	$E^*_0V+\cdots + E^*_iV$  & $ E^*_iV + \cdots + E^*_dV$   \\
	$d^*0^*$ &
	$ E^*_{d}V + \cdots + E^*_{d-i}V$  
 &	
	$E^*_{d-i}V+\cdots + E^*_{0}V$ 
	\end{tabular}}
\medskip

\end{theorem}

\section{$24$ bases for $V$}

\medskip
\noindent Let $\Phi$ denote the Leonard system in
(\ref{eq:ourstartingpt}),
and let $V$ denote the 
irreducible left ${\cal A}$-module.
For each element $g \in S_4$, we display 
a basis for $V$, denoted $\lbrack g \rbrack $.
To describe our procedure, we use the following notation.

\medskip
\noindent 
Let $u_0, u_1, \ldots, u_d$ denote a basis for $V$, and
set $L_i = \hbox{Span}(u_i)$ for $0 \leq i\leq d$.
Observe the sequence $L_0, L_1, \ldots, L_d$ is a decomposition
of $V$. We say this decomposition is {\it induced} by
 $u_0, u_1, \ldots, u_d$.

\medskip
\noindent 
Let the set  $\Omega$ be as in Definition
\ref{def:S4interp},
and let $yz$ denote an ordered pair of distinct elements
of $\Omega$.
Consider the corresponding
 decomposition of $V$, denoted $\lbrack yz \rbrack$.
 We define two bases for $V$, both of  which induce 
 $\lbrack yz \rbrack $.
We denote  these bases
by  $\lbrack wxyz \rbrack $ and $\lbrack xwyz \rbrack$,
where $w$ and $x$ denote
the elements in  $\Omega $
 other than $y,z$.
Apparently this procedure yields, for each $g\in S_4$, 
a basis 
$\lbrack g \rbrack $ for $V$.
These 24 bases are displayed below.
\begin{theorem}
\label{thm:bases}
Let $\Phi$ denote the Leonard
system in 
(\ref{eq:ourstartingpt}), 
and let $V$ denote the irreducible left $\cal A$-module.
Let $\eta_0$, $\eta_d$, $\eta^*_0$, $\eta^*_d$ denote nonzero
vectors in $V$ such that
\begin{equation}
\eta_0 \in E_0V,\qquad 
\eta_d \in E_dV,\qquad 
\eta^*_0 \in E^*_0V,\qquad 
\eta^*_d \in E^*_dV.
\label{eq:videf}
\end{equation}
With reference to Definition
\ref{def:S4interp},
let $g$ denote an element of $S_4$ and consider
row $g$ of the table below.
For $0 \leq i \leq d$,  the vector  $v_i$ given in
that row is a basis for 
the subspace given to its right. Moreover,  the sequence
$v_0, v_1, \ldots, v_d$ is
a basis for $V$. We denote this basis by $\lbrack g \rbrack $.  

\medskip

\centerline{
\begin{tabular}[t]{c|c|c}
        $g$ &$v_i$ & $v_i$ is basis for \\ \hline  \hline
        $d^*00^*d$ & 
	$(A-\theta_0)\cdots (A-\theta_{i-1})\eta^*_0$ &
	$(E^*_0V + \cdots + E^*_iV)\cap (E_iV + \cdots + E_dV) $
	\\
        $0d^*0^*d$ & 
	$(A^*-\theta^*_d)\cdots (A^*-\theta^*_{i+1})\eta_d$ & 
	$(E^*_0V + \cdots + E^*_iV)\cap (E_iV + \cdots + E_dV) $
	\\
        $d^*0d0^*$ &
	$(A-\theta_0)\cdots (A-\theta_{d-i-1})\eta^*_0$ &
	$(E^*_0V + \cdots + E^*_{d-i}V)\cap (E_{d-i}V + \cdots + E_dV) $
	\\
        $0d^*d0^*$ & 
	$(A^*-\theta^*_d)\cdots (A^*-\theta^*_{d-i+1})\eta_d$ &
	$(E^*_0V + \cdots + E^*_{d-i}V)\cap (E_{d-i}V + \cdots + E_dV) $
        \\
  \hline	
        $d0^*0d^*$ & 
	$(A^*-\theta^*_0)\cdots (A^*-\theta^*_{i-1})\eta_0$ &
	$(E_0V + \cdots + E_iV)\cap (E^*_iV + \cdots + E^*_dV) $
	\\
        $0^*d0d^*$ &
	$(A-\theta_d)\cdots (A-\theta_{i+1})\eta^*_d$ &
	$(E_0V + \cdots + E_iV)\cap (E^*_iV + \cdots + E^*_dV) $
	\\
        $d0^*d^*0$ &
	$(A^*-\theta^*_0)\cdots (A^*-\theta^*_{d-i-1})\eta_0$ & 
	$(E_0V + \cdots + E_{d-i}V)\cap (E^*_{d-i}V + \cdots + E^*_dV) $
	\\
        $0^*dd^*0$ & 
	$(A-\theta_d)\cdots (A -\theta_{d-i+1})\eta^*_d$ &
	$(E_0V + \cdots + E_{d-i}V)\cap (E^*_{d-i}V + \cdots + E^*_dV) $
        \\
  \hline	
        $dd^*00^*$ &
	$(A^*-\theta^*_d)\cdots (A^*-\theta^*_{d-i+1})\eta_0$ & 
	$(E_0V + \cdots + E_iV)\cap (E^*_{d-i}V + \cdots + E^*_0V) $ 
	\\
	$d^*d00^*$ & 
	$(A-\theta_d)\cdots (A-\theta_{i+1})\eta^*_0$ &
	$(E_0V + \cdots + E_iV)\cap (E^*_{d-i}V + \cdots + E^*_0V) $
	\\
        $dd^*0^*0$ & 
	$(A^*-\theta^*_d)\cdots (A^*-\theta^*_{i+1})\eta_0$ &
	$(E_0V + \cdots + E_{d-i}V)\cap (E^*_iV + \cdots + E^*_0V) $
	\\
	$d^*d0^*0$ &
	$(A-\theta_d)\cdots (A-\theta_{d-i+1})\eta^*_0$ & 
	$(E_0V + \cdots + E_{d-i}V)\cap (E^*_iV + \cdots + E^*_0V) $ 
	\\
  \hline	
        $00^*dd^*$ &
	$(A^*-\theta^*_0)\cdots (A^*-\theta^*_{i-1})\eta_d$ &
	$(E_dV + \cdots + E_{d-i}V)\cap (E^*_iV + \cdots + E^*_dV) $
	\\
	$0^*0dd^*$ &
	$(A-\theta_0)\cdots (A-\theta_{d-i-1})\eta^*_d$ & 
	$(E_dV + \cdots + E_{d-i}V)\cap (E^*_iV + \cdots + E^*_dV) $ 
	\\
        $00^*d^*d$ & 
	$(A^*-\theta^*_0)\cdots (A^*-\theta^*_{d-i-1})\eta_d$ &
	$(E_dV + \cdots + E_{i}V)\cap (E^*_{d-i}V + \cdots + E^*_dV) $ 
	\\
	$0^*0d^*d$ & 
	$(A-\theta_0)\cdots (A-\theta_{i-1})\eta^*_d$ & 
	$(E_dV + \cdots + E_{i}V)\cap (E^*_{d-i}V + \cdots + E^*_dV) $
	\\
       \hline 
	$d^* 0^* 0d$ & $E_i\eta^*_0$  & $E_iV$ \\ 
        $0^* d^* 0d$ & $E_i\eta^*_d$  & $E_iV$ \\ 
        $d^* 0^* d0$ & $E_{d-i}\eta^*_0$  & $E_{d-i}V$ \\ 
        $0^* d^* d0$ & $E_{d-i}\eta^*_d$  & $E_{d-i}V$ \\ 
       \hline 
        $d00^*d^*$ & $E^*_i\eta_0$  & $E^*_iV$ \\ 
        $0d0^*d^*$ & $E^*_i\eta_d$  & $E^*_iV$ \\ 
        $d0d^*0^*$ & $E^*_{d-i}\eta_0$  & $E^*_{d-i}V$ \\ 
        $0dd^*0^*$ & $E^*_{d-i}\eta_d$  & $E^*_{d-i}V$  
\end{tabular}}
\medskip

\end{theorem}

\begin{proof} Concerning the first row of the above table,
our assertions follow from the lines preceeding
(\ref{eq:basis1}). Concerning the
 third row of the above table, our assertions follow
 upon replacing $i$
by $d-i$ in the first row.
We have now proved our assertions for
the first and third rows of the table. Applying
these assertions to the relatives of $\Phi$, we obtain the first 16 
rows of the table. 
Consider the next remaining row, where $g$ equals
$d^*0^*0d$. For this row, our assertions
are immediate from
Lemma \ref{lem:stbasis}.
Applying this result to the relatives of $\Phi$,
we obtain the remaining rows of the table.

\end{proof}
\noindent
We record a few observations.

\begin{lemma}
\label{lem:invers}
Referring to Theorem
\ref{thm:bases}, for all elements $wxyz$ in $S_4$, the basis
$\lbrack wxyz \rbrack $ 
from
Theorem
\ref{thm:bases}
induces the decomposition
$\lbrack yz \rbrack$ of $V$ from Definition
\ref{def:byzb}.

\end{lemma}
\begin{proof} Compare the data in 
Theorem \ref{thm:bases} with the data in
Theorem
\ref{thm:sixdecp}.

\end{proof}

\begin{lemma}
\label{lem:reconcile}
Let $\Phi$ denote the Leonard system
in
(\ref{eq:ourstartingpt}), and let $V$ denote the irreducible left
$\cal A$-module. In the table below,
each basis for $V$ contained
the first column (resp. second column) (resp. third
column) is a $\Psi$-standard basis (resp. $\Psi$-split basis)
(resp. $\Psi$-inverted split
basis), where $\Psi$ is the 
relative of $\Phi$ given to the left of this basis.

\bigskip
\centerline{
\begin{tabular}[t]{c|ccc}
$\Psi$ & $\Psi$-standard basis  & $\Psi$-split basis & $\Psi$-inv. split basis  
\\
\hline  \hline
\\
$\Phi$ 
 &
$ \lbrack d^*0^*0d \rbrack $
 &
$\lbrack d^*00^*d \rbrack $
&
$\lbrack d^*0d0^* \rbrack $
\\
\\
$\Phi^\downarrow $
&
$ \lbrack 0^*d^*0d \rbrack  $
 &
$\lbrack 0^*0d^*d \rbrack $
&
$\lbrack 0^*0dd^* \rbrack $
\\
\\
$\Phi^\Downarrow $
&
$ \lbrack d^*0^*d0 \rbrack  $
 &
$\lbrack d^*d0^*0 \rbrack $
&
$\lbrack d^*d00^* \rbrack $
\\
\\
$\Phi^{\downarrow \Downarrow} $
&
$ \lbrack 0^*d^*d0 \rbrack  $
 &
$\lbrack 0^*dd^*0 \rbrack $
&
$\lbrack 0^*d0d^* \rbrack $
\\
\\
\hline
\\
$\Phi^*$ 
 &
$ \lbrack d00^*d^* \rbrack $
 &
$\lbrack d0^*0d^* \rbrack $
&
$\lbrack d0^*d^*0 \rbrack $
\\
\\
$\Phi^{\downarrow *} $
&
$ \lbrack d0d^*0^* \rbrack  $
 &
$\lbrack dd^*00^* \rbrack $
&
$\lbrack dd^*0^*0 \rbrack $
\\
\\
$\Phi^{\Downarrow *}$
&
$ \lbrack 0d0^*d^* \rbrack  $
 &
$\lbrack 00^*dd^* \rbrack $
&
$\lbrack 00^*d^*d \rbrack $
\\
\\
$\Phi^{\downarrow \Downarrow *} $
&
$ \lbrack 0dd^*0^* \rbrack  $
 &
$\lbrack 0d^*d0^* \rbrack $
&
$\lbrack 0d^*0^*d \rbrack $
\end{tabular}}
\medskip
\end{lemma}

\begin{proof} Immediate from inspecting the table in
Theorem
\ref{thm:bases}.

\end{proof}

\noindent Later in the paper, we will compute, for each ordered pair
$g, h$ of adjacent elements in $S_4$, the entries in the transition
matrix from the basis $\lbrack g \rbrack $ to the basis $\lbrack h \rbrack$.
Before going that far,
we  say something
about the general nature of these transition matrices. First 
we recall  our terms.

\medskip
\noindent
Let $\Phi$ denote the Leonard system in 
(\ref{eq:ourstartingpt}), and let $V$ denote the irreducible
left $\cal A$-module.
Suppose we are given two bases for $V$, written
$u_0, u_1, \ldots, u_d$ and
$v_0, v_1, \ldots, v_d$. By the {\it transition matrix} 
from 
$u_0, u_1, \ldots, u_d$ to 
$v_0, v_1, \ldots, v_d$, we mean the matrix 
$T$ in $\hbox{Mat}_{d+1}(\K)$ satisfying
\begin{equation}
v_j = \sum_{i=0}^d T_{ij}u_i      \qquad \qquad (0 \leq j\leq d).
\label{eq:transdefpre}
\end{equation}
We recall a few properties of transition matrices.
Let $T$ denote the transition matrix from
$u_0, u_1, \ldots, u_d$ to 
$v_0, v_1, \ldots, v_d$. Then $T^{-1}$ exists, and equals
the transition matrix from
$v_0, v_1, \ldots, v_d$ to
$u_0, u_1, \ldots, u_d$.  
Let 
$w_0, w_1, \ldots, w_d$ denote a basis for $V$, and let
$S$ denote the transition matrix from 
$v_0, v_1, \ldots, v_d$ to
$w_0, w_1, \ldots, w_d$. Then $TS$ is the transition 
matrix from 
$u_0, u_1, \ldots, u_d$  to 
$w_0, w_1, \ldots, w_d$.  

\begin{lemma}
\label{lem:transshape}
Let $\Phi$ denote the Leonard system in 
(\ref{eq:ourstartingpt}), and let $V$ denote the irreducible
left $\cal A$-module. With reference to 
Definition
\ref{def:S4interp}, 
let  $g,h$ denote adjacent elements in $S_4$,
and consider the corresponding bases $\lbrack g \rbrack$,
$\lbrack h \rbrack $ for $V$ given in Theorem
\ref{thm:bases}.
Then (i)--(iii) hold
below.
\begin{enumerate}
\item Suppose $g,h$ are $1$-adjacent. Then the transition matrix
from $\lbrack g\rbrack $ to $\lbrack h\rbrack $ is diagonal.
\item Suppose $g,h$ are $2$-adjacent. Then the transition matrix
from $\lbrack g\rbrack $ to $\lbrack h\rbrack $ is lower triangular.
\item Suppose $g,h$ are $3$-adjacent.  Then 
 $\lbrack g\rbrack $ is the inversion of $\lbrack h\rbrack $.
\end{enumerate}

\end{lemma}

\begin{proof} For notational convenience we write $g=wxyz$. \hfil\break
\noindent 
(i) In this case  $h=xwyz$. Observe $\lbrack g \rbrack $
and $\lbrack h \rbrack $
both induce the decomposition $\lbrack yz\rbrack$
by Lemma 
\ref{lem:invers},
so 
the transition matrix 
from $\lbrack g \rbrack $ to $\lbrack h \rbrack $ is diagonal.

\noindent (ii) In this case $h=wyxz$.  
By Lemma 
\ref{lem:invers},
the bases $\lbrack g \rbrack$
and $\lbrack h \rbrack $ induce
 the decompositions $\lbrack yz\rbrack $ and 
 $\lbrack xz \rbrack $, respectively.
When we consider how the decompositions 
  $\lbrack yz\rbrack $ and $\lbrack xz \rbrack $ are related,
  we find 
the transition matrix 
from $\lbrack g \rbrack $ to $\lbrack h \rbrack $ is lower triangular.

\noindent (iii) 
In this case $h=wxzy$.
 In the table of Theorem
\ref{thm:bases},
for each block we compare rows 1,3 and rows 2,4. We find in all
cases $\lbrack g \rbrack $ is the inversion of $\lbrack h \rbrack $.

\end{proof}

\section{Some scalars}

\noindent Our next goal is to compute the matrices 
representing $A$ and $A^*$ with respect to each of the
bases in Theorem 
\ref{thm:bases}. To describe the entries of these matrices,
we will use the following parameters.

\begin{definition}
\label{def:aidef} Let $\Phi$ denote the Leonard system
in
(\ref{eq:ourstartingpt}).
We define
\begin{eqnarray}
a_i = \hbox{tr}\, AE^*_i, \qquad  \quad
a^*_i = \hbox{tr}\, A^*E_i, \qquad \qquad  
 (0 \leq i \leq d), \qquad
\label{eq:defofaiS99}
\end{eqnarray}
where $tr$ means trace.
\end{definition}

\medskip
\noindent The scalars $a_i, a^*_i$ have the following
interpretation.

\begin{lemma}
\label{lem:sumae}
With reference to  Definition
\ref{def:aidef},
\begin{eqnarray}
E^*_iAE^*_i &=& a_i E^*_i \qquad \qquad (0 \leq i \leq d),
\label{eq:aimeaning}
\\
E_iA^*E_i &=& a^*_i E_i \qquad \qquad (0 \leq i \leq d).
\label{eq:aismeaning}
\end{eqnarray}

\end{lemma}

\begin{proof} Concerning 
(\ref{eq:aimeaning}), let $i$ be given. Since $E^*_i$ is a rank 1
idempotent, there exists a scalar $\alpha_i\in \K$ such that
\begin{equation}
E^*_iAE^*_i = \alpha_i E^*_i.
\label{eq:findalph}
\end{equation}
Taking the trace of both sides of (\ref{eq:findalph}),  and
recalling 
$ XY$,  $YX$ have the same trace, we routinely find $\alpha_i = a_i$.
We have now proved
(\ref{eq:aimeaning}).
Applying this to $\Phi^*$, we obtain
(\ref{eq:aismeaning}).

\end{proof}

\begin{lemma}
\label{thm:aiform} 
Let $\Phi$ denote the Leonard system
in
(\ref{eq:ourstartingpt}).
Then for $0 \leq i \leq d$ the scalar $a_i$ equals
both
\begin{eqnarray}
 \theta_i + {{\varphi_i}\over {\theta^*_i-\theta^*_{i-1}}}
+ {{\varphi_{i+1}}\over {\theta^*_i-\theta^*_{i+1}}},
\qquad \qquad 
 \theta_{d-i} + {{\phi_i}\over {\theta^*_i-\theta^*_{i-1}}}
+ {{\phi_{i+1}}\over {\theta^*_i-\theta^*_{i+1}}},
\label{eq:aivarphi}
\end{eqnarray}
 where  $\theta^*_{-1}$, $\theta^*_{d+1}$ denote 
 indeterminants.
Moreover, the scalar $a^*_i$ equals both
\begin{eqnarray}
 \theta^*_i + {{\varphi_i}\over {\theta_i-\theta_{i-1}}}
+ {{\varphi_{i+1}}\over {\theta_i-\theta_{i+1}}},
\qquad \qquad 
 \theta^*_{d-i} + {{\phi_{d-i+1}}\over {\theta_i-\theta_{i-1}}}
+ {{\phi_{d-i}}\over {\theta_i-\theta_{i+1}}},
\label{eq:aisvarphi} 
\end{eqnarray}
where
 $\theta_{-1}$, $\theta_{d+1}$ denote 
 indeterminants.
\end{lemma}

\begin{proof} Let the integer $i$ be given.
The scalar $a_i$ equals the expression on
the left in 
(\ref{eq:aivarphi}) by 
\cite[Lem. 5.1]{LS99}. Applying this fact to $\Phi^\Downarrow$,
and using Theorem
\ref{thm:phimod}(iii), 
we find $a_i$ equals
the expression on the right in 
(\ref{eq:aivarphi}). 
We have now shown $a_i$ equals the two expressions in 
(\ref{eq:aivarphi}). 
Applying this to $\Phi^*$, and using
Theorem
\ref{thm:phimod}(i), we find
 $a^*_i$ equals the two expressions in
(\ref{eq:aisvarphi}). 

\end{proof}

\section{The $24$ bases; matrices representing $A$ and $A^*$}

\medskip
\noindent 
In this section, we return to the 24 bases
in Theorem
\ref{thm:bases}.
For each $g \in S_4$, we  
compute the matrices representing  $A$ and $A^*$ with respect
to the basis $\lbrack g \rbrack $.

\medskip
\noindent We use the following notation.

\begin{definition}
\label{def:matrixrep} Let $\Phi$ denote the Leonard system 
in 
(\ref{eq:ourstartingpt}), and let $V$ denote the irreducible
left $\cal A$-module. With reference to Definition
\ref{def:S4interp}, 
let $g$ denote an element in $S_4$.
For all $X \in \cal A $,
we let $X^g$ denote the matrix in 
$\hbox{Mat}_{d+1}(\K)$ that represents $X$ with respect to
the basis $\lbrack g \rbrack $, where  
$\lbrack g \rbrack $ is from
Theorem
\ref{thm:bases}.
Denoting this basis by $v_0, v_1, \ldots, v_d$ we have
\beast
Xv_j = \sum_{i=0}^d X^g_{ij} v_i \qquad \qquad (0 \leq j \leq d).
\eeast
We observe the map $X \rightarrow X^g$
is a $\K$-algebra isomorphism from $\cal A $ to 
$\hbox{Mat}_{d+1}(\K)$.

\end{definition}

\begin{theorem}
\label{thm:repa}
Let $g$ denote an element of $S_4$.
With reference to Definition 
\ref{def:matrixrep},
the entries of $A^g$ and $A^{*g}$ 
are given in the tables below. Any entry not displayed  
is zero. 

\bigskip

\centerline{
\begin{tabular}[t]{c|ccc|ccc}
g & $A^{g}_{i,i-1}$ & $A^{g}_{ii}$  & $A^{g}_{i-1,i}$   &
 $A^{*g}_{i,i-1}$ & $A^{*g}_{ii}$  & $A^{*g}_{i-1,i}$   
\\
\hline  \hline
 $d^*00^*d$ 
 &
1
&
$\theta_i $
& 
0
 &
 0 
 &
$\theta^*_i  $
 &
 $\varphi_i $ 
\\
 $0d^*0^*d$ 
 &
$\varphi_i $
&
$\theta_{i} $
& 
0
 &
 0
 &
$\theta^*_{i}$ 
 &
1
\\
$d^*0d0^*$ 
 &
0
&
$\theta_{d-i} $
& 
1
 &
$\varphi_{d-i+1} $
 &
$\theta^*_{d-i}  $
 &
0
\\
 $0d^*d0^*$ 
 &
0
&
$\theta_{d-i} $
& 
$\varphi_{d-i+1} $
 &
 1
 &
$\theta^*_{d-i}  $
 &
0
\\
\hline
 $d0^*0d^*$ 
 &
0
&
$\theta_{i} $
& 
$\varphi_{i} $
 &
 1
 &
$\theta^*_{i} $ 
 &
0
\\
$0^*d0d^*$ 
 &
0
&
$\theta_{i} $
& 
1 
 &
$\varphi_{i} $
 &
$\theta^*_{i} $ 
 &
0
\\
$d0^*d^*0$ 
 &
$\varphi_{d-i+1} $
&
$\theta_{d-i} $
& 
0
 &
 0
 &
$\theta^*_{d-i}  $
 &
1
\\
 $0^*dd^*0$ 
 &
1
&
$\theta_{d-i} $
& 
0
 &
 0
 &
$\theta^*_{d-i} $ 
 &
$\varphi_{d-i+1} $
\\
\hline
 $dd^*00^*$ 
 &
0
&
$\theta_{i} $
& 
$\phi_{d-i+1} $
 &
1 
 &
$\theta^*_{d-i} $ 
 &
0
\\
$d^*d00^*$ 
 &
0
&
$\theta_{i} $
& 
1 
 &
$\phi_{d-i+1} $
 &
$\theta^*_{d-i} $ 
 &
0
\\
$dd^*0^*0$ 
 &
$\phi_i $
&
$\theta_{d-i} $
& 
0 
 &
0 
 &
$\theta^*_{i}  $
 &
1
\\
$d^*d0^*0$ 
 &
1
&
$\theta_{d-i} $
& 
0 
 &
0 
 &
$\theta^*_{i} $ 
 &
$\phi_i $
\\
\hline
 $00^*dd^*$ 
 &
0
&
$\theta_{d-i} $
& 
$\phi_{i} $
 &
1 
 &
$\theta^*_{i}  $
 &
0
\\
$0^*0dd^*$ 
 &
0
&
$\theta_{d-i} $
& 
1 
 &
$\phi_{i} $
 &
$\theta^*_{i}  $
 &
0
\\
$00^*d^*d$ 
 &
$\phi_{d-i+1} $
&
$\theta_{i} $
& 
0 
 &
0 
 &
$\theta^*_{d-i} $
 &
1
\\
$0^*0d^*d$ 
 &
1
&
$\theta_{i} $
& 
0 
 &
0 
 &
$\theta^*_{d-i}  $
 &
$\phi_{d-i+1} $
\end{tabular}}
\medskip

\centerline{
\begin{tabular}[t]{c|c|ccc}
g & $A^g_{ii}$ & $A^{*g}_{i,i-1}$ & $A^{*g}_{ii}$  & $A^{*g}_{i-1,i}$  
\\
\hline
&&&&
\\
 $d^* 0^* 0d$ & $\theta_i $  &
 $\phi_{d-i+1} \frac{(\theta_i -\theta_d)\cdots (\theta_i - \theta_{i+1})}{(\theta_{i-1}-\theta_d)\cdots (\theta_{i-1}-\theta_i)}$
 &
 $a^*_i $ & $\varphi_i \frac{(\theta_{i-1}-\theta_0)\cdots (\theta_{i-1}-\theta_{i-2})}{(\theta_{i}-\theta_0)\cdots (\theta_{i}-\theta_{i-1})}$ 
\\
 &&&&
\\ 
 $0^*d^* 0d$ & $\theta_{i} $  &
 $\varphi_{i} \frac{(\theta_{i} -\theta_d)\cdots (\theta_{i} - \theta_{i+1})}{(\theta_{i-1}-\theta_d)\cdots (\theta_{i-1}-\theta_i)}$
 &
 $a^*_i $ 
 &
 $\phi_{d-i+1} \frac{(\theta_{i-1}-\theta_0)\cdots (\theta_{i-1}-\theta_{i-2})}{(\theta_{i}-\theta_0)\cdots (\theta_{i}-\theta_{i-1})}$ 
\\
 &&&&
\\
 $d^* 0^* d0$ & $\theta_{d-i} $  &
 $\varphi_{d-i+1} \frac{(\theta_{d-i} -\theta_0)\cdots (\theta_{d-i} - \theta_{d-i-1})}{(\theta_{d-i+1}-\theta_0)\cdots (\theta_{d-i+1}-\theta_{d-i})}$
 &
 $a^*_{d-i} $ 
 &
 $\phi_i \frac{(\theta_{d-i+1}-\theta_d)\cdots (\theta_{d-i+1}-\theta_{d-i+2})}{(\theta_{d-i}-\theta_d)\cdots (\theta_{d-i}-\theta_{d-i+1})}$ 
\\
&&&&
\\ 
 $0^*d^* d0$ & $\theta_{d-i} $  &
 $\phi_{i} \frac{(\theta_{d-i} -\theta_0)\cdots (\theta_{d-i} - \theta_{d-i-1})}{(\theta_{d-i+1}-\theta_0)\cdots (\theta_{d-i+1}-\theta_{d-i})}$
 &
 $a^*_{d-i} $ 
 &
 $\varphi_{d-i+1} \frac{(\theta_{d-i+1}-\theta_d)\cdots (\theta_{d-i+1}-\theta_{d-i+2})}{(\theta_{d-i}-\theta_d)\cdots (\theta_{d-i}-\theta_{d-i+1})}$ 
\end{tabular}}

\bigskip

\centerline{
\begin{tabular}[t]{c|ccc|c}
g & $A^{g}_{i,i-1}$ & $A^{g}_{ii}$  & $A^{g}_{i-1,i}$   &
$A^{*g}_{ii}$
\\
\hline 
 &&&&
 \\
$d00^*d^*$ 
 &
 $\phi_{i} \frac{(\theta^*_i -\theta^*_d)\cdots (\theta^*_i - \theta^*_{i+1})}{(\theta^*_{i-1}-\theta^*_d)\cdots (\theta^*_{i-1}-\theta^*_i)}$
 &
 $a_i $ 
 &
 $\varphi_i \frac{(\theta^*_{i-1}-\theta^*_0)\cdots (\theta^*_{i-1}-\theta^*_{i-2})}{(\theta^*_{i}-\theta^*_0)\cdots (\theta^*_{i}-\theta^*_{i-1})}$ 
 &
 $\theta^*_i $ 
\\
 &&&&
\\ 
 $0d 0^*d^*$ 
 &
 $\varphi_{i} \frac{(\theta^*_{i} -\theta^*_d)\cdots (\theta^*_{i} - \theta^*_{i+1})}{(\theta^*_{i-1}-\theta^*_d)\cdots (\theta^*_{i-1}-\theta^*_i)}$
 &
 $a_i $ 
 &
 $\phi_{i} \frac{(\theta^*_{i-1}-\theta^*_0)\cdots (\theta^*_{i-1}-\theta^*_{i-2})}{(\theta^*_{i}-\theta^*_0)\cdots (\theta^*_{i}-\theta^*_{i-1})}$ 
 &
 $\theta^*_{i} $ 
\\
&&&&
\\
 $d 0 d^*0^*$
 &
 $\varphi_{d-i+1} \frac{(\theta^*_{d-i} -\theta^*_0)\cdots (\theta^*_{d-i} - \theta^*_{d-i-1})}{(\theta^*_{d-i+1}-\theta^*_0)\cdots (\theta^*_{d-i+1}-\theta^*_{d-i})}$
 &
 $a_{d-i} $ 
 &
 $\phi_{d-i+1} \frac{(\theta^*_{d-i+1}-\theta^*_d)\cdots (\theta^*_{d-i+1}-\theta^*_{d-i+2})}{(\theta^*_{d-i}-\theta^*_d)\cdots (\theta^*_{d-i}-\theta^*_{d-i+1})}$ 
& 
 $\theta^*_{d-i} $ 
\\
&&&&
\\ 
 $0d d^*0^*$
 &
 $\phi_{d-i+1} \frac{(\theta^*_{d-i} -\theta^*_0)\cdots 
 (\theta^*_{d-i} - \theta^*_{d-i-1})}{(\theta^*_{d-i+1}-\theta^*_0)\cdots (\theta^*_{d-i+1}-\theta^*_{d-i})}$
 &
 $a_{d-i} $ 
 &
 $\varphi_{d-i+1} \frac{(\theta^*_{d-i+1}-\theta^*_d)\cdots (\theta^*_{d-i+1}-\theta^*_{d-i+2})}{(\theta^*_{d-i}-\theta^*_d)\cdots (\theta^*_{d-i}-\theta^*_{d-i+1})}$ 
 &
 $\theta^*_{d-i} $ 
\end{tabular}}
\medskip

\end{theorem}

\begin{proof} Consider  the first row of the first table,
where $g$ equals $d^*00^*d$. As indicated in 
the table of  
Lemma
\ref{lem:reconcile}, row 1, column 2,
the basis  
 $\lbrack g \rbrack $ 
is a $\Phi$-split basis. 
From the line above  
(\ref{eq:matrepaastar}), we find
 $A^g$ (resp. $A^{*g}$) is 
given on the left (resp. right) in
(\ref{eq:matrepaastar}).  From this we obtain our
results for the first row of the first table.
Now consider the third row of the first table, where
$g$ equals $d^*0d0^*$.
From the table of 
Lemma
\ref{lem:reconcile}, row 1, column 3, 
the basis $\lbrack d^*0d0^* \rbrack $
is a $\Phi$-inverted split basis.
From the line above
(\ref{eq:invmat}) 
 we find $A^g$ (resp. $A^{*g}$) is 
given on the left (resp. right) in
(\ref{eq:invmat}).
From this we obtain our
results for the third row of the first table.
We have now proved our assertions for rows 1 and 3
of the first table.
Applying this result to the relatives of $\Phi$,
and using Theorem
\ref{thm:phimod},
we obtain the remaining rows of the first table.
Consider the first row of the second table, where 
$g$ equals $d^*0^*0d$.
From the table of Theorem
\ref{thm:bases}, row 17,
we find the corresponding basis $\lbrack g \rbrack $ is 
\begin{equation}
E_0\eta^*_0, E_1\eta^*_0,\ldots, E_d\eta^*_0.
\label{eq:stbasis}
\end{equation}
For $0 \leq i \leq d$, the vector $E_i\eta^*_0 $ is an eigenvector
for $A$, with eigenvalue $\theta_i$. Therefore
\begin{equation}
A^g = \hbox{diag} (\theta_0,\theta_1,\ldots, \theta_d).
\label{eq:asrep}
\end{equation}
We now find $A^{*g}$. From the construction,
and since $A,A^*$ is a Leonard pair, the matrix $A^{*g}$ is irreducible
tridiagonal. From 
(\ref{eq:aismeaning})
we find the diagonal entries 
$A^{*g}_{ii} = a^*_i$ for  $0 \leq i  \leq d $.
We show
\begin{eqnarray}
A^{*g}_{i-1,i} &=& \varphi_{i}
\frac{(\theta_{i-1}-\theta_0)(\theta_{i-1}-\theta_1)\cdots (\theta_{i-1}-\theta_{i-2})}
{(\theta_{i}-\theta_0)(\theta_{i}-\theta_1)\cdots
(\theta_{i}-\theta_{i-1})}
\label{eq:biform}
\end{eqnarray}
for $1 \leq i \leq d$.
To see 
(\ref{eq:biform}),
we momentarily return to the basis 
$\lbrack d^*00^*d\rbrack $.
From the table of 
Theorem \ref{thm:bases}, row 1,
we find that 
 for $0 \leq j \leq d$, the $j^{\hbox{th}}$ vector
in the basis
$\lbrack d^*00^*d\rbrack $
is  given by
\begin{equation}
(A-\theta_0 I)(A- \theta_1 I) \cdots (A-\theta_{j-1} I)\eta^*_0.
\label{eq:ith}
\end{equation}
We write 
(\ref{eq:ith}) in terms of 
(\ref{eq:stbasis}).
Recall the sum $E_0+E_1+\cdots + E_d$ equals
the identity $I$. Applying this sum to the vector  
(\ref{eq:ith}) and simplifying the result using
(\ref{eq:primid1S99}),
we find
the vector (\ref{eq:ith})
equals
\begin{equation}
\sum_{i=0}^d (\theta_i-\theta_0)(\theta_i-\theta_1)\cdots
(\theta_i-\theta_{j-1})E_i\eta^*_0.
\label{eq:basech}
\end{equation}
Let $L$ denote the matrix in $\hbox{Mat}_{d+1}(\K)$
with $ij^{\hbox{th}}$ entry 
$(\theta_i-\theta_0)\cdots
(\theta_i-\theta_{j-1})$, for $0 \leq i,j\leq d$. 
Apparently $L$
is the transition matrix from the basis 
$\lbrack d^*0^*0d \rbrack $
to the  basis 
$\lbrack d^*00^*d\rbrack $. By linear algebra,
we obtain 
\begin{equation}
A^{*g}L=LA^{*h},
\label{eq:atta}
\end{equation}
where we recall $g=d^*0^*0d$ and we abbreviate $h=d^*00^*d$.
For $1 \leq i \leq d$, we compute the $i-1,i$ entry
in
(\ref{eq:atta}).
Since 
$A^{*g}$ is tridiagonal,  and since $L$ is lower triangular,
we find the  $i-1,i$ entry of $A^{*g}L$ equals
$A^{*g}_{i-1,i}L_{ii}$ or in other words  
\begin{equation}
A^{*g}_{i-1,i} 
(\theta_{i}-\theta_0)(\theta_{i}-\theta_1)\cdots
(\theta_{i}-\theta_{i-1}).
\label{eq:bside1}
\end{equation}
We mentioned above the matrix 
 $A^{*h}$ is given on the right in
(\ref{eq:matrepaastar}).
Since $A^{*h}$ is upper bidiagonal, 
and since $L$ is lower triangular,
we find the  $i-1,i$ entry of $LA^{*h}$ equals
$L_{i-1,i-1}A^{*h}_{i-1,i}$ or in other words 
\begin{equation}
(\theta_{i-1}-\theta_0)(\theta_{i-1}-\theta_1)\cdots
(\theta_{i-1}-\theta_{i-2})
\varphi_{i}. 
\label{eq:bside2}
\end{equation}
Equating 
(\ref{eq:bside1}),  
(\ref{eq:bside2}), we obtain
(\ref{eq:biform}).
Applying 
(\ref{eq:biform})
to $\Phi^\Downarrow $
 and using
Theorem
\ref{thm:phimod},
we routinely find 
\beast
A^{*g}_{i,i-1} = \phi_{d-i+1} 
\frac{(\theta_i-\theta_d)(\theta_i-\theta_{d-1})\cdots (\theta_i-\theta_{i+1})}
{(\theta_{i-1}-\theta_d)(\theta_{i-1}-\theta_{d-1})\cdots
(\theta_{i-1}-\theta_i)}
\eeast
for $1\leq i \leq d$.
We have now proved our assertions for 
the first row of the second table.
Applying these facts to the relatives of $\Phi$, and using
Theorem
\ref{thm:phimod}, we obtain the remaining 
rows of the second table and all rows of the third table.

\end{proof}

\noindent Summarizing the data from Theorem
\ref{thm:repa},
we have the following.

\begin{lemma}
\label{lem:shapeofa}
Referring to Theorem
\ref{thm:repa},  pick any $g\in S_4$, 
and 
consider the  form of
$A^g$ and $A^{*g}$.
Writing $g=wxyz$, this
form is given as follows.

\bigskip

\centerline{
\begin{tabular}[t]{c|c|c|c}
$y \in \lbrace  0^*,d^*\rbrace $  
& $z \in \lbrace  0^*,d^*\rbrace $  
& $A^g$  
 & $A^{*g}$ 
\\
\hline  
$ No $ & $No$  & diagonal & irred. tridiagonal 
\\
$Yes$ & $No$ & lower bidiagonal & upper bidiagonal  
\\
$No$ & $Yes$  & upper bidiagonal & lower bidiagonal  
\\
$Yes $  & $Yes$  &irred. tridiagonal &  diagonal 
\end{tabular}}

\bigskip
\noindent
We remark the number of elements in $S_4$ satisfying
each of the above four cases is  $4, 8,8,4$, respectively.

\end{lemma}

\begin{proof} Follows from the data in
Theorem
\ref{thm:repa}.

\end{proof}

\section{The eigenvalues and dual eigenvalues}

\medskip
\noindent  
Our next goal is to compute, for each ordered pair $g,h$ of 
adjacent elements in $S_4$, the entries in the  transition 
matrix from the basis $\lbrack g \rbrack $ to the basis
$\lbrack h \rbrack $. In order to describe these entries,
we  make some comments about eigenvalues, and
 define some expressions.  In this section, we focus
on eigenvalues.

\medskip
\noindent
Let $\beta $ denote a scalar in $\K$.
Let $d$ denote a nonnegative integer, and 
let $\sigma_0, \sigma_1, \ldots, \sigma_d $ denote a sequence
of scalars taken from $\K$. We say this sequence
 is $\beta $-{\it recurrent}  
whenever $\sigma_{i-1}-\beta \sigma_i + \sigma_{i+1}$ is independent
of $i$ for $1 \leq i \leq d-1$. Let $\Phi$ denote the Leonard
system in 
Theorem
\ref{thm:classls}. Then by
condition (v) of that theorem, 
 the eigenvalue
sequence and the dual eigenvalue sequence of $\Phi$ are 
$\beta$-recurrent, where $\beta +1 $ is the common value of 
(\ref{eq:betaplusone}). These two sequences are the ones
we wish to discuss  in this section,
but since what we have to say 
about them
applies to all $\beta $-recurrent sequences, we keep things
general.

\medskip
\noindent
We begin by mentioning some well known formula concerning $\beta$-recurrent
sequences. 
Recall ${\tilde {\K}}$ denotes the algebraic closure of 
the field $\K$.

\begin{lemma}
\label{lem:closedf} Let 
$d$ denote a nonnegative integer,
and 
 let $\sigma_0, \sigma_1, \ldots, \sigma_d $ denote a 
sequence of scalars taken from $\K$. 
Let $\beta $ denote a scalar in $\K$, and assume
 $\sigma_0, \sigma_1, \ldots, \sigma_d $  
is 
 $\beta$-recurrent.
Let $q$ denote a nonzero scalar in $\tilde \K$
such that $q+q^{-1} = \beta$. 
\begin{enumerate}
\item
Suppose $q\not=1, q\not=-1$.
Then there exists scalars $a,b,c$ in $\tilde \K$ such that
\begin{equation}
\sigma_i = a + bq^i + cq^{-i} \qquad \qquad (0 \leq i \leq d).
\label{eq:maincase}  
\end{equation}
\item
Suppose $ q=1$. Then there exists scalars $a,b,c $ in $\K$
such that
\beast
\sigma_i = a + bi + ci(i-1)/2 \qquad \qquad (0 \leq i \leq d).
\eeast
\item
Suppose $q=-1$, and that the characteristic of $\K$ is not $2$.
Then 
 there exists scalars $a,b,c $ in $\K$
such that
\beast
\sigma_i = a + b(-1)^i + ci(-1)^i \qquad \qquad (0 \leq i \leq d).
\eeast
\end{enumerate}
Referring to case (ii) above, 
if $\K $ has characteristic 2, we interpret the  expression
$i(i-1)/2 $ as $0$ if $i=0$ or $i=1$ (mod 4), and as $1$ if
$i=2$ or $i=3$ (mod 4).
\end{lemma}

\begin{definition}
\label{def:nbracknv}
Let $q$ denote a nonzero scalar in $\tilde \K$, 
and let $n$ denote an integer. We let 
 $\lbrack n \rbrack_q$ denote the following 
scalar in $\tilde \K$.

\medskip
\noindent First assume $n$ is odd. In this case we define
\begin{eqnarray}
\lbrack n \rbrack_q= 
\cases{{\displaystyle{\frac{q^{n/2}-q^{-n/2}}{q^{1/2}-q^{-1/2}}}}, 
& $\qquad $
if $\quad q\not=1$;\cr
\;\; n, 
& $\qquad $
if $\quad q=1$.\cr}
\label{eq:odddef}
\end{eqnarray}
We observe
\beast
\lbrack n \rbrack_q = q^{(n-1)/2} + q^{(n-3)/2} + \cdots
+ q^{(3-n)/2} +  q^{(1-n)/2} \qquad \qquad (\hbox{if} \quad n> 0)
\eeast
and that $\lbrack -n \rbrack_q = -\lbrack n \rbrack_q$.
For example,
\beast
&&
\lbrack -5 \rbrack_q = -q^2-q-1-q^{-1}-q^{-2}, \qquad 
\lbrack -3 \rbrack_q = -q-1-q^{-1}, \qquad 
\lbrack -1 \rbrack_q = -1, \qquad 
\\
&&
\qquad \lbrack 1 \rbrack_q = 1, \qquad 
\lbrack 3 \rbrack_q = q+1+q^{-1}, \qquad 
\lbrack 5 \rbrack_q = q^2+q+1+q^{-1}+q^{-2}. 
\eeast
Next assume $n$ is even. In this case we define
\begin{eqnarray}
\lbrack n \rbrack_q=
\cases{{\displaystyle{\frac{q^{n/2}-q^{-n/2}}{q-q^{-1}}}}, 
& $\qquad $
if $\quad q\not=1, \quad q\not=-1$;\cr
\;\; n/2, 
& $\qquad $
if $\quad q=1$;\cr
\;\;(-1)^{n/2-1}n/2, 
& $\qquad $ 
if $\quad q=-1$.\cr
}
\label{eq:evenndef}
\end{eqnarray}
We observe
\beast
\lbrack n \rbrack_q = q^{{n/2}-1} + q^{{n/2}-3} + \cdots
+ q^{3-{n/2}} +  q^{1-{n/2}} \qquad \qquad (\hbox{if} \quad n\geq 0)
\eeast
and that $\lbrack -n \rbrack_q = -\lbrack n \rbrack_q$.
For example,
\beast
&&
\lbrack -6 \rbrack_q = -q^2-1-q^{-2}, \qquad 
\lbrack -4 \rbrack_q = -q-q^{-1}, \qquad  
\lbrack -2 \rbrack_q = -1, \qquad  
\lbrack 0 \rbrack_q = 0, \qquad  
\\
&&
\qquad  \lbrack 2 \rbrack_q = 1, \qquad  
\lbrack 4 \rbrack_q = q+q^{-1}, \qquad  
\lbrack 6 \rbrack_q = q^2+1+q^{-2}. 
\eeast
Referring to the cases $q=1,q=-1$ of  
(\ref{eq:evenndef}),
if $\K$ has characteristic $2$, we interpret 
$n/2$ as $1$ if $n=2$  (mod 4), and as $0$ if $n=0$ (mod 4).
\end{definition}

\noindent We mention a handy recursion.

\begin{lemma} 
\label{lem:recnq}
Let $q$ denote a nonzero scalar in $\tilde \K$.
Then for all integers $n$,
\begin{equation}
(q+q^{-1})\lbrack n \rbrack_q = 
\lbrack n+2 \rbrack_q + \lbrack n-2\rbrack_q.
\label{eq:recnq}
\end{equation}
\end{lemma}

\begin{proof}  
Routine calculation using
(\ref{eq:odddef}) and 
(\ref{eq:evenndef}).

\end{proof}

\begin{corollary}
\label{cor:ink}
Let $q$ denote a nonzero element of $\tilde \K$ such that
 $q+q^{-1} \in \K$. 
Then $\lbrack n \rbrack_q \in \K$ for all integers $n$.
\end{corollary}

\begin{proof}
The scalars 
 $\lbrack 0 \rbrack_q$ and 
 $\lbrack 2 \rbrack_q$ are contained in $\K$, since 
 these equal $0$ and $1$, respectively.
By this and a routine induction using
Lemma \ref{lem:recnq},
we find $\lbrack n \rbrack_q $ is contained in $\K$ for
all even integers $n$.
The scalars 
 $\lbrack -1 \rbrack_q$ and 
 $\lbrack 1 \rbrack_q$ are contained in $\K$, since 
 these equal $-1$ and $1$, respectively.
By this and a routine induction using
Lemma \ref{lem:recnq},
we find $\lbrack n \rbrack_q $ is contained in $\K$ for
all odd integers $n$.

\end{proof}

\begin{lemma} 
\label{lem:eigformnv} Let $d$ denote a nonnegative integer,
and let 
  $\sigma_0, \sigma_1, \ldots, \sigma_d $ denote a 
sequence of scalars taken from $\K$.
Let $\beta$ denote a scalar in $\K$, and assume
  $\sigma_0, \sigma_1, \ldots, \sigma_d $ is $\beta$-recurrent. 
Let $q$ denote a nonzero 
scalar in $\tilde \K $ such that $q+q^{-1}=\beta $.
Then for 
 $0 \leq i,j,r,s\leq d$ we have 
\begin{equation}
\lbrack r-s \rbrack_q
(\sigma_i-\sigma_j)
= 
\lbrack i-j\rbrack_q
(\sigma_r-\sigma_s),
\label{eq:handyf}
\end{equation}
provided $i+j=r+s$.
\end{lemma}

\begin{proof} 
Let the integers $i,j,r,s$ be given, and assume
 $i+j=r+s$.  
First suppose $q\not=1$,  $q\not=-1$. Let $n$ denote
the common value of $i+j$, $r+s$, and for convenience
set $e = q^{1/2}-q^{-1/2}$ (if $n$ is odd)
and $e=q-q^{-1}$ (if $n$ is even). Observe $r-s$ and
$r+s=n$ have the same parity, so by Definition
\ref{def:nbracknv},
\begin{eqnarray}
\lbrack r-s \rbrack_q &=& (q^{(r-s)/2} - q^{(s-r)/2})e^{-1}
\nonumber
\\
&=& (q^r-q^s)e^{-1} q^{-n/2}.
\label{eq:iminj}
\end{eqnarray}
Similarly
\begin{eqnarray}
\lbrack i-j \rbrack_q &=&
 (q^i-q^j)e^{-1} q^{-n/2}.
\label{eq:rmins}
\end{eqnarray}
By Lemma 
\ref{lem:closedf}(i),  there exists scalars $a,b,c$ in ${\tilde \K} $
such that $\sigma_0, \sigma_1, \ldots, \sigma_d$
are given by
(\ref{eq:maincase}). Observe
\begin{eqnarray}
\sigma_i - \sigma_j &=& b(q^i-q^j)+c(q^{-i}-q^{-j})
\nonumber
\\
&=& (q^i-q^j)(b-cq^{-n}).
\label{eq:iminj2}
\end{eqnarray}
Similarly
\begin{eqnarray}
\sigma_r - \sigma_s &=& 
(q^r-q^s)(b-cq^{-n}).
\label{eq:rmins2}
\end{eqnarray}
Combining 
(\ref{eq:iminj})--(\ref{eq:rmins2}) we obtain
(\ref{eq:handyf}).
We have now proved the lemma for the case $q\not=1$,
$q\not=-1$. The proof for the cases $q=1$, $q=-1$
is similar, and omitted.
\end{proof}

\medskip
\noindent
Let $q$ denote a nonzero scalar in ${\tilde \K}$,
and let $r,s,t$ denote
 nonnegative integers. A bit later in the paper,
we will define some expressions $\lbrack r,s,t \rbrack_q$ that
make sense under the assumption 
$\lbrack i \rbrack_q \not=0$  
for $1 \leq i \leq r+s+t$. 
We comment on this assumption.
First observe  
$\lbrack 1 \rbrack_q $ and $\lbrack 2 \rbrack_q $ are nonzero,
since these scalars both equal 1. For $i\geq 3$,
it could happen
that 
$\lbrack i \rbrack_q=0 $; we explain how in the rext result.

\begin{lemma}
\label{lem:branotz}
Let $q$ denote a nonzero scalar in $\tilde \K$, and 
let $i$ denote a positive integer.  
Then (i)--(vi) hold below.  
\begin{enumerate}
\item Assume $q\not=1$, $q\not=-1$. Then  $\lbrack i \rbrack_q = 0$
if and only if $q^i=1 $.
\item Assume $q=1$ and that $\K$ has characteristic 0. Then 
$\lbrack i \rbrack_q \not=0 $. 
\item Assume $q=1$ and that $\K$ has characteristic $p$, $p\geq 3$. Then 
$\lbrack i \rbrack_q =0 $  
if and only if $p$ divides $i$.
\item Assume $q=-1$ and that $\K$ has characteristic $0$. Then 
$\lbrack i \rbrack_q \not=0 $. 
\item Assume $q=-1$ and that $\K$ has characteristic $p$, $p\geq 3$. 
Then 
$\lbrack i \rbrack_q =0 $  
if and only if $ 2p$ divides $i$.
\item Assume $q=1$ and that $\K$ has characteristic $2$. Then 
$\lbrack i \rbrack_q =0 $   
if and only if $4$ divides $i$. 
\end{enumerate}
\end{lemma}

\begin{proof}
First assume $q\not=1$, $q\not=-1$. Then 
$\lbrack i \rbrack_q $ is a nonzero scalar multiple of
$q^i-1$ by
Definition
\ref{def:nbracknv},
and assertion (i) follows.
Next assume $q=1$ and that the characteristic 
of  $\K$ is not 2. Then the sequence 
$\lbrack 1 \rbrack_q, \lbrack 2 \rbrack_q, \ldots $ is given
by $1,1,3,2,5,3,7,4\ldots $ and
assertions (ii), (iii) follow. 
Next assume 
$q=-1$ and that the characteristic 
of  $\K$ is not 2. Then the sequence 
$\lbrack 1 \rbrack_q, \lbrack 2 \rbrack_q, \ldots $ is given
by
 $1,1,-1,-2,1,3,-1,-4,\ldots $
and assertions (iv), (v) follow.
Now assume $q=1$ and that $\K$ has characteristic 2.
Then 
the sequence 
$\lbrack 1 \rbrack_q, \lbrack 2 \rbrack_q, \ldots $ is given
by
 $1,1,1,0,1,1,1,0,\ldots $ and
assertion (vi) follows.

\end{proof}

\begin{lemma}
\label{lem:restnv} 
Let $d$ denote an integer at least 3.
Let  $\sigma_0, \sigma_1, \ldots, \sigma_d $  
 denote a sequence of distinct scalars taken from $\K$,
 and assume 
\begin{equation} 
 \frac{\sigma_{i-2}-\sigma_{i+1}}{\sigma_{i-1}-\sigma_i}
\label{eq:sigfrac}
\end{equation}
is independent of $i$ for $2 \leq i \leq d-1$.
Let $q$ denote a nonzero scalar in $\tilde \K$ such that
$q+q^{-1}+1$ equals the common value of  
(\ref{eq:sigfrac}).
Then 
$\lbrack i \rbrack_q \not=0$ for $1 \leq i \leq d$.
\end{lemma}

\begin{proof} Abbreviate $\beta =q+q^{-1}$, and observe
 $\sigma_0, \sigma_1, \ldots, \sigma_d $   is $\beta$-recurrent.
First suppose $q\not=1$ and $q\not=-1$. Then for $1 \leq i \leq d$
we have $q^i\not=1$; otherwise $\sigma_i = \sigma_0$
by Lemma 
\ref{lem:closedf}(i). The result now follows
by
Lemma \ref{lem:branotz}(i).
Next suppose $q=1$ and that $\K$ has characteristic 0.
Then the result holds by
Lemma \ref{lem:branotz}(ii).
Next suppose $q=1$ and that $\K$ has characteristic 
$p$, $p\geq 3$. Then $d <p$; otherwise $\sigma_p=\sigma_0$
in view of Lemma
\ref{lem:closedf}(ii). The result now follows
by Lemma \ref{lem:branotz}(iii).
Next suppose $q=-1$ and that $\K$ has characteristic 0.
Then the result holds by
Lemma \ref{lem:branotz}(iv).
Next suppose $q=-1$ and that $\K$ has characteristic 
$p$, $p\geq 3$. Then $d <2p$; otherwise $\sigma_{2p}=\sigma_0$
in view of Lemma
\ref{lem:closedf}(iii). The result now follows
by Lemma \ref{lem:branotz}(v).
Now suppose $q=1$ and that $\K$ has characteristic 2.
Then  $d\leq 3$; otherwise $\sigma_4 = \sigma_0$ by
Lemma \ref{lem:closedf}(iii) and the comment at the end of
that lemma. The result now follows by Lemma
\ref{lem:branotz}(vi).

\end{proof}

\begin{corollary}
\label{lem:dennonz}
Let $\Phi$ denote the Leonard system in
(\ref{eq:ourstartingpt}), and assume $d\geq 3$.
Let $q$  denote a 
nonzero scalar in $\tilde \K $ such that
$q+q^{-1} + 1$ equals the common value of 
(\ref{eq:betaplusone}).  
Then $\lbrack i \rbrack_q \not=0$ for $1 \leq i \leq d$.
\end{corollary}

\begin{proof} 
Apply Lemma 
\ref{lem:restnv} to the eigenvalue sequence of $\Phi$. 

\end{proof}

\noindent We finish this section with a definition. 

\begin{definition}
Let $q$ denote a nonzero scalar in $\tilde \K$.
For each nonnegative integer $n$ we define
\begin{equation}
\lbrack n \rbrack!_q = 
\lbrack 1 \rbrack_q \lbrack 2 \rbrack_q \cdots \lbrack n \rbrack_q.
\label{eq:nqfact}
\end{equation}
We interpret $\lbrack 0 \rbrack !_q=1$.
\end{definition}

\section{The  scalars $\lbrack r,s,t \rbrack_q$}

\medskip
\noindent  
A bit later in the paper we will compute,
for each ordered pair $g,h$ of 
adjacent elements in $S_4$, the entries in the  transition 
matrix from the basis $\lbrack g \rbrack $ to the basis
$\lbrack h \rbrack $.
Among the entries in these transition matrices,
we will encounter an expression that occurs so often
we will give it a name. The details are in the following
definition. 

\begin{definition}
\label{def:bra} Let $q$ denote a nonzero scalar in $\tilde \K$ and
let $r,s,t$ denote nonnegative integers. 
We define the expressions $(r,s,t)_q$ and $\lbrack r,s,t\rbrack_q $
as follows.
We set
\begin{equation}
(r,s,t)_q = 
\cases{ q+q^{-1}+2 
&$\qquad $if each of $r,s,t$ is odd;\cr
\;\;1,  & $\qquad $if at least one of $r,s,t$ is even.\cr}
\label{eq:curvedef}
\end{equation}
Next assume $\lbrack i \rbrack_q \not=0 $ for $1 \leq i \leq r+s+t $.
Then we set 
\begin{equation}
\lbrack r,s,t \rbrack_q   
= \frac{\lbrack r+s \rbrack !_q
\lbrack r+t \rbrack !_q 
\lbrack s+t \rbrack !_q
(r,s,t)_q}{\lbrack r\rbrack !_q  
\lbrack s \rbrack !_q \lbrack t \rbrack !_q  \lbrack r+s+t \rbrack !_q}.
\label{eq:triplemain}
\end{equation}
We remark $\lbrack r,s,t\rbrack_q \in \K$ provided $q+q^{-1} \in \K$. 
Moreover, 
$\lbrack r,s,t\rbrack_q=1$ if at least one of $r,s,t$ equals $0$. 
\end{definition}

\medskip
\noindent Referring to the above definition,
to get a better appreciation for
$\lbrack r,s,t\rbrack_q$  
we now evaluate the expression
on the right in
(\ref{eq:triplemain}) 
using
Definition
\ref{def:nbracknv}. To express our results, 
we use the following notation.
For all $a, q \in  
{\tilde \K}$
we define
\beast
(a;q)_n = (1-a)(1-aq)\cdots (1-aq^{n-1}) \qquad \qquad n=0,1,2,\ldots
\eeast
and interpret $(a;q)_0=1$.

\begin{lemma}
\label{lem:branv}
Let $q$ denote a nonzero scalar in $\tilde \K$,
let $r,s,t$ denote nonnegative integers, and assume
$\lbrack i \rbrack_q \not=0$ for $1 \leq i \leq r+s+t$.
\begin{enumerate}
\item Suppose $q\not=1$, $q\not=-1$. Then
\begin{equation}
\lbrack r,s,t \rbrack_q = \frac{(q;q)_{r+s}(q;q)_{r+t}(q;q)_{s+t}}{(q;q)_r(q;q)_s(q;q)_t(q;q)_{r+s+t}}.
\label{eq:nextlevel}
\end{equation}

\item  Suppose $q= 1$ and that the characteristic of $\K$ is not 2. Then 
\begin{equation}
\lbrack r,s,t \rbrack_q= {{(r+s)!\,(r+t)!\,(s+t)!}\over{r!\,s!\,t!\,(r+s+t)!}}.
\label{eq:nextlevelone}
\end{equation}

\item Suppose $q=-1$ 
 and that the characteristic of $\K$ is not 2. If 
each of $r,s,t$ is odd, then  
$\lbrack r,s,t\rbrack_q =0$. If at least one of $r,s,t$ is even,
then
\begin{equation}
\lbrack r,s,t \rbrack_q = 
\frac{\lfloor \frac{r+s}{2}\rfloor ! 
\lfloor \frac{r+t}{2}\rfloor !
\lfloor \frac{
s+t}{2}\rfloor !}{\lfloor \frac{r}{2}\rfloor ! 
\lfloor \frac{s}{2}\rfloor ! 
\lfloor \frac{t}{2}\rfloor !  \lfloor \frac{r+s+t}{2} \rfloor !} .
\label{eq:nextlevelminone}
\end{equation}
The expression $\lfloor n \rfloor $ denotes the greatest integer less
than or equal to $n$.

\item Suppose $q=1$ and that $\K$ has characteristic 2. Recall in
this case  $r+s+t\leq 3$
by
Lemma
\ref{lem:branotz}(vi).
If each of $r,s,t$ equals 1, then $\lbrack r,s,t\rbrack_q =0$.
If at least one of $r,s,t$ equals 0 then $\lbrack r,s,t \rbrack_q=1$. 
\end{enumerate}
Concerning the expressions on the right in
(\ref{eq:nextlevel}),
(\ref{eq:nextlevelone}),
(\ref{eq:nextlevelminone}), the
 denominator is nonzero by Lemma
\ref{lem:branotz}.
\end{lemma}

\begin{proof} Evaluate 
(\ref{eq:triplemain}) using
Definition
\ref{def:nbracknv},
(\ref{eq:nqfact}), and
(\ref{eq:curvedef}).

\end{proof}

\medskip
\noindent We will need the following identity.

\begin{lemma}
\label{lem:ident} Let $q$ denote a nonzero scalar in $\tilde \K$, 
and let $r,s,t$ denote positive integers. Assume
$\lbrack i \rbrack_q \not=0$ for $1 \leq i <r+s+t$. 
Then with reference to 
Definition
\ref{def:bra} 
 we have
\begin{equation}
\lbrack r-t\rbrack_q \lbrack r+t\rbrack^{-1}_q
\lbrack r,s-1,t\rbrack_q
=\lbrack r-1,s,t\rbrack_q - \lbrack r,s,t-1 \rbrack_q.
\label{eq:ident}
\end{equation}
\end{lemma}

\begin{proof} First assume $q\not=1$ and $q\not=-1$. By 
Definition
\ref{def:nbracknv},
and since the integers $r+t, r-t$ have the same parity, we find
\begin{eqnarray}
\frac{\lbrack r-t\rbrack_q }{ \lbrack r+t\rbrack_q} &=&
\frac{q^{(r-t)/2} - q^{(t-r)/2}}{q^{(r+t)/2}-q^{-(r+t)/2}}
\nonumber
\\
&=& \frac{q^t-q^r}{1-q^{r+t}}.
\label{eq:rptrmt}
\end{eqnarray}
Using 
(\ref{eq:nextlevel}), we obtain
\begin{eqnarray}
\lbrack r-1,s,t\rbrack_q &=& x(1-q^{s+t})(1-q^r),
\label{eq:rmone}
\\
\lbrack r,s-1,t\rbrack_q &=& x(1-q^{r+t})(1-q^s),
\label{eq:smone}
\\
\lbrack r,s,t-1\rbrack_q &=& x(1-q^{r+s})(1-q^t),
\label{eq:tmone}
\end{eqnarray}
where
\beast
x = \frac{(q;q)_{r+s-1}(q;q)_{r+t-1}(q;q)_{s+t-1}}{(q;q)_r(q;q)_s(q;q)_t(q;q)_{r+s+t-1}}.
\eeast
One readily verifies
\begin{equation}
(q^t-q^r)(1-q^s) = (1-q^{s+t})(1-q^r)- (1-q^{r+s})(1-q^t).
\label{eq:fif}
\end{equation}
Multiplying both sides of 
(\ref{eq:fif}) by $x$, and evaluating the result using
(\ref{eq:rptrmt})--(\ref{eq:tmone}), we routinely obtain
(\ref{eq:ident}).
We have now proved the result for the case $q\not=1$, $q\not=-1$.
The proof for the  cases $q=1$, $q=-1$ are similar, and omitted.
\end{proof}

\section{The scalars $\varepsilon_0,\varepsilon_d, \varepsilon^*_0, 
\varepsilon^*_d$ }

\medskip
\noindent  
In the next section we will compute,
for each  ordered pair $g,h$ of 
adjacent elements in $S_4$, the entries in the  transition 
matrix from the basis $\lbrack g \rbrack $ to the basis
$\lbrack h \rbrack $.
Recall our 24 bases
are constructed using four vectors
$\eta_0, \eta_d, \eta^*_0, \eta^*_d$, and each of these vectors 
is determined only up to multiplication by a nonzero scalar.
To account for this, we introduce four scalars
$\varepsilon_0,\varepsilon_d, \varepsilon^*_0, 
\varepsilon^*_d$.

\medskip
\noindent
For convenience, we make the following definition.
\begin{definition}
\label{def:etilde}
Let $\Phi$ denote the Leonard system in
(\ref{eq:ourstartingpt}).
We define
\begin{eqnarray}
{\tilde E}_0 &=& (A-\theta_1I)(A-\theta_2I) \cdots (A-\theta_d I),
\label{eq:tildeedef}
\\
{\tilde E}_d &=& (A-\theta_0I)(A-\theta_1I) \cdots (A-\theta_{d-1} I),
\\
{\tilde E}^*_0 &=& (A^*-\theta^*_1I)(A^*-\theta^*_2I) \cdots 
(A^*-\theta^*_d I),
\\
{\tilde E}^*_d &=& (A^*-\theta^*_0I)(A^*-\theta^*_1I) \cdots
(A^*-\theta^*_{d-1} I),
\end{eqnarray}
where the $\theta_i, \theta^*_i$ are from Definition
\ref{def:evseq}.
\end{definition}

\begin{lemma}
\label{lem:etildee}
Let $\Phi$ denote the Leonard system in
(\ref{eq:ourstartingpt}).
Then with reference to 
Definition 
\ref{def:etilde},
\begin{enumerate}
\item $
{\tilde E}_0 = E_0(\theta_0-\theta_1)(\theta_0-\theta_2) 
\cdots (\theta_0-\theta_d), $
\item
${\tilde E}_d = E_d(\theta_d-\theta_0)(\theta_d-\theta_1) \cdots 
(\theta_d-\theta_{d-1}),
$
\item
$
{\tilde E}^*_0 = E^*_0(\theta_0^*-\theta^*_1)(\theta_0^*-\theta^*_2) \cdots 
(\theta_0^*-\theta^*_d) ,
$
\item
$
{\tilde E}^*_d = E^*_d(\theta_d^*-\theta^*_0)(\theta_d^*-\theta^*_1) \cdots
(\theta_d^*-\theta^*_{d-1}).
$
\end{enumerate}
\end{lemma}

\begin{proof}
To get (i), set $i=0$ in  (\ref{eq:primiddef}) and compare
the result with
(\ref{eq:tildeedef}).
Assertions (ii)-(iv) are similarily proved.

\end{proof}

\begin{lemma}
\label{lem:tildeaction}
Let $\Phi$ denote a Leonard system in
(\ref{eq:ourstartingpt}).
Let $g$ denote the element $d^*00^*d$
of $S_4$ and recall by 
Lemma \ref{lem:reconcile} that 
$\lbrack g \rbrack $ 
is a  $\Phi$-split basis. 
For $0 \leq i,j\leq d$,
the $ij^{\hbox{th}}$ entry of the matrices
${\tilde E}^g_0$,
${\tilde E}^g_d$,
${\tilde E}^{*g}_0$,
${\tilde E}^{*g}_d$ are given as follows.
\begin{enumerate}
\item 
 The $ij^{\hbox{th}}$ entry of 
${\tilde E}^g_0 $ is
\beast
(\theta_0-\theta_{i+1})(\theta_0-\theta_{i+2}) \cdots 
(\theta_0-\theta_d)
\eeast
if $j=0$, and $0$ if $j\not=0$.
\item
The 
 $ij^{\hbox{th}}$ entry of 
${\tilde E}^g_d $ is 
\beast
(\theta_d-\theta_0)(\theta_d-\theta_1) \cdots 
(\theta_d-\theta_{j-1}) 
\eeast
if $i=d$, and $0$ if
$i\not=d$.
\item The $ij^{\hbox{th}}$ entry
of 
 ${\tilde E}^{*g}_0 $ is   
\beast
(\theta^*_0-\theta^*_{j+1})(\theta^*_0-\theta^*_{j+2}) \cdots 
(\theta^*_0-\theta^*_d)\varphi_1 \varphi_2 \cdots \varphi_j 
\eeast
if $i=0$, and $0$ if $i\not=0$. 
\item 
 The $ij^{\hbox{th}}$ entry of 
${\tilde E}^{*g}_d$ is 
\beast
(\theta^*_d-\theta^*_0)(\theta^*_d-\theta^*_1) \cdots 
(\theta^*_d-\theta^*_{i-1}) \varphi_{i+1} \varphi_{i+2} \cdots \varphi_d 
\eeast
if $j=d$, and $0$ if $j\not=d$. 
\end{enumerate}
\end{lemma}

\begin{proof} The entries of $E^g_0$, $E^g_d$, $E^{*g}_0$,
$E^{*g}_d$ are given in 
\cite[Thm. 4.8]{LS99}.  Using these entries and Lemma
\ref{lem:etildee}, we routinely obtain the assertions of the present
lemma.

\end{proof}

\noindent For notational convenience, we introduce the following
notation.

\begin{definition} 
\label{def:phinot}
Let $\Phi$ denote the Leonard system in
(\ref{eq:ourstartingpt}).
We set
\begin{equation}
\varphi:=\varphi_1 \varphi_2 \cdots \varphi_d, \qquad \qquad 
\phi:=\phi_1 \phi_2 \cdots \phi_d,
\label{eq:phiabb}
\end{equation}
where $\varphi_1, \varphi_2, \ldots, \varphi_d$ denotes the
first split sequence of $\Phi$ and 
where $\phi_1, \phi_2, \ldots, \phi_d$ denotes the
second split sequence of $\Phi$. We observe by
Theorem \ref{thm:classls}(i) that  $\varphi\not=0$,
$\phi\not=0$.

\end{definition}

\begin{lemma} 
\label{lem:tractildee}
Let $\Phi$ denote the Leonard system in
(\ref{eq:ourstartingpt}).
Then 
with reference to 
Definition \ref{def:etilde}, 
the trace of each of 
${\tilde E}_d{\tilde E}^*_0 $, 
${\tilde E}_0{\tilde E}^*_d $ equals 
$\varphi $.
Moreover,  the trace of each of ${\tilde E}_0{\tilde E}^*_0 $, 
${\tilde E}_d{\tilde E}^*_d $ equals $\phi $.
\end{lemma}

\begin{proof} Using the data in Lemma
\ref{lem:tildeaction}, we routinely find 
the trace of ${\tilde E}_d{\tilde E}^*_0 $ equals 
$\varphi $.
To obtain the remaining assertions, apply this result to 
the relatives of $\Phi$, and use 
Theorem
\ref{thm:phimod}.

\end{proof}

\begin{lemma}
\label{lem:backandforth}
Let $\Phi$ denote the Leonard system in
(\ref{eq:ourstartingpt}).
Then with reference to  
Definition \ref{def:etilde},
\begin{eqnarray}
&&{\tilde E}^*_0 {\tilde E}_d {\tilde E}^*_0 = 
\varphi {\tilde E}^*_0,
\qquad \qquad  
{\tilde E}_d {\tilde E}^*_0 {\tilde E}_d = 
\varphi {\tilde E}_d,
\label{eq:bkf1}
\\
&&{\tilde E}_0 {\tilde E}^*_d {\tilde E}_0 = 
\varphi {\tilde E}_0,
\qquad \qquad  
{\tilde E}^*_d {\tilde E}_0 {\tilde E}^*_d = 
\varphi {\tilde E}^*_d,
\label{eq:bkf2}
\\
&&{\tilde E}_0 {\tilde E}^*_0 {\tilde E}_0 = 
\phi {\tilde E}_0,
\qquad \qquad  
{\tilde E}^*_0 {\tilde E}_0 {\tilde E}^*_0 = 
\phi {\tilde E}^*_0,
\label{eq:bkf3}
\\
&& {\tilde E}_d {\tilde E}^*_d {\tilde E}_d = 
\phi {\tilde E}_d,
\qquad  \qquad  
{\tilde E}^*_d {\tilde E}_d {\tilde E}^*_d = 
\phi {\tilde E}^*_d.
\label{eq:bkf4}
\end{eqnarray}
\end{lemma}

\begin{proof} We first prove  the equation on the left
in (\ref{eq:bkf1}). Since $E^*_0$ is a rank one idempotent,
and since ${\tilde E}^*_0 $ is a nonzero scalar multiple
of 
 $E^*_0$, 
 there exists a scalar
$\alpha \in \K $ such that
${\tilde E}^*_0 {\tilde E}_d {\tilde E}^*_0 = 
\alpha {\tilde E}^*_0 $. 
We show $\alpha  = \varphi$. 
We mentioned 
 ${\tilde E}^*_0 $ is a nonzero scalar multiple
of 
 $E^*_0$, 
so
\begin{equation}
E^*_0 {\tilde E}_d {\tilde E}^*_0 = 
\alpha E^*_0.
\label{eq:middle}
\end{equation}
We take the trace of each side of    
(\ref{eq:middle}).
Observe the trace of $E^*_0$ equals 1, so the trace of the right
side of 
(\ref{eq:middle})
equals $\alpha$. Since $XY$ and $YX$ have the same trace, and using   
${\tilde E}^*_0 E^*_0 ={\tilde E}^*_0$, we find in view of 
Lemma \ref{lem:tractildee} that the trace of the left side of
(\ref{eq:middle})
equals $\varphi$. Apparently $\alpha = \varphi$,
and this implies the equation on the left in
in (\ref{eq:bkf1}).
Applying this result to the relatives of $\Phi$, we
obtain the remaining assertions.
\end{proof}

\begin{lemma}
\label{lem:trick}
Let $\Phi$ denote the Leonard system in
(\ref{eq:ourstartingpt}).
Then with reference to  
Definition \ref{def:etilde},
we have  the following.
\begin{eqnarray}
&& {\tilde E}^*_d {\tilde E}_0 {\tilde E}^*_0 = 
 {\tilde E}^*_d {\tilde E}_d {\tilde E}^*_0, \qquad \qquad 
{\tilde E}_d {\tilde E}^*_0 {\tilde E}_0 = 
 {\tilde E}_d {\tilde E}^*_d {\tilde E}_0,
\label{eq:tripleprod}
\\
&&
{\tilde E}_0 {\tilde E}^*_0 {\tilde E}_d = 
 {\tilde E}_0 {\tilde E}^*_d {\tilde E}_d,
\qquad \qquad 
{\tilde E}^*_0 {\tilde E}_0 {\tilde E}^*_d = 
 {\tilde E}^*_0 {\tilde E}_d {\tilde E}^*_d.
\label{eq:tripleprod2}
\end{eqnarray}
\end{lemma}

\begin{proof} The equation on the left in
(\ref{eq:tripleprod})
is readily obtained
 using the matrix representations
given in Lemma
\ref{lem:tildeaction}.
Applying this equation to the relatives of $\Phi$,
we 
obtain the remaining equations
in 
(\ref{eq:tripleprod}), 
(\ref{eq:tripleprod2}).

\end{proof}

\begin{lemma}
\label{lem:newvi}
Let $\Phi$ denote the Leonard system in
(\ref{eq:ourstartingpt}), and 
let $V$ denote the irreducible left $\cal A$-module.
Let 
$\eta_0, \eta_d, \eta^*_0,  \eta^*_d$ denote nonzero vectors in
$V$ that satisfy
(\ref{eq:videf}). Then there exists 
nonzero scalars 
$\varepsilon_0, \varepsilon_d, \varepsilon^*_0,\varepsilon^*_d $
in $\K$ 
such that
\begin{eqnarray}
&&{\tilde E}_d\eta^*_0/\varepsilon^*_0 = \eta_d/\varepsilon_d, \qquad  \qquad 
{\tilde E}_d\eta^*_d/\varepsilon^*_d = \eta_d/\varepsilon_d, 
\label{eq:ed}
\\
&&
{\tilde E}^*_0\eta_0/\varepsilon_0 =  \eta^*_0/\varepsilon^*_0,
\qquad \qquad 
{\tilde E}^*_0\eta_d/\varepsilon_d = \varphi  
\eta^*_0/\varepsilon^*_0,  
\label{eq:e0s}
\\
&&{\tilde E}^*_d\eta_0/\varepsilon_0 =  \eta^*_d/\varepsilon^*_d,
\qquad \qquad 
{\tilde E}^*_d\eta_d/\varepsilon_d = \phi 
\eta^*_d/\varepsilon^*_d,  
\label{eq:eds}
\\
&&{\tilde E}_0\eta^*_0/\varepsilon^*_0 = \phi 
\eta_0/\varepsilon_0,  
\qquad \quad \;\,
{\tilde E}_0\eta^*_d/\varepsilon^*_d = \varphi  
\eta_0/\varepsilon_0. 
\label{eq:e0}
\end{eqnarray}
\end{lemma}

\begin{proof}
Let $\varepsilon^*_0$ denote an arbitrary nonzero
scalar in $\K$. To obtain $\varepsilon_0$,
consider the basis 
 $\lbrack d^*d00^* \rbrack $  from  
 the table of Theorem
\ref{thm:bases}, row 10.
Using  
(\ref{eq:tildeedef}), we recognize  
 the vector ${\tilde E}_0 \eta^*_0$ is
 the $0^{\hbox{th}}$ 
vector in this basis.
By Theorem
\ref{thm:bases}, we find  
 ${\tilde E}_0 \eta^*_0$ is a basis for $E_0V$.
By the construction $\eta_0$ is a basis for $E_0V$,
so 
 ${\tilde E}_0 \eta^*_0$ is  a nonzero scalar multiple of $\eta_0$. 
Apparently, there exists a nonzero scalar
 $\varepsilon_0 \in \K$ 
 that satisfies the equation on the
 left in 
(\ref{eq:e0}).
Similarly,
there exists 
 nonzero scalars
 $\varepsilon_d, 
 \varepsilon^*_d $
in $\K$ 
 that satisfy the equations on the
 left in 
(\ref{eq:ed}),
(\ref{eq:eds}), respectively.
To obtain the equation on the right in
(\ref{eq:eds}), apply the equation on the left
in 
(\ref{eq:tripleprod}) to $\eta^*_0/\varepsilon^*_0$,
and evaluate the result using
$E^*_0 \eta^*_0 = \eta^*_0$, Lemma
\ref{lem:etildee}(iii), and the equations on the left
in 
(\ref{eq:ed}),
(\ref{eq:eds}),
(\ref{eq:e0}).
To obtain the equation on the left in 
(\ref{eq:e0s}), apply the equation on the right in
(\ref{eq:bkf3}) to $\eta^*_0/\varepsilon^*_0$, and evaluate
the result using 
$E^*_0 \eta^*_0 = \eta^*_0$, Lemma
\ref{lem:etildee}(iii), and the equation on the left
in 
(\ref{eq:e0}).
The equations on the right in 
(\ref{eq:ed}), (\ref{eq:e0s}), (\ref{eq:e0}) are similarly 
obtained.

\end{proof}

\medskip
\begin{note}  The scalars 
$\varepsilon_0, \varepsilon_d, \varepsilon^*_0,\varepsilon^*_d $
from  Lemma
\ref{lem:newvi} are ``free'' in the following sense.
Let $\Phi$ denote the Leonard system
in
(\ref{eq:ourstartingpt}), and let 
$V$ denote the irreducible left $\cal A$-module.
Let 
 $\varepsilon_0, \varepsilon_d, \varepsilon^*_0,\varepsilon^*_d $
denote arbitrary nonzero scalars in $\K$. Then
there exists nonzero vectors $\eta_0, \eta_d, \eta^*_0, \eta^*_d$
in $V$ that satisfy 
(\ref{eq:videf}) and 
(\ref{eq:ed})--(\ref{eq:e0}).

\end{note}

\begin{note}
\label{note:nosym}
The reader may notice a certain lack of
symmetry in the definition of
$\varepsilon_0, \varepsilon_d, \varepsilon^*_0$, $\varepsilon^*_d $.
We accept this asymmetry to avoid introducing the square
roots of  $\varphi$ and $\phi$. We remark these square roots
may not be in $\K$.
To display the underlying symmetry in 
(\ref{eq:ed})--(\ref{eq:e0}), make the following  change of
variables:
\beast
\varepsilon_0 = \varepsilon'_0,  \qquad 
\varepsilon_d = \varepsilon^{\prime}_d \varphi^{-1/2}\phi^{-1/2} \qquad
\varepsilon^*_0 = \varepsilon^{*\prime}_0 \phi^{-1/2},  \qquad 
\varepsilon^*_d = \varepsilon^{*\prime}_d \varphi^{-1/2}.  
\eeast

\end{note}

\noindent The following  equations will be useful.

\begin{lemma}
\label{lem:eadjust}
Let $\Phi$ denote the Leonard system in
(\ref{eq:ourstartingpt}), and 
 let  $V$ denote the irreducible left $\cal A$-module. Let
$\eta_0, \eta_d, \eta^*_0,  \eta^*_d$ denote nonzero vectors in
$V$ that satisfy
(\ref{eq:videf}). 
Let
the scalars
$\varepsilon_0, \varepsilon_d, \varepsilon^*_0,\varepsilon^*_d $
be as in Lemma 
\ref{lem:newvi}.
Then
\begin{eqnarray}
&&E_d \eta^*_d /\varepsilon^*_d = E_d\eta^*_0/\varepsilon^*_0,
\qquad \qquad 
E^*_0 \eta_d/\varepsilon_d = \varphi  E^*_0\eta_0/\varepsilon_0,
\label{eq:dds}
\\
&&E^*_d \eta_d/\varepsilon_d = \phi  
E^*_d\eta_0/\varepsilon_0, \qquad \quad  \,
E_0 \eta^*_d/\varepsilon^*_d =  
\varphi / \phi
 E_0\eta^*_0/\varepsilon^*_0.
\label{eq:dds2}
\end{eqnarray}

\end{lemma}

\begin{proof} 
First consider the equation on the left in
(\ref{eq:dds}). Comparing the two equations 
in 
(\ref{eq:ed}),
we find 
${\tilde E}_d \eta^*_d /\varepsilon^*_d 
= {\tilde E}_d\eta^*_0/\varepsilon^*_0$. Recall $E_d$
is a scalar multiple of ${\tilde E}_d$, so 
$ E_d \eta^*_d /\varepsilon^*_d 
=  E_d\eta^*_0/\varepsilon^*_0$. We now have the equation
on the left in 
(\ref{eq:dds}). The remaining equations
in 
(\ref{eq:dds}), 
(\ref{eq:dds2}) are similarly proved.

\end{proof}

\section{The $24$ bases; transition matrices}

\noindent 
Let $\Phi$ denote the Leonard system in 
(\ref{eq:ourstartingpt}), and let $V$ denote the irreducible
left $\cal A$-module. 
For each element $g\in S_4$, 
we displayed in Theorem \ref{thm:bases} a basis for $V$,
denoted $\lbrack g \rbrack$.
In this section we compute, for each ordered pair $g,h$
of adjacent elements of $S_4$, the entries in
the transition matrix from the basis 
$\lbrack g \rbrack $ to the basis $\lbrack h \rbrack $.

\medskip
\noindent 
We mention a few points from linear algebra.
In line 
(\ref{eq:transdefpre})
we recalled the notion of a transition matrix.
We now recall the closely related concept  of an intertwining matrix.  
Let $g, h$ denote elements of $S_4$, and consider
the corresponding bases $\lbrack  g \rbrack $, $\lbrack h \rbrack $
of $V$.
By an {\it intertwining matrix} from 
$\lbrack  g \rbrack $ to  $\lbrack h \rbrack $,
 we mean a nonzero matrix
$S \in \hbox{Mat}_{d+1}(\K) $ satisfying
\beast
X^g S = S X^h \qquad \qquad (\forall X \in \alg).
\eeast
We observe a matrix in  
$\hbox{Mat}_{d+1}(\K) $ is an intertwining matrix
from
$\lbrack  g \rbrack $ to  $\lbrack h \rbrack $ if and only
if 
 it is a nonzero
scalar multiple of the transition matrix from
$\lbrack  g \rbrack $ to  $\lbrack h \rbrack $.

\medskip
\noindent
The following matrix will play a role in our discussion. 
We let $Z$ denote the matrix in 
$\hbox{Mat}_{d+1}(\K)$ with  entries 
\begin{equation}
Z_{ij} = \cases{1, &if $\;i+j=d$;\cr
0, &if $\;i+j \not=d$\cr}
\qquad \qquad 
(0 \leq i,j\leq d).
\label{eq:zmat}
\end{equation}
We observe $Z^2=I$.

\begin{lemma}
\label{lem:whentrans}
Let $\Phi$ denote the Leonard system in
(\ref{eq:ourstartingpt}), and let $g, h$ denote elements
in $S_4$. Then for all $S \in 
 \hbox{Mat}_{d+1}(\K)$, the following are equivalent.
\begin{enumerate}
\item $S$ is an intertwining matrix from $\lbrack g \rbrack $ to 
$\lbrack h \rbrack $.
\item $S$ is nonzero and both
\begin{equation}
A^g S = S A^h, \qquad \qquad  A^{*g}S = SA^{*h}.
\label{eq:twocond}
\end{equation}
\end{enumerate}
\end{lemma}
\begin{proof} The implication $(i)\rightarrow (ii)$ is clear,
so consider the implication $(ii)\rightarrow (i)$. Let $T$ denote
the transition matrix from $\lbrack g \rbrack $ to $\lbrack h \rbrack $.
We show $S$ is a nonzero
scalar multiple of $T$. Since $T$ is the transition matrix
from $\lbrack g \rbrack $ to $\lbrack h \rbrack $, it is an intertwining matrix
from $\lbrack g \rbrack $ to $ \lbrack h \rbrack $. Therefore
\begin{equation}
A^g T = T A^h, \qquad \qquad  A^{*g}T = TA^{*h}.
\label{eq:twoconds}
\end{equation}
Combining 
(\ref{eq:twocond}), 
(\ref{eq:twoconds}), we find $ST^{-1}$ commutes with both
$A^g$ and $A^{*g}$.
We mentioned the map 
$X \rightarrow X^g$ from $\cal A$ to $\hbox{Mat}_{d+1}(\K)$
is an isomorphism of $\K$-algebras. 
Combining this with our previous comment and using
Corollary
\ref{cor:rig},
we see  
  $ST^{-1}$ is a scalar multiple of the identity.
 Denoting this scalar  by $\alpha$ we have
  $S=\alpha T$. We observe $\alpha \not=0$
 since $S\not=0$. Apparently $S$ is a nonzero scalar multiple
 of $T$, so 
 $S$ is an intertwining matrix from 
 $\lbrack g \rbrack $ to 
$\lbrack h \rbrack $.
\end{proof}

\begin{theorem}
\label{thm:trans} Let $\Phi$ denote the Leonard system in
(\ref{eq:ourstartingpt}). With reference to
Definition
\ref{def:S4interp},
let $wxyz$ denote  an element of $S_4$, and consider
the transition matrices  from  
the basis $\lbrack wxyz \rbrack $ to the  bases
\begin{equation} 
 \lbrack xwyz \rbrack, \qquad 
 \lbrack wyxz \rbrack, \qquad 
 \lbrack wxzy \rbrack.
\label{eq:threetargets}
\end{equation}
The first and second transition matrices are diagonal and lower triangular,
respectively, and their entries are given in the following tables.
The third transition matrix is the matrix $Z$
from
(\ref{eq:zmat}).
In the tables below,
$\theta_0, \theta_1,\ldots, \theta_d $  
(resp. $\theta^*_0, \theta^*_1,\ldots, \theta^*_d $)
denotes the eigenvalue sequence
(resp. dual eigenvalue sequence) for $\Phi$. Moreover
$\varphi_1, \varphi_2,\ldots, \varphi_d $ (resp. 
$\phi_1, \phi_2,\ldots, \phi_d $) denotes the 
first split sequence (resp. second split sequence)
for $\Phi$. 
 The scalars $\varphi, \phi$ are from 
(\ref{eq:phiabb}), and 
the scalars $\varepsilon_0, \varepsilon_d,
\varepsilon^*_0, \varepsilon^*_d$ are from Lemma 
\ref{lem:newvi}.

\medskip

\centerline{
\begin{tabular}[t]{c|c|c}
$wxyz$ &  $\lbrack wxyz \rbrack \rightarrow \lbrack xwyz \rbrack $ & 
 $\lbrack wxyz \rbrack \rightarrow \lbrack wyxz \rbrack $ \\
 & $ii$ entry & $ij$ entry $(i\geq j)$
\\
\hline \hline 
&&
\\
$d^*00^*d$
&
$ \frac{1}{\varphi_1 \cdots \varphi_i} 
\;\frac{\varepsilon_d \varphi}{\varepsilon^*_0}$
&
$\frac{1}{(\theta_{j}-\theta_{0})\cdots (\theta_j - \theta_{j-1})}\;
\frac{1}{(\theta_{j}-\theta_{j+1})\cdots (\theta_j - \theta_{i})} $
\\
&&
\\
$0d^*0^*d$ 
&
$\varphi_1 \cdots \varphi_{i}
\frac{\varepsilon^*_0 }{\varepsilon_d \varphi}
$
&
$(\theta^*_d-\theta^*_0)\cdots (\theta^*_d-\theta^*_{i-j-1})
\lbrack j,i-j,d-i \rbrack_q
$
\\
&&
\\
$d^*0d0^*$ 
&
$\varphi_d \cdots \varphi_{d-i+1}
\frac{\varepsilon_d }{\varepsilon^*_0}$
&
$(\theta_0-\theta_d)\cdots (\theta_0-\theta_{d-i+j+1})
\lbrack j,i-j,d-i \rbrack_q
$
\\
&&
\\
$0d^*d0^*$
&
$ \frac{1}{\varphi_d \cdots \varphi_{d-i+1}} 
\;\frac{\varepsilon^*_0 }{\varepsilon_d}
$
&
$\frac{1}{(\theta^*_{d-j}-\theta^*_{d})\cdots (\theta^*_{d-j} - \theta^*_{d-j+1})}\;
\frac{1}{(\theta^*_{d-j}-\theta^*_{d-j-1})\cdots (\theta^*_{d-j} - \theta^*_{d-i})} $
\\
&&
\\
\hline
&&
\\
$d0^*0d^*$
&
$ \frac{1}{\varphi_1 \cdots \varphi_i} 
\;\frac{\varepsilon^*_d \varphi }{\varepsilon_0}
$
&
$\frac{1}{(\theta^*_{j}-\theta^*_{0})\cdots (\theta^*_j - \theta^*_{j-1})}\;
\frac{1}{(\theta^*_{j}-\theta^*_{j+1})\cdots (\theta^*_j - \theta^*_{i})} $
\\
&&
\\
$0^*d0d^*$ 
&
$\varphi_1 \cdots \varphi_{i}
\frac{\varepsilon_0}{\varepsilon^*_d \varphi}
$
&
$(\theta_d-\theta_0)\cdots (\theta_d-\theta_{i-j-1})
\lbrack j,i-j,d-i \rbrack_q
$
\\
&&
\\
$d0^*d^*0$ 
&
$\varphi_d \cdots \varphi_{d-i+1}
\frac{\varepsilon^*_d}{\varepsilon_0}
$
&
$(\theta^*_0-\theta^*_d)\cdots (\theta^*_0-\theta^*_{d-i+j+1})
\lbrack j,i-j,d-i \rbrack_q
$
\\
&&
\\
$0^*dd^*0$
&
$ \frac{1}{\varphi_d \cdots \varphi_{d-i+1}} 
\;\frac{\varepsilon_0}{\varepsilon^*_d}
$
&
$\frac{1}{(\theta_{d-j}-\theta_{d})\cdots (\theta_{d-j} - \theta_{d-j+1})}\;
\frac{1}{(\theta_{d-j}-\theta_{d-j-1})\cdots (\theta_{d-j} - \theta_{d-i})} $
\\
&&
\\
\hline
&&
\\
$dd^*00^*$
&
$ \frac{1}{\phi_d \cdots \phi_{d-i+1}} 
\;\frac{\varepsilon^*_0 \phi}{\varepsilon_0}
$
&
$\frac{1}{(\theta^*_{d-j}-\theta^*_{d})\cdots (\theta^*_{d-j} - \theta^*_{d-j+1})}\;
\frac{1}{(\theta^*_{d-j}-\theta^*_{d-j-1})\cdots (\theta^*_{d-j} - \theta^*_{d-i})} $
\\
&&
\\
$d^*d00^*$ 
&
$\phi_d \cdots \phi_{d-i+1}
\frac{\varepsilon_0 }{\varepsilon^*_0 \phi}
$
&
$(\theta_d-\theta_0)\cdots (\theta_d-\theta_{i-j-1})
\lbrack j,i-j,d-i \rbrack_q
$
\\
&&
\\
$dd^*0^*0$ 
&
$\phi_1 \cdots \phi_{i}
\frac{\varepsilon^*_0 }{\varepsilon_0 }
$
&
$(\theta^*_d-\theta^*_0)\cdots (\theta^*_d-\theta^*_{i-j-1})
\lbrack j,i-j,d-i \rbrack_q
$
\\
&&
\\
$d^*d0^*0$
&
$ \frac{1}{\phi_1 \cdots \phi_{i}} 
\;\frac{\varepsilon_0 }{\varepsilon^*_0 }
$
&
$\frac{1}{(\theta_{d-j}-\theta_{d})\cdots (\theta_{d-j} - \theta_{d-j+1})}\;
\frac{1}{(\theta_{d-j}-\theta_{d-j-1})\cdots (\theta_{d-j} - \theta_{d-i})} $
\\
&&
\\
\hline
&&
\\
$00^*dd^*$
&
$ \frac{1}{\phi_1 \cdots \phi_{i}} 
\;\frac{\varepsilon^*_d }{\varepsilon_d }
$
&
$\frac{1}{(\theta^*_{j}-\theta^*_{0})\cdots (\theta^*_{j} - \theta^*_{j-1})}\;
\frac{1}{(\theta^*_{j}-\theta^*_{j+1})\cdots (\theta^*_{j} - \theta^*_{i})} $
\\
&&
\\
$0^*0dd^*$ 
&
$\phi_1 \cdots \phi_{i}
\frac{\varepsilon_d }{\varepsilon^*_d }
$
&
$(\theta_0-\theta_d)\cdots (\theta_0-\theta_{d-i+j+1})
\lbrack j,i-j,d-i \rbrack_q
$
\\
&&
\\
$00^*d^*d$ 
&
$\phi_d \cdots \phi_{d-i+1}
\frac{\varepsilon^*_d }{\varepsilon_d \phi }
$
&
$(\theta^*_0-\theta^*_d)\cdots (\theta^*_0-\theta^*_{d-i+j+1})
\lbrack j,i-j,d-i \rbrack_q
$
\\
&&
\\
$0^*0d^*d$
&
$ \frac{1}{\phi_d \cdots \phi_{d-i+1}}
\;\frac{\varepsilon_d \phi}{\varepsilon^*_d}
$
&
$\frac{1}{(\theta_{j}-\theta_{0})\cdots (\theta_{j} - \theta_{j-1})}\;
\frac{1}{(\theta_{j}-\theta_{j+1})\cdots (\theta_{j} - \theta_{i})} $
\end{tabular}}

\medskip
\noindent In the above table, $q$ denotes a scalar 
in the algebraic closure of  $\K$ such that 
$q+q^{-1}+1$ is the common value of
(\ref{eq:betaplusone}).
\bigskip

\centerline{
\begin{tabular}[t]{c|c|c}
$wxyz$ &  $\lbrack wxyz \rbrack \rightarrow \lbrack xwyz \rbrack $ & 
 $\lbrack wxyz \rbrack \rightarrow \lbrack wyxz \rbrack $ \\
 & $ii$ entry & $ij$ entry $(i\geq j)$
\\
\hline
\hline
&&
\\
$d^*0^*0d$
& 
$\frac{\phi_d \cdots \phi_{d-i+1} }{\varphi_1 \cdots \varphi_i} 
\;\frac{\varepsilon^*_d \varphi}{\varepsilon^*_0 \phi}
$ 
&
$(\theta_i-\theta_0)\cdots (\theta_i-\theta_{j-1})$
\\
&&
\\
$0^*d^*0d$
& 
$\frac{\varphi_1 \cdots \varphi_i }{\phi_d \cdots \phi_{d-i+1}}  
\;\frac{\varepsilon^*_0 \phi}{\varepsilon^*_d \varphi}
$
&
$(\theta_i-\theta_0)\cdots (\theta_i-\theta_{j-1})$
\\
&&
\\
$d^*0^*d0$
& 
$\frac{\varphi_d \cdots \varphi_{d-i+1} }{\phi_1 \cdots \phi_{i}}  
\;\frac{\varepsilon^*_d }{\varepsilon^*_0 }
$
&
$(\theta_{d-i}-\theta_d)\cdots (\theta_{d-i}-\theta_{d-j+1})$
\\
&&
\\
$0^*d^*d0$
& 
$\frac{\phi_1 \cdots \phi_{i} }{\varphi_d \cdots \varphi_{d-i+1}}  
\;\frac{\varepsilon^*_0 }{\varepsilon^*_d }
$
&
$(\theta_{d-i}-\theta_d)\cdots (\theta_{d-i}-\theta_{d-j+1})$
\\
&&
\\
\hline
&&
\\
$d00^*d^*$
& 
$\frac{\phi_1 \cdots \phi_i }{\varphi_1 \cdots \varphi_i} 
\;\frac{\varepsilon_d \varphi }{\varepsilon_0 }
$
&
$(\theta^*_i-\theta^*_0)\cdots (\theta^*_i-\theta^*_{j-1})$
\\
&&
\\
$0d0^*d^*$
& 
$\frac{\varphi_1 \cdots \varphi_i }{\phi_1 \cdots \phi_i}  
\;\frac{\varepsilon_0 }{\varepsilon_d  \varphi}
$
&
$(\theta^*_i-\theta^*_0)\cdots (\theta^*_i-\theta^*_{j-1})$
\\
&&
\\
$d0d^*0^*$
& 
$\frac{\varphi_d \cdots \varphi_{d-i+1} }{\phi_d \cdots \phi_{d-i+1}} 
\;\frac{\varepsilon_d \phi }{\varepsilon_0 }
$
&
$(\theta^*_{d-i}-\theta^*_d)\cdots (\theta^*_{d-i}-\theta^*_{d-j+1})$
\\
&&
\\
$0dd^*0^*$
& 
$\frac{\phi_d \cdots \phi_{d-i+1} }{\varphi_d \cdots \varphi_{d-i+1}} 
\;\frac{\varepsilon_0 }{\varepsilon_d \phi }
$
&
$(\theta^*_{d-i}-\theta^*_d)\cdots (\theta^*_{d-i}-\theta^*_{d-j+1})$
\end{tabular}}

\end{theorem}
\begin{proof}
The basis 
 $\lbrack wxzy \rbrack  $, which is 
on the right in
(\ref{eq:threetargets}),
is the inversion of 
$\lbrack wxyz \rbrack  $  by 
Lemma \ref{lem:transshape}(iii).
Apparently  $Z$ is the transition matrix  from
 $ \lbrack wxyz  \rbrack $ to
 $ \lbrack wxzy \rbrack $.
We now consider the other two bases in
(\ref{eq:threetargets}).
For these
we prove our assertions case by case.
We begin with the  first row of the first table,
where $wxyz$ equals $d^*00^*d$. 
We consider the transition matrix from 
 $ \lbrack d^*00^*d \rbrack $   to
 $ \lbrack 0d^*0^*d \rbrack $.
We denote this matrix by $T$ and 
let $D$ denote the diagonal matrix in $\hbox{Mat}_{d+1}(\K)$
with $ii^{\hbox{th}}$ entry
\begin{equation}
D_{ii} = \frac{1}{\varphi_1  \varphi_2 \cdots \varphi_i} 
\;\frac{\varepsilon_d \varphi}{\varepsilon^*_0}
\qquad \qquad (0 \leq i \leq d).
\label{eq:ddef}
\end{equation}
We show $D=T$. 
Recall $\varepsilon_d\not=0$ by Lemma
\ref{lem:newvi}, and $\varphi \not=0$
by Definition
\ref{def:phinot}, so
 $D\not=0$.
Using the data  in the
first table in Theorem
\ref{thm:repa}, rows 1 and 2, we routinely find 
$ A^g D = DA^h $ and $ A^{*g} D = D A^{*h}$,
where we abbreviate $g$ for 
 $  d^*00^*d  $  and $h$ for 
 $  0d^*0^*d $. 
Applying
Lemma
\ref{lem:whentrans}, we find 
$D $ is an 
intertwining matrix from
 $ \lbrack d^*00^*d \rbrack $   to
 $ \lbrack 0d^*0^*d \rbrack $.
Therefore $D$ is a scalar multiple of $T$.  We show
this scalar is $1$.
To do this, we compare the $dd^{\hbox{th}}$ entry
of $D$ and $T$.
Setting $i=d$ in 
(\ref{eq:ddef}), and recalling 
$\varphi=\varphi_1 \varphi_2\cdots \varphi_d$, we find
 the
 $dd^{\hbox{th}}$ entry of $D$ equals 
$\varepsilon_d /\varepsilon^*_0 $.
We now 
 find the 
 $dd^{\hbox{th}}$ entry of $T$.
From the table in 
Theorem
\ref{thm:bases}, row $1$, we find the 
  $d^{\hbox{th}}$ vector in the
 basis 
 $ \lbrack d^*00^*d \rbrack $ is ${\tilde E}_d \eta^*_0$.
From the same table, row $2$,
  we find the $d^{\hbox{th}}$ vector in the
 basis 
 $ \lbrack 0d^*0^*d \rbrack $ is $ \eta_d$.
From the equation on the left in
(\ref{eq:ed}), 
 we find  
  $\eta_d= 
\varepsilon_d /
  \varepsilon^*_0 
  {\tilde E}_d \eta^*_0$,
 and it follows
the  $dd^{\hbox{th}}$ entry of $T$
is 
$\varepsilon_d /\varepsilon^*_0 $. 
 We now see $D$ and  $T$ have the same 
$dd^{\hbox{th}}$ entry, so $D=T$. In particular,
$D$ is the transition matrix from
 $ \lbrack d^*00^*d \rbrack $   to
 $ \lbrack 0d^*0^*d \rbrack $.

\medskip
\noindent 
We now consider the transition matrix from 
 $ \lbrack d^*00^*d \rbrack $   to
 $ \lbrack d^*0^*0d \rbrack $.
  We found the transition
 matrix from 
 $ \lbrack d^*0^*0d \rbrack $ to
 $ \lbrack d^*00^*d \rbrack $
  in the proof of 
 Theorem
\ref{thm:repa}.
To summarize, let $L$  
 denote the matrix in $\hbox{Mat}_{d+1}(\K)$
with $ij^{\hbox{th}}$ entry 
\begin{eqnarray}
L_{ij} = (\theta_i-\theta_0)(\theta_i - \theta_1)\cdots
(\theta_i-\theta_{j-1})
\qquad 
(0 \leq i,j\leq d).
\label{eq:ldefrem}
\end{eqnarray}
Then $L$
is the transition matrix from 
 $ \lbrack d^*0^*0d \rbrack $ to
 $ \lbrack d^*00^*d \rbrack $.
To get the transition matrix from 
 $ \lbrack d^*00^*d \rbrack $ to
 $ \lbrack d^*0^*0d \rbrack $,
 we find the inverse
of $L$. Observe $L$ is lower triangular.  Let $K$ denote the 
lower triangular matrix in 
  $\hbox{Mat}_{d+1}(\K)$ with $ij^{\hbox{th}}$ entry
\begin{equation}
K_{ij} = \frac{1}{(\theta_j-\theta_0)\cdots (\theta_j-\theta_{j-1})}
\; \frac{1}{(\theta_j-\theta_{j+1})\cdots (\theta_j-\theta_{i})}
\label{eq:kent}
\end{equation}
for $0 \leq j\leq i\leq d$. 
We recall $\theta_0, \theta_1,\ldots, \theta_d$ are mutually distinct,
so the denominator in 
(\ref{eq:kent}) is nonzero.
We claim $K$ is the inverse of $L$. To prove this, 
we show 
$LK = I$.
The matrices $L$ and $K$ are both lower triangular, so $LK$
is lower triangular. By
(\ref{eq:ldefrem}),
(\ref{eq:kent})
we find that for 
$0 \leq i \leq d$,
\beast
K_{ii} &=& \frac{1}{(\theta_i-\theta_0)\cdots (\theta_i-\theta_{i-1})}
\\
 &=& L_{ii}^{-1}
\eeast
so $(LK)_{ii} = 1$.
We  now show $(LK)_{ij}=0$ for $0 \leq j< i \leq d$.
Let $i,j$ be given.  
It suffices to show $(\theta_i-\theta_j)(LK)_{ij} = 0$, since
$\theta_0, \theta_1, \ldots, \theta_d$ are mutually distinct.
Observe
\begin{eqnarray}
(\theta_i-\theta_j)(LK)_{ij}
&=& (\theta_i-\theta_j)\sum_{h=0}^d L_{ih}K_{hj}    \nonumber
\\
&=& (\theta_i-\theta_j)\sum_{h=j}^i L_{ih}K_{hj}     \nonumber
\\
&=& \sum_{h=j}^i L_{ih}K_{hj}(\theta_i-\theta_h+\theta_h-\theta_j)
\nonumber
\\
&=& \sum_{h=j}^{i-1} L_{ih}(\theta_i-\theta_h)K_{hj}- 
\sum_{h=j+1}^i L_{ih}K_{hj}(\theta_j-\theta_h)
\nonumber
\\
&=& \sum_{h=j}^{i-1} L_{i,h+1}K_{hj}  -
\sum_{h=j+1}^i L_{ih}K_{h-1,j}
\label{eq:twosum}
\\
&=& 0
\nonumber
\end{eqnarray}
since the two sums in 
(\ref{eq:twosum}) are one and the same.
We have now shown $(LK)_{ij} = 0 $ for $0 \leq j<i\leq d$.
Combining our above arguments, we find $LK=I$ so
 $K$ is the inverse of $L$.
Now apparently $K$ is the transition matrix from 
 $ \lbrack d^*00^*d \rbrack $ to
 $ \lbrack d^*0^*0d \rbrack $.

\medskip
\noindent
We have now proved our assertions concerning the first row
of the first table. Applying these assertions to
the relatives of $\Phi$, and using
both
Theorem 
\ref{thm:phimod} and
Note
\ref{note:nosym},
we obtain our assertions concerning the
first and fourth
rows of each block of the first table.

\medskip
\noindent
We now consider the second row of the first table,
where $wxyz$ equals $0d^*0^*d$.
We find the transition matrix from 
$\lbrack 0d^*0^*d \rbrack$ to  
$\lbrack d^*00^*d \rbrack$.
Referring to the diagonal matrix
 $D$
from
(\ref{eq:ddef}) we showed
$D $ is the transition matrix from 
$\lbrack d^*00^*d \rbrack$ to
$\lbrack 0d^*0^*d \rbrack$. Therefore $D^{-1}$ 
is the transition matrix from 
$\lbrack 0d^*0^*d \rbrack$ to 
$\lbrack d^*00^*d \rbrack$.

\medskip
\noindent
We now consider the transition matrix from 
$\lbrack 0d^*0^*d \rbrack$ to  
$\lbrack 00^*d^*d \rbrack$.
Let $q$ denote a nonzero scalar in $\tilde  \K$ such that
$q+q^{-1}+1$ is the common value of
(\ref{eq:betaplusone}).
Let $H$ denote the lower triangular matrix
  in $\hbox{Mat}_{d+1}(\K)$  
with $ij^{\hbox{th}}$ entry  
\begin{equation}
 H_{ij} = (\theta^*_d-\theta^*_0)(\theta^*_d-\theta^*_1)\cdots
(\theta^*_d-\theta^*_{i-j-1})\lbrack j,i-j,d-i\rbrack_q
\label{eq:tentry}
\end{equation}
for $0 \leq j\leq i \leq d$.
The expression 
$\lbrack j,i-j,d-i\rbrack_q $ is given
in 
(\ref{eq:triplemain}). We remark each of $\lbrack 1 \rbrack_q,
\lbrack 2 \rbrack_q, \ldots , \lbrack d \rbrack_q$ is nonzero
by Corollary 
\ref{lem:dennonz},
so the denominator in 
$\lbrack j,i-j,d-i\rbrack_q $
 is nonzero.
We show $H$ is the  transition matrix from
$\lbrack 0d^*0^*d \rbrack$ to  
$\lbrack 00^*d^*d \rbrack$.
Observe $H_{ii} = 1 $ for $0 \leq i \leq d$, so $H$ is invertible.
 We show
 $A^{*g}H=HA^{*h}$, where we abbreviate
$g$ for $0d^*0^*d $ and  
$h$ for  $00^*d^*d $.
 The entries of
$A^{*g}$ and $A^{*h}$ are given in  the first table of
Theorem
\ref{thm:repa}, rows 2 and 15. Using this information, we find
that for $0 \leq i,j\leq d$, the $ij^{th}$
entry of $A^{*g}{H}$ is given by
\begin{equation}
\theta^*_i{H}_{ij}+{H}_{i+1,j},
\label{eq:ent1}
\end{equation}
where we interpret 
${H}_{i+1,j}=0$ if $i=d$. Similarily,
the 
$ij^{th}$
entry of ${H}A^{*h}$ is given by
\begin{equation}
{H}_{i,j-1}+
\theta^*_{d-j}{H}_{ij},
\label{eq:ent2}
\end{equation}
where we interpret 
${H}_{i,j-1}=0$ if $j=0$. 
We show 
(\ref{eq:ent1}) equals 
(\ref{eq:ent2}) or in other words
\begin{equation}
(\theta^*_i-\theta^*_{d-j}){H}_{ij} =
{H}_{i,j-1} - {H}_{i+1,j} .
\label{eq:needed}
\end{equation}
To prove
(\ref{eq:needed}),
first suppose $j-i>1$. Then each of
${H}_{ij}$,
${H}_{i,j-1}$, ${H}_{i+1,j}$ is zero
since $H$ is lower triangular, so
both sides of (\ref{eq:needed}) are zero.
 Next suppose $j-i=1$. Then 
${H}_{ij}=0$ since
$H$ is lower triangular. Moreover
${H}_{i,j-1}=H_{ii} =1$  and 
${H}_{i+1,j}=H_{jj} = 1$, 
 so both sides of 
(\ref{eq:needed}) are zero.
Next suppose $i=d$ and $j=0$. Then both sides
of 
(\ref{eq:needed}) are zero. 
Next suppose $i=d$ and $1\leq j\leq d$.
 Then using 
(\ref{eq:tentry}) we find both sides of 
(\ref{eq:needed}) equal 
$(\theta^*_d-\theta^*_0)(\theta^*_d-\theta^*_1)
\cdots (\theta^*_d -\theta^*_{d-j})$.
Next suppose $0 \leq i<d$ and $j=0$. Then using 
(\ref{eq:tentry}) we find both sides of 
(\ref{eq:needed}) equal the opposite of
$(\theta^*_d-\theta^*_0)(\theta^*_d-\theta^*_1)
\cdots (\theta^*_d -\theta^*_{i})$.
Finally suppose 
$1 \leq j\leq i\leq d-1$.
To verify 
(\ref{eq:needed}) in this case,
we use Lemma
\ref{lem:ident}. 
Set $r=j$, $s=i-j+1$, $t=d-i$, and
observe
each of $r,s,t$ is positive.
Since $r+s+t=d+1$,
and since each of
 $\lbrack 1 \rbrack_q,
\lbrack 2 \rbrack_q, \ldots , \lbrack d \rbrack_q$ is nonzero,
we find
$\lbrack h \rbrack_q \not=0$ for $1 \leq h <r+s+t$.
Apparently
our choice of $r,s,t $ satisfy 
the conditions of 
Lemma 
\ref{lem:ident}.  
Applying that lemma
we find
\begin{eqnarray}
&&\frac{\lbrack i-d+j\rbrack_q}{\lbrack d-i+j\rbrack_q}
\lbrack j,i-j,d-i\rbrack_q \nonumber \\ 
&& \qquad \qquad \qquad \qquad  =\;\; \lbrack j-1,i-j+1,d-i\rbrack_q - 
\lbrack j,i-j+1,d-i-1\rbrack_q. \qquad \qquad 
\label{eq:neededz}
\end{eqnarray}
Applying Lemma
\ref{lem:eigformnv}
 to the sequence $\theta^*_0, \theta^*_1,\ldots,\theta^*_d$,
and recalling each of 
 $\lbrack 1 \rbrack_q,  
 \lbrack 2 \rbrack_q,  
\ldots ,
 \lbrack d \rbrack_q  $ is  nonzero, 
we find
\begin{eqnarray}
\frac{\theta^*_i-\theta^*_{d-j}}{\theta^*_d-\theta^*_{i-j}}
=\frac{\lbrack i-d+j\rbrack_q}{\lbrack d-i+j\rbrack_q}.
\label{eq:neededt2}
\end{eqnarray}
Combining
(\ref{eq:neededz}),
(\ref{eq:neededt2}) we obtain 
\begin{eqnarray}
\frac{\theta^*_i-\theta^*_{d-j}}{\theta^*_d-\theta^*_{i-j}}
\lbrack j,i-j,d-i\rbrack_q  
 = \lbrack j-1,i-j+1,d-i\rbrack_q - 
\lbrack j,i-j+1,d-i-1\rbrack_q.
\label{eq:neededt}
\end{eqnarray}
Multiplying both sides of 
(\ref{eq:neededt}) by 
$(\theta^*_d-\theta^*_0)(\theta^*_d-\theta^*_1)\cdots
(\theta^*_d-\theta^*_{i-j})$, and evaluating the result using 
(\ref{eq:tentry}),
we routinely obtain 
(\ref{eq:needed}).
We have now shown
(\ref{eq:needed}) holds for $0 \leq i,j\leq d$,
and it follows
$ A^{*g}{H}={H}A^{*h}$.
Recall we are trying to show $H$ is the transition matrix
from 
$\lbrack 0d^*0^*d \rbrack$ to   
$\lbrack 00^*d^*d \rbrack$.
 Let $N$ denote this 
 transition matrix.
To  
show $H=N$,
we proceed in two steps. We first show
$H$ is a scalar multiple of $N$. We then show
this scalar equals $1$.
Proceeding with the first step,
 we define  
$S:={N}{H}^{-1}$ and show $S$ is a scalar multiple
of the identity. 
By Lemma
\ref{lem:transshape}(ii),
we find $N$ is lower triangular.
Recall $H$ is lower triangular, so $S$ is lower triangular.
Since $N$ is
the transition  matrix from
$\lbrack 0d^*0^*d \rbrack$ to  
$\lbrack 00^*d^*d \rbrack$ we find 
$N$ is an 
 intertwining matrix from
$\lbrack 0d^*0^*d \rbrack$ to  
$\lbrack 00^*d^*d \rbrack$. Therefore 
$A^{*g}{N}={N}A^{*h}$.
Combining this with 
$ A^{*g}{H}={ H}A^{*h}$, we  find
$SA^{*g}=A^{*g}S$.
We claim  $S$ is diagonal.
Suppose not. Then there exists a pair of integers $i,j$ $(0 \leq j<i\leq d)$
such that $S_{ij}\not=0$. Of all such pairs $i,j$ pick one with
$i-j$ maximal.
We compute the $ij^{\hbox{th}}$ entry in 
$SA^{*g}=A^{*g}S$.
Observe the $ij^{\hbox{th}}$ entry of $SA^{*g}$ is
$S_{ij}\theta^*_j$ and that of $A^{*g}S$ is $\theta^*_iS_{ij}$,
so $(\theta^*_i-\theta^*_j)S_{ij}=0$.
Observe $\theta^*_i\not=\theta^*_j$,
so $S_{ij}=0$, a contradiction. We have now shown $S$ is diagonal.
Computing entries just above the main diagonal in
 $SA^{*g}=A^{*g}S$, we find $S$ is a scalar multiple of the identity.
Apparently
 $H$ is a scalar  multiple of $N$.
We now show this scalar equals $1$. To do this,
we compare the 
 $dd^{\hbox{th}}$ entry of 
$H$ and $N$.
We saw above that the 
 $dd^{\hbox{th}}$ entry of $H$ equals $1$.
We find the
 $dd^{\hbox{th}}$ entry of 
 $N$.
From the table in Theorem
\ref{thm:bases}, row 2, we find the 
  $d^{\hbox{th}}$ vector in the
 basis 
 $ \lbrack 0d^*0^*d \rbrack $ is  $\eta_d$.
From the same table, row 15,
 we find the 
  $d^{\hbox{th}}$ vector in the
 basis 
 $ \lbrack 00^*d^*d \rbrack $ is  $\eta_d$.
Apparently 
 the $dd^{\hbox{th}}$ entry of $N$ equals $1$.
We now see $H$ and  $N$ have the same
$dd^{\hbox{th}}$ entry,
so $H=N$.
	In particular
 $H$ is the transition matrix from
$\lbrack 0d^*0^*d \rbrack$ to 
$\lbrack 00^*d^*d \rbrack$.

\medskip
\noindent
We have now proved our assertions concerning the second row
of the first table. Applying these assertions to
the relatives of $\Phi$,
and using both
Theorem
\ref{thm:phimod} and
Note \ref{note:nosym},
we obtain our assertions concerning the second and third rows 
of each block 
of the first table. We have now verified all our assertions concerning
the first table.

\medskip
\noindent
Consider the first row of the second table, where 
$wxyz $ equals $d^*0^*0d$. We find the transition matrix
from 
 $\lbrack d^*0^*0d \rbrack $ to
 $\lbrack 0^*d^*0d \rbrack $. Let $P$ denote this matrix
 and let $F$ denote the diagonal
 matrix in $\hbox{Mat}_{d+1}(\K)$ with diagonal entries
 \begin{equation}
 F_{ii} = \frac{\phi_d \phi_{d-1} \cdots \phi_{d-i+1}}
 {\varphi_1 \varphi_2 \cdots \varphi_i}\; 
 \frac{\varepsilon^*_d\varphi}{\varepsilon^*_0 \phi} 
 \qquad \qquad (0 \leq i \leq d).
\label{eq:fdef}
\end{equation}
We show $F=P$. 
To do this, we first show
$F$ is
an intertwining matrix
from 
 $\lbrack d^*0^*0d \rbrack $ to
 $\lbrack 0^*d^*0d \rbrack $.  
Clearly $F\not=0$. We 
show $A^gF = FA^h$, $A^{*g}F = FA^{*h}$, where
we abbreviate $g$ for 
 $ d^*0^*0d$ and  $h$ for  
 $ 0^*d^*0d$. 
The matrices representing 
$A$ and $A^*$ with respect to 
 $\lbrack d^*0^*0d \rbrack $ and 
 $\lbrack 0^*d^*0d \rbrack $
are given in  the 
second table of Theorem \ref{thm:repa}, rows 1 and 2.
Using the
data in these rows, we routinely find
$A^gF = FA^h$, $A^{*g}F = FA^{*h}$.
  Applying Lemma
\ref{lem:whentrans}, we find $F$ is an intertwining matrix
from 
 $\lbrack d^*0^*0d \rbrack $ to
 $\lbrack 0^*d^*0d \rbrack $.
Now apparently $F$ is a scalar multiple of $P$.
We show this scalar equals $1$.
To do this, we compare the 
 $dd^{\hbox{th}}$ entry
of $F$ and $P$.
Setting $i=d$ in 
(\ref{eq:fdef}),
and recalling $\varphi=\varphi_1\varphi_2\cdots \varphi_d$,
 $\phi=\phi_1\phi_2\cdots \phi_d$,
we find 
 the $dd^{\hbox{th}}$ entry of $F$ equals 
 $\varepsilon^*_d/\varepsilon^*_0$.
We now find the 
  $dd^{\hbox{th}}$ entry of $P$. 
From the table in 
Theorem \ref{thm:bases}, row 17, we find the 
  $d^{\hbox{th}}$ vector in the
 basis 
 $ \lbrack d^*0^*0d \rbrack $ is  $E_d\eta^*_0$.
From the same table, row 18,
 we find the 
  $d^{\hbox{th}}$ vector in the
 basis 
 $ \lbrack 0^*d^*0d \rbrack $ is  $E_d\eta^*_d$.
From
the equation on the left in 
(\ref{eq:dds}), we find 
 $ E_d\eta^*_d =
 \varepsilon^*_d/\varepsilon^*_0
  E_d\eta^*_0 $.
Apparently 
 the $dd^{\hbox{th}}$ entry of $P$ equals
 $\varepsilon^*_d/\varepsilon^*_0 $.
We now see $F$ and $P$ have the same  
 $dd^{\hbox{th}}$ entry, so  $F=P$.
In particular, $F$ is the 
transition matrix from
 $\lbrack d^*0^*0d \rbrack $ to
 $\lbrack 0^*d^*0d \rbrack $.

\medskip
\noindent
We  already found the transition matrix from
 $\lbrack d^*0^*0d \rbrack $ to
 $\lbrack d^*00^*d \rbrack $. This is the matrix $L$ from
(\ref{eq:ldefrem}).

\medskip
\noindent
We have now obtained
our assertions concerning the first row of the second table.
Applying these assertions to the relatives of $\Phi$,
and using both
Theorem
\ref{thm:phimod} and
Note
 \ref{note:nosym},
we
obtain all our assertions concerning  the
second table. This completes
the proof.

\end{proof}
\noindent We finish this section with some comments on
transition matrices.
Let $\Phi $ denote the Leonard system in 
(\ref{eq:ourstartingpt}), and let $g, h$ denote elements
in $S_4$. 
Consider the transition
 matrix
from the basis $\lbrack g \rbrack $ to the basis
 $\lbrack h \rbrack $. If $g$ and $h$ are adjacent
 in the sense of
 Definition
\ref{def:S4interp},
 then this transition matrix is given in 
 Theorem
\ref{thm:trans}.
If the above restriction  on $g, h$ is removed, then
this transition matrix can be computed as follows.
To explain the idea, we use the following notation. 
By an {\it edge} in $S_4$, we mean an ordered pair 
consisting of adjacent elements of $S_4$.
Let $r$ denote a nonnegative integer. 
By a {\it walk of length $r$} in $S_4$,
we mean a sequence $g_0, g_1, \ldots, g_r$ of elements of $S_4$
such that $g_{i-1},g_i$ is an edge for $1 \leq i \leq r$.
The above walk is said to be {\it from } $g_0$ {\it to } $g_r$.
let $gh$ denote an edge in $S_4$. By the {\it weight}
of that  edge, we mean the transition matrix from $\lbrack g \rbrack
$ to $\lbrack h \rbrack$.
Let 
 $g_0, g_1, \ldots, g_r$ denote a walk in $S_4$.
By the {\it weight} of this  walk, we mean the product
$W_1W_2\cdots W_r$, where $W_i$ is the weight of the edge
$g_{i-1},g_i$ for $1 \leq i \leq r$.
Let $g,h$ denote elements in $S_4$. Then the transition matrix
from 
 $\lbrack g \rbrack
$ to $\lbrack h \rbrack$
is given by the weight of any walk from $g$ to $h$.

%
%
%


\section{Remarks}
\medskip
\noindent In the introduction to this paper, we mentioned  
that Leonard pairs are related to 
certain orthogonal polynomials contained in 
the Askey scheme. 
One significance of the polynomials
is that they give the entries in the transition
matrices relating certain pairs of bases among our set of 24.
In this section, we illustrate what is going on  with some examples.
For related work, see 
\cite{GYZnature},
 \cite{GYLZmut},
\cite{GYZlinear},
\cite{Zhidd} and 
\cite{Koelink3},
\cite{Koelink1},
\cite{Koelink2}, 
\cite{Koelink4},
\cite{koo3},
\cite[ch. 4]{Hjal}.

\medskip
\noindent Throughout this section, 
we let $\Phi$ denote the Leonard system in
(\ref{eq:ourstartingpt}),
with eigenvalue sequence
 $\theta_0, \theta_1,\ldots, \theta_d $,
dual eigenvalue sequence 
$\theta^*_0, \theta^*_1,\ldots, \theta^*_d $,
first split sequence 
$\varphi_1, \varphi_2,\ldots, \varphi_d $  
and 
second split sequence
$\phi_1, \phi_2,\ldots, \phi_d $.
For $0 \leq i,j\leq d$
we define
\begin{eqnarray}
{\cal P}_{ij}=
\sum_{n=0}^d
\frac{(\theta_i-\theta_0)(\theta_i-\theta_1)\cdots (\theta_i-\theta_{n-1})
(\theta^*_j-\theta^*_0)(\theta^*_j-\theta^*_1)
\cdots (\theta^*_j-\theta^*_{n-1})
}{\varphi_1 \varphi_2 \cdots \varphi_n}.
\label{eq:sumpart}
\end{eqnarray}
We observe ${\cal P}_{ij}$
is
 a polynomial of degree $ j$ in $\theta_i$ and a 
polynomial of degree $i$ in $\theta^*_j$.
These are the polynomials of interest.

\medskip
\noindent The ${\cal P}_{ij}$ arise in the following context.
Let
$V$ denote the irreducible left $\cal A$-module.
In Theorem 
\ref{thm:bases}, we presented 24 bases for $V$.
Of these, we focus on the following two:
\begin{eqnarray}
&&\lbrack d^*0^*0d \rbrack: \qquad \qquad \qquad 
E_0 \eta^*_0, E_1 \eta^*_0, \ldots, E_d \eta^*_0, \qquad \qquad \qquad 
\label{eq:stbasiscom}
\\
&&\lbrack d00^*d^* \rbrack:
\qquad \qquad \qquad  
E^*_0 \eta_0, E^*_1 \eta_0, \ldots, E^*_d \eta_0. \qquad \qquad \qquad
\label{eq:dstbasiscom}
\end{eqnarray}
 We recall the basis 
(\ref{eq:stbasiscom}) is a $\Phi$-standard basis. With respect
to this basis,
 the matrix representing $A$ is diagonal, and
the matrix representing $A^*$ is irreducible tridiagonal.
We denote these matrices by $H$ and $B^*$, respectively.
Their entries  are given in 
the second table
   of Theorem
\ref{thm:repa}, row 1.
The basis
(\ref{eq:dstbasiscom}) is a $\Phi^*$-standard basis. 
With respect to this basis,
the matrix 
 representing $A^*$ is  diagonal
 and
the matrix representing $A$ is  
 irreducible tridiagonal.
We denote these matrices by $H^*$ and $B$, respectively.
Their entries  are given in 
the third table
   of Theorem
\ref{thm:repa}, row 1.
Let $P$ denote the 
transition matrix from  
(\ref{eq:stbasiscom}) to
(\ref{eq:dstbasiscom}), with the  vectors
$\eta_0, \eta^*_0$ chosen so that
\begin{equation}
\eta^*_0=E^*_0\eta_0.
\label{eq:goodchoice}
\end{equation}
The effect of 
(\ref{eq:goodchoice}) is that $P_{i0}=1$ for $0 \leq i \leq d$.
We let $P^*$ 
 denote the 
transition matrix from  
(\ref{eq:dstbasiscom}) to 
(\ref{eq:stbasiscom}),
this time with the 
$\eta_0, \eta^*_0$ chosen so that
\begin{equation}
\eta_0=E_0\eta^*_0.
\label{eq:goodchoice2}
\end{equation}
As expected 
$P^*_{i0}=1$ for $0 \leq i \leq d$.
From the construction of $P$ and $P^*$ we find
there exists a nonzero scalar $\nu \in \K$ such that
\begin{eqnarray}
PP^*=\nu I.
\label{eq:ppn}
\end{eqnarray}
Moreover by Lemma 
\ref{lem:whentrans}
we have 
\begin{eqnarray}
B^*P=PH^*, \qquad \qquad BP^*=P^*H.
\label{eq:prel}
\end{eqnarray}
We compute the entries of $P$. For this
we use the method outlined in the last paragraph of the previous section.
The following is a walk in $S_4$ from
$d^*0^*0d $ to $ d00^*d^* $.
\begin{eqnarray}
d^*0^*0d, \quad d^*00^*d, \quad 0d^*0^*d, \quad 0d^*d0^*, \quad 0dd^*0^*,
\quad d0d^*0^*, \quad d00^*d^*.
\label{eq:walk1}
\end{eqnarray}
Apparently $P$ equals the weight of the walk 
(\ref{eq:walk1}). 
Computing this weight using the data in Theorem
\ref{thm:trans},
we find 
\begin{eqnarray}
P_{ij}=k_j{\cal P}_{ij} \qquad \qquad  
(0 \leq i,j\leq d),
\label{eq:pcalp}
\end{eqnarray}
where ${\cal P}_{ij}$ is from
(\ref{eq:sumpart}), and where
 $k_j$ equals
\begin{eqnarray}
\frac{\varphi_1 \varphi_2 \cdots \varphi_j}{\phi_1 \phi_2 \cdots \phi_j}
\label{eq:kjpart1}
\end{eqnarray}
times
\begin{eqnarray}
\frac{(\theta^*_0-\theta^*_1)(\theta^*_0-\theta^*_2)\cdots
(\theta^*_0-\theta^*_d)}
{(\theta^*_j-\theta^*_0)\cdots (\theta^*_j-\theta^*_{j-1})
(\theta^*_j-\theta^*_{j+1})\cdots (\theta^*_j-\theta^*_d)}
\label{eq:kjpart2}
\end{eqnarray}
for $0 \leq j \leq d$.
We now compute $P^*$. Replacing $\Phi$ by $\Phi^*$ in the above
discussion, and using
Theorem \ref{thm:phimod},
we routinely find
\begin{eqnarray}
P^*_{ij}=k^*_j{\cal P}_{ji} \qquad \qquad  
(0 \leq i,j\leq d),
\label{eq:pcalps}
\end{eqnarray}
where ${\cal P}_{ji}$ is from
(\ref{eq:sumpart}), and where
$k^*_j$ equals
\begin{eqnarray}
\frac{\varphi_1 \varphi_2 \cdots \varphi_j}{\phi_d \phi_{d-1} \cdots \phi_{d-j+1}}
\label{eq:kjspart1}
\end{eqnarray}
times
\begin{eqnarray}
\frac{(\theta_0-\theta_1)(\theta_0-\theta_2)\cdots
(\theta_0-\theta_d)}
{(\theta_j-\theta_0)\cdots (\theta_j-\theta_{j-1})
(\theta_j-\theta_{j+1})\cdots (\theta_j-\theta_d)}
\label{eq:kjspart2}
\end{eqnarray}
for $0 \leq j \leq d$.
We now compute the scalar $\nu$ from 
(\ref{eq:ppn}).
From the construction of $P$ and $P^*$ we routinely find
$\nu E_0E^*_0E_0=E_0$.  Taking the trace in this equation we find
\begin{eqnarray}
\hbox{trace}\,E_0E^*_0 = \nu^{-1}.
\label{eq:tracen}
\end{eqnarray}
Evaluating the left side in
(\ref{eq:tracen}) using
Lemma 
\ref{lem:etildee} and
Lemma \ref{lem:tractildee}, we routinely find
\begin{eqnarray}
\nu = 
\frac{
(\theta_0-\theta_1)(\theta_0-\theta_2)\cdots
(\theta_0-\theta_d)
(\theta^*_0-\theta^*_1)(\theta^*_0-\theta^*_2)\cdots
(\theta^*_0-\theta^*_d)}
{\phi_1 \phi_2 \cdots \phi_d}.
\label{eq:nclform}
\end{eqnarray}
From (\ref{eq:ppn}) we obtain the following orthogonality relations
for the ${\cal P}_{ij}$.
Expanding the left side of $PP^*=\nu I$ using matrix 
multiplication, and evaluating the result using
(\ref{eq:pcalp}),
(\ref{eq:pcalps}) we find
\begin{eqnarray}
\sum_{n=0}^d {\cal P}_{in}{\cal P}_{jn} k_n = \delta_{ij} \nu k^{*-1}_j
\qquad \qquad (0 \leq i,j \leq d).
\end{eqnarray}
Doing something similar with the equation
 $P^*P=\nu I$ we find 
\begin{eqnarray}
\sum_{n=0}^d {\cal P}_{ni}{\cal P}_{nj} k^*_n = \delta_{ij} \nu k^{-1}_j
\qquad \qquad (0 \leq i,j \leq d).
\end{eqnarray}
We remark the equations 
(\ref{eq:prel})  express several three-term recurrences satisfied
by the ${\cal P}_{ij}$.

\medskip
\noindent We now indicate how the ${\cal P}_{ij}$ fit into
the Askey scheme. 
Instead of  giving a complete treatment,
we content ourselves with two 
examples.

\medskip
\noindent Our first example
is associated with the Leonard pair  
from (\ref{eq:fam1}). For this example the ${\cal P}_{ij}$
will turn out to be Krawtchouk polynomials.
Let $d$ denote a nonnegative integer, and 
consider the following elements of $\K$.
\begin{eqnarray}
&&\theta_i =  d-2i, \qquad \qquad \theta^*_i = d-2i 
\qquad \qquad (0 \leq i \leq d),
\label{eq:thsol}
\\
&&
\varphi_i = -2i(d-i+1), \qquad \qquad \phi_i = 2i(d-i+1) 
\qquad  \qquad (1 \leq i \leq d).
\label{eq:vpsol}
\end{eqnarray}
To avoid degenerate situations, we assume the characteristic
of $\K$ is zero or an odd prime greater than $d$.
It is routine to show 
(\ref{eq:thsol}), 
(\ref{eq:vpsol}) satisfy the conditions (i)--(v) of 
Theorem 
\ref{thm:classls}.
Let us assume $\Phi$ is the 
 corresponding Leonard system
from that theorem.
For this $\Phi$, we routinely find
 $B$ and $B^*$ both equal the matrix
on the left in 
(\ref{eq:fam1}). Moreover $H$ and $H^*$ both equal the matrix
on the right in 
(\ref{eq:fam1}).
Pick any integers $i,j$ $(0 \leq i,j\leq d)$.
Evaluating the right side of
(\ref{eq:sumpart}) using
(\ref{eq:thsol}), 
(\ref{eq:vpsol}), we find 
 ${\cal P}_{ij}$ equals 
\begin{eqnarray}
\sum_{n=0}^d \frac{(-i)_n (-j)_n 2^n}{(-d)_n n! }, 
\label{eq:2F1expand}
\end{eqnarray}
 where
\beast
(a)_n:=a(a+1)(a+2)\cdots (a+n-1) \qquad \qquad n=0,1,2,\ldots
\eeast
Hypergeometric series are defined in \cite[p. 3]{gasperrahmanbk}.
From this definition we find 
(\ref{eq:2F1expand}) is the hypergeometric series
\begin{eqnarray}
{{}_2}F_1\Biggl({{-i, -j}\atop {-d}}\;\Bigg\vert \;2\Biggr).
\label{eq:2F1not}
\end{eqnarray}
A definition of the Krawtchouk polynomials can be found in 
 \cite{AAR} or 
\cite{KoeSwa}. Comparing this definition with 
(\ref{eq:2F1not}),
we find ${\cal P}_{ij}$  
is a  
Krawtchouk polynomial of degree $j$ in $\theta_i$ and
a Krawtchouk polynomial 
of degree $i$ in $\theta^*_j$.
Pick an integer $j$ $(0 \leq j \leq d)$.
Evaluating
(\ref{eq:kjpart1}), 
(\ref{eq:kjpart2}) and 
(\ref{eq:kjspart1}),
(\ref{eq:kjspart2}) using
(\ref{eq:thsol}), 
(\ref{eq:vpsol}), we find 
$k_j$ and $k^*_j$ both equal the binomial coefficient 
\beast
\Biggl({{ d }\atop {j}}\Biggr). 
\eeast
Evaluating 
(\ref{eq:nclform}) using 
(\ref{eq:thsol}),
(\ref{eq:vpsol})
we find $\nu=2^d$. We comment that for this example
$P=P^*$, so $P^2=2^dI$.

\medskip
\noindent
We now give our second example. For this
example the 
 ${\cal P}_{ij}$ will turn out to be  
$q$-Racah polynomials. To begin, 
let  $d$ denote a nonnegative integer,   and consider the 
following elements in $\K$.
\begin{eqnarray}
\theta_i &=& \theta_0 + h(1-q^i)(1-sq^{i+1})/q^i,
\label{eq:thdefend}
\\
\theta^*_i &=& \theta^*_0 + h^*(1-q^i)(1-s^*q^{i+1})/q^i
\label{eq:thsdefend}
\end{eqnarray}
for $0 \leq i \leq d$, and
\begin{eqnarray}
\varphi_i &=& hh^*q^{1-2i}(1-q^i)(1-q^{i-d-1})(1-r_1q^i)(1-r_2q^i),
\label{eq:varphidefend}
\\
\phi_i &=& hh^*q^{1-2i}(1-q^i)(1-q^{i-d-1})(r_1-s^*q^i)(r_2-s^*q^i)/s^*
\label{eq:phidefend}
\end{eqnarray}
for $1 \leq i \leq d$. We assume 
$q, h, h^*, s, s^*, r_1, r_2$ are nonzero scalars
in the algebraic closure $\tilde K$, and that $r_1r_2 = s s^*q^{d+1}$.
It is routine to show 
(\ref{eq:thdefend})--(\ref{eq:phidefend})
give
a parametric solution to 
Theorem 
\ref{thm:classls}(iii)--(v).
Let us assume conditions (i),(ii) of 
Theorem \ref{thm:classls} are satisfied as well,
so that 
(\ref{eq:thdefend})--(\ref{eq:phidefend})
correspond to a Leonard system.
We assume $\Phi$ is the
corresponding Leonard system from 
Theorem \ref{thm:classls}.
For this $\Phi$ we find $B$, $B^*$, 
${\cal P}_{ij}$,
$k_j$, $k^*_j$,
 $\nu $.  
Recall the entries of $B$ are given in 
the third table of Theorem \ref{thm:repa}, row 1.
Evaluating these entries using
(\ref{eq:thdefend})--(\ref{eq:phidefend}),
 we find
\beast
B_{01} &=& \frac{h(1-q^{-d})(1-r_1q)(1-r_2q)}
{1-s^*q^2},
\\
B_{i-1,i} &=& \frac{h(1-q^{i-d-1})(1-s^*q^i)(1-r_1q^i)(1-r_2q^i)}
{(1-s^*q^{2i-1})(1-s^*q^{2i})} \qquad \quad (2 \leq i \leq d),
\\
B_{i,i-1} &=& \frac{h(1-q^i)(1-s^*q^{i+d+1})(r_1-s^*q^i)(r_2-s^*q^i)}
{s^*q^d(1-s^*q^{2i})(1-s^*q^{2i+1})} \qquad (1 \leq i \leq d-1),
\\
B_{d,d-1} &=& \frac{h(1-q^d)(r_1-s^*q^d)(r_2-s^*q^d)}
{s^*q^d(1-s^*q^{2d})},
\\
B_{ii} &=& \theta_0 - B_{i,i-1}-B_{i,i+1} \qquad (0 \leq i \leq d),
\eeast
where we define $B_{0,-1}:=0$, $B_{d,d+1}:=0$.
The entries
of $B^*$ are similarly obtained.
To get the  entries of $B^*$, in the above formulae 
exhange $(\theta_0,h,s) $ and $(\theta^*_0,h^*,s^*)$,
and preserve $(r_1, r_2,q )$.
Pick  integers $i,j$ $(0 \leq i,j\leq d)$.
Evaluating the right side of 
(\ref{eq:sumpart}) using
(\ref{eq:thdefend})--(\ref{eq:phidefend}),
we find  ${\cal P}_{ij}$ equals
\begin{eqnarray}
\sum_{n=0}^d \frac{(q^{-i};q)_n (sq^{i+1};q)_n 
(q^{-j};q)_n (s^*q^{j+1};q)_n q^n}
{(r_1q;q)_n(r_2q;q)_n (q^{-d};q)_n(q;q)_n},
\label{eq:uihyper}
\end{eqnarray}
where 
\beast
(a;q)_n := (1-a)(1-aq)(1-aq^2)\cdots (1-aq^{n-1})\qquad \qquad n=0,1,2\ldots 
\eeast
Basic hypergeometric series are defined in \cite[p. 4]{gasperrahmanbk}.
From that definition we find
(\ref{eq:uihyper}) is the basic hypergeometric series
\begin{eqnarray}
 {}_4\phi_3 \Biggl({{q^{-i}, \;sq^{i+1},\;q^{-j},\;s^*q^{j+1}}\atop
{r_1q,\;\;r_2q,\;\;q^{-d}}}\;\Bigg\vert \; q,\;q\Biggr).
\label{eq:qrac}
\end{eqnarray}
A definition of the  $q$-Racah polynomials can be found in 
\cite{Ask}, \cite{AskAW}, or 
\cite{KoeSwa}. Comparing this definition with 
(\ref{eq:qrac}),
 and
recalling $r_1r_2=s s^*q^{d+1}$,
we find  ${\cal P}_{ij}$ 
is a  
$q$-Racah polynomial of degree $j$ in $\theta_i$ and
a $q$-Racah polynomial 
of degree $i$ in $\theta^*_j$.
Pick an integer $j$ $(0 \leq j \leq d)$.
Evaluating  
(\ref{eq:kjpart1}),
(\ref{eq:kjpart2}) using 
(\ref{eq:thdefend})--(\ref{eq:phidefend}),
we find
\begin{eqnarray}
k_j = \frac{(r_1q;q)_j(r_2q;q)_j(q^{-d};q)_j(s^*q;q)_j(1-s^*q^{2j+1})}
{s^jq^j(q;q)_j(s^*q/r_1;q)_j(s^*q/r_2;q)_j(s^*q^{d+2};q)_j(1-s^*q)}. 
\label{eq:qrack}
\end{eqnarray}
The scalar $k^*_j$ is similarly found.
To get $k^*_j$, in 
(\ref{eq:qrack})
exchange $s$ and $s^*$, and preserve $(r_1, r_2,q )$.
Evaluating 
(\ref{eq:nclform}) using 
(\ref{eq:thdefend})--(\ref{eq:phidefend}),
 we find
\beast
\nu = \frac{(sq^2;q)_d (s^*q^2;q)_d}{r^d_1q^d(sq/r_1;q)_d(s^*q/r_1;q)_d}. 
\eeast
%

%

\begin{thebibliography}{10}

\bibitem{AAR}
George~E. Andrews, Richard Askey, and Ranjan Roy.
\newblock {\em Special functions}.
\newblock Cambridge University Press, Cambridge, 1999.

\bibitem{Ask}
R.~Askey and J.~Wilson.
\newblock A set of orthogonal polynomials that generalize the {R}acah
  coefficients or $6-j$\ symbols.
\newblock {\em SIAM J. Math. Anal.}, 10(5):1008--1016, 1979.

\bibitem{AskAW}
R.~Askey and J.~Wilson.
\newblock Some basic hypergeometric orthogonal polynomials that generalize
  {J}acobi polynomials.
\newblock {\em Mem. Amer. Math. Soc.}, 54(319):iv+55, 1985.

\bibitem{BanIto}
E.~Bannai and T.~Ito.
\newblock {\em Algebraic Combinatorics I: Association Schemes}.
\newblock Benjamin/Cummings, London, 1984.

\bibitem{bcn}
A.~E. Brouwer, A.~M. Cohen, and A.~Neumaier.
\newblock {\em Distance-Regular Graphs}.
\newblock Springer-Verlag, Berlin, 1989.

\bibitem{Cau}
J.~S. Caughman~{I}{V}.
\newblock The {T}erwilliger algebras of bipartite ${P}$- and ${Q}$-polynomial
  schemes.
\newblock {\em Discrete Math.}, 196(1-3):65--95, 1999.

\bibitem{CurNom}
B.~Curtin and K.~Nomura.
\newblock Distance-regular graphs related to the quantum enveloping algebra of
  $sl(2)$.
\newblock {\em J. Algebraic Combin.}, to appear.

\bibitem{Curspin}
Brian Curtin.
\newblock Distance-regular graphs which support a spin model are thin.
\newblock {\em Discrete Math.}, 197/198:205--216, 1999.
\newblock 16th British Combinatorial Conference (London, 1997).

\bibitem{CR}
Charles~W. Curtis and Irving Reiner.
\newblock {\em Methods of representation theory. {V}ol. {I}}.
\newblock John Wiley \& Sons Inc., New York, 1990.
\newblock With applications to finite groups and orders, Reprint of the 1981
  original, A Wiley-Interscience Publication.

\bibitem{gasperrahmanbk}
G.~Gasper and M.~Rahman.
\newblock {\em Basic hypergeometric series}, volume~35 of {\em Encyclopedia of
  Mathematics and its Applications}.
\newblock Cambridge University Press, Cambridge, 1990.

\bibitem{go}
J.~Go.
\newblock The {T}erwilliger algebra of the {H}ypercube ${Q_D}$.
\newblock {\em European J. Combin.}, to appear.

\bibitem{GYZnature}
Ya.~A. Granovski{\u\i} and A.~S. Zhedanov.
\newblock Nature of the symmetry group of the $6j$-symbol.
\newblock {\em Zh. \`Eksper. Teoret. Fiz.}, 94(10):49--54, 1988.

\bibitem{GYLZmut}
Ya.~I. Granovski{\u\i}, I.~M. Lutzenko, and A.~S. Zhedanov.
\newblock Mutual integrability, quadratic algebras, and dynamical symmetry.
\newblock {\em Ann. Physics}, 217(1):1--20, 1992.

\bibitem{GYZTwisted}
Ya.~I. Granovski{\u\i} and A.~S. Zhedanov.
\newblock ``{T}wisted'' {C}lebsch-{G}ordan coefficients for ${\rm {s}{u}}\sb
  q(2)$.
\newblock {\em J. Phys. A}, 25(17):L1029--L1032, 1992.

\bibitem{GYZlinear}
Ya.~I. Granovski{\u\i} and A.~S. Zhedanov.
\newblock Linear covariance algebra for ${\rm {s}{l}}\sb q(2)$.
\newblock {\em J. Phys. A}, 26(7):L357--L359, 1993.

\bibitem{GYZspherical}
Ya.~I. Granovski{\u\i} and A.~S. Zhedanov.
\newblock Spherical $q$-functions.
\newblock {\em J. Phys. A}, 26(17):4331--4338, 1993.

\bibitem{GH4}
F.~Alberto Gr{\"u}nbaum.
\newblock Some bispectral musings.
\newblock In {\em The bispectral problem (Montreal, PQ, 1997)}, pages 31--45.
  Amer. Math. Soc., Providence, RI, 1998.

\bibitem{GH5}
F.~Alberto Gr{\"u}nbaum and Luc Haine.
\newblock Bispectral {D}arboux transformations: an extension of the {K}rall
  polynomials.
\newblock {\em Internat. Math. Res. Notices}, 1997(8):359--392.

\bibitem{GH7}
F.~Alberto Gr{\"u}nbaum and Luc Haine.
\newblock The $q$-version of a theorem of {B}ochner.
\newblock {\em J. Comput. Appl. Math.}, 68(1-2):103--114, 1996.

\bibitem{GH6}
F.~Alberto Gr{\"u}nbaum and Luc Haine.
\newblock Some functions that generalize the {A}skey-{W}ilson polynomials.
\newblock {\em Comm. Math. Phys.}, 184(1):173--202, 1997.

\bibitem{GH1}
F.~Alberto Gr{\"u}nbaum and Luc Haine.
\newblock On a $q$-analogue of the string equation and a generalization of the
  classical orthogonal polynomials.
\newblock In {\em Algebraic methods and $q$-special functions (Montr\'eal, QC,
  1996)}, pages 171--181. Amer. Math. Soc., Providence, RI, 1999.

\bibitem{GH3}
F.~Alberto Gr{\"u}nbaum and Luc Haine.
\newblock The {W}ilson bispectral involution: some elementary examples.
\newblock In {\em Symmetries and integrability of difference equations
  (Canterbury, 1996)}, pages 353--369. Cambridge Univ. Press, Cambridge, 1999.

\bibitem{GH2}
F.~Alberto Gr{\"u}nbaum, Luc Haine, and Emil Horozov.
\newblock Some functions that generalize the {K}rall-{L}aguerre polynomials.
\newblock {\em J. Comput. Appl. Math.}, 106(2):271--297, 1999.

\bibitem{HobIto}
S.~Hobart and T.~Ito.
\newblock The structure of nonthin irreducible ${T}$-modules of endpoint 1:
  ladder bases and classical parameters.
\newblock {\em J. Algebraic Combin.}, 7(1):53--75, 1998.

\bibitem{TD00}
T.~Ito, K.~Tanabe, and P.~Terwilliger.
\newblock Some algebra related to ${P}$- and ${Q}$-polynomial association
  schemes.
\newblock In {\em Codes and Association Schemes (Piscataway, NJ, 1999)}. Amer.
  Math. Soc., Providence, RI, 2000.

\bibitem{KoeSwa}
R.~Koekoek and R.~Swarttouw.
\newblock {\em The Askey-scheme of hypergeometric orthogonal polyomials and its
  $q$-analog}, volume 98-17 of {\em Reports of the faculty of Technical
  Mathematics and Informatics}.
\newblock Delft, The Netherlands, 1998.

\bibitem{Koelink3}
H.~T. Koelink.
\newblock Askey-{W}ilson polynomials and the quantum ${\rm {s}{u}}(2)$ group:
  survey and applications.
\newblock {\em Acta Appl. Math.}, 44(3):295--352, 1996.

\bibitem{Koelink1}
H.~T. Koelink.
\newblock $q$-{K}rawtchouk polynomials as spherical functions on the {H}ecke
  algebra of type ${B}$.
\newblock {\em Trans. Amer. Math. Soc.}, 352:4789--4813, 2000.

\bibitem{Koelink2}
H.~T. Koelink and J.~Van Der~Jeugt.
\newblock Convolutions for orthogonal polynomials from {L}ie and quantum
  algebra representations.
\newblock {\em SIAM J. Math. Anal.}, 29(3):794--822 (electronic), 1998.

\bibitem{Koelink4}
H.~T. Koelink and J.~Van~der Jeugt.
\newblock Bilinear generating functions for orthogonal polynomials.
\newblock {\em Constr. Approx.}, 15(4):481--497, 1999.

\bibitem{koo3}
Tom~H. Koornwinder.
\newblock Askey-{W}ilson polynomials as zonal spherical functions on the ${\rm
  {s}{u}}(2)$ quantum group.
\newblock {\em SIAM J. Math. Anal.}, 24(3):795--813, 1993.

\bibitem{Leodual}
D.~A. Leonard.
\newblock Orthogonal polynomials, duality and association schemes.
\newblock {\em SIAM J. Math. Anal.}, 13(4):656--663, 1982.

\bibitem{Leopandq}
D.~A. Leonard.
\newblock Parameters of association schemes that are both ${P}$- and
  ${Q}$-polynomial.
\newblock {\em J. Combin. Theory Ser. A}, 36(3):355--363, 1984.

\bibitem{Ronan}
Mark Ronan.
\newblock {\em Lectures on buildings}.
\newblock Academic Press Inc., Boston, MA, 1989.

\bibitem{Hjal}
H.~Rosengren.
\newblock {\em Multivariable orthogonal polynomials as coupling coefficients
  for {L}ie and quantum algebra representations}.
\newblock Centre for {M}athematical {S}ciences, Lund University, Sweden, 1999.

\bibitem{Tan}
K.~Tanabe.
\newblock The irreducible modules of the {T}erwilliger algebras of {D}oob
  schemes.
\newblock {\em J. Algebraic Combin.}, 6(2):173--195, 1997.

\bibitem{Tercharpq}
P.~Terwilliger.
\newblock A characterization of ${P}$- and ${Q}$-polynomial association
  schemes.
\newblock {\em J. Combin. Theory Ser. A}, 45(1):8--26, 1987.

\bibitem{TersubI}
P.~Terwilliger.
\newblock The subconstituent algebra of an association scheme. {I}.
\newblock {\em J. Algebraic Combin.}, 1(4):363--388, 1992.

\bibitem{TersubII}
P.~Terwilliger.
\newblock The subconstituent algebra of an association scheme. {I}{I}.
\newblock {\em J. Algebraic Combin.}, 2(1):73--103, 1993.

\bibitem{TersubIII}
P.~Terwilliger.
\newblock The subconstituent algebra of an association scheme. {I}{I}{I}.
\newblock {\em J. Algebraic Combin.}, 2(2):177--210, 1993.

\bibitem{Ternew}
P.~Terwilliger.
\newblock A new inequality for distance-regular graphs.
\newblock {\em Discrete Math.}, 137(1-3):319--332, 1995.

\bibitem{Terint}
P.~Terwilliger.
\newblock An introduction to {L}eonard pairs and {L}eonard systems.
\newblock {\em S\=urikaisekikenky\=usho K\=oky\=uroku}, (1109):67--79, 1999.
\newblock Algebraic combinatorics (Kyoto, 1999).

\bibitem{qSerre}
P.~Terwilliger.
\newblock Two relations that generalize the $q$-{S}erre relations and the
  {D}olan {G}rady relations.
\newblock In {\em Proceedings of Nagoya 1999 Workshop on Physics and
  Combinatorics (Nagoya, Japan, 1999)}. World Scientific Publishing Co., Inc,
  River Edge, NJ, Providence, RI, 2000.

\bibitem{LS99}
P.~Terwilliger.
\newblock Two linear transformations each tridiagonal with respect to an
  eigenbasis of the other.
\newblock {\em Linear Algebra Appl.}, 330(1--3):149--203, 2001.

\bibitem{Zhidd}
A.~S. Zhedanov.
\newblock ``{H}idden symmetry'' of {A}skey-{W}ilson polynomials.
\newblock {\em Teoret. Mat. Fiz.}, 89(2):190--204, 1991.

\bibitem{ZheCart}
A.~S. Zhedanov.
\newblock Quantum ${\rm {s}{u}}\sb q(2)$ algebra: ``{C}artesian'' version and
  overlaps.
\newblock {\em Modern Phys. Lett. A}, 7(18):1589--1593, 1992.

\bibitem{Zhidden}
A.~S. Zhedanov.
\newblock Hidden symmetry algebra and overlap coefficients for two ring-shaped
  potentials.
\newblock {\em J. Phys. A}, 26(18):4633--4641, 1993.

\end{thebibliography}


\medskip
\noindent Paul Terwilliger, Department of Mathematics, University of
Wisconsin, 480 Lincoln Drive, Madison, Wisconsin, 53706, USA \hfil\break
email: terwilli@math.wisc.edu \hfil\break
\end{document}